%% file: 20241201.tex
\documentclass[draftclsnofoot,onecolumn,journal,12pt]{IEEEtran}

\usepackage{amssymb,graphicx}
\usepackage{latexsym}
\usepackage{mathtools}
\usepackage{cite}
\usepackage{color}
\usepackage{subfigure}
\usepackage{bm,bbm}
\usepackage{tabularx,booktabs}
\usepackage{ltablex}
\usepackage{mathrsfs,dsfont}
\usepackage{enumitem}
\usepackage{cases}
\usepackage{yhmath}

\newcolumntype{M}[1]{>{\centering\arraybackslash}m{#1}}

\usepackage[linesnumbered,lined,boxed]{algorithm2e}
\usepackage{setspace}
\SetAlCapSkip{1em}
\SetKwComment{Comment}{$\triangleright$\ }{}
\RequirePackage[OT1]{fontenc}
\RequirePackage{amsthm,amsmath}
\RequirePackage[numbers]{natbib}
\RequirePackage[colorlinks,citecolor=blue,urlcolor=blue]{hyperref}
\newcommand*{\fullref}[2]{\hyperref[{#1}]{\ref*{#1} {\it(#2)}}}
\newcommand*{\partialref}[2]{\hyperref[{#1}]{#2}}

\renewcommand{\proof}[1]{\emph{#1}}
\renewcommand{\endproof}{{\hfill $\square$}}

\newtheorem{theorem}{\bf Theorem}
\newtheorem{lemma}{\bf Lemma}
\newtheorem{proposition}{\bf Proposition}
\newtheorem{corollary}{\bf Corollary}

\newcommand{\bA}{\bm{A}}
\newcommand{\bB}{\bm{B}}
\newcommand{\bF}{\mathbf F}
\newcommand{\bR}{\mathbf R}
\newcommand{\bS}{\bm{S}}

\newcommand{\bh}{\bm{h}}
\newcommand{\bs}{\bm{s}}
\newcommand{\balpha}{\bm{\alpha}}

\newcommand{\sPhi}{\Phi^{\text{stat}}}
\newcommand{\db}[1]{\Bar{\Bar{#1}}}
\newcommand{\dbphi}{\db \phi}
\newcommand{\dbpsi}{\db\psi}
\newcommand{\dbvtheta}{\db\vartheta}
\newcommand{\dba}{\db a}
\newcommand{\dbsPhiz}{\db{\Phi}^{\text{stat}}_{0}}

\newcommand{\bbR}{\mathbb R}
\newcommand{\bbP}{\mathbb P}
\newcommand{\ext}{\operatorname{ext}}

\newcommand{\bbE}{\mathbb E}

\newcommand{\oalpha}{\alpha}
\newcommand{\obmalpha}{\bm{\alpha}}
\newcommand{\oPhiLocal}{\Phi^{\text{LOC}}_0}

\newcommand{\sS}{\mathsf S}
\newcommand{\sA}{\mathsf A}
\newcommand{\sH}{\mathsf H}

\newcommand{\hatnu}{\hat{\nu}}
\newcommand{\hatnuVec}{\hat{\bm{\nu}}}
\newcommand{\Ia}{\mathfrak{I}}
\newcommand{\La}{\wideparen{a}}
\newcommand{\LaVec}{\wideparen{\bm{a}}}
\newcommand{\Lipschitza}{a}
\newcommand{\LenM}{\wideparen{M}}
\newcommand{\rank}{\mathfrak{r}}
\newcommand{\lPhi}{\Phi}  
\newcommand{\dblPhi}{\db{\Phi}}
\newcommand{\dblPhiz}{\db{\Phi}^0}

\newcommand{\actionphi}{\phi}
\newcommand{\actionpsi}{\psi}
\newcommand{\actionvarphi}{\varphi}
\newcommand{\actiondbphi}{\dbphi}
\newcommand{\actionbarphi}{\bar{\phi}}

\newcommand{\functionQ}{Q}

\newcommand{\indicator}{1}

\newcommand{\scrd}{\mathfrak{d}}

\newcommand{\scri}{\mathfrak{i}}

\newcommand{\scrz}{\mathfrak{z}}

\usepackage[deletedmarkup=xout]{changes}

\graphicspath{{./figs/}}

\newcounter{condition}
\newenvironment{condition}[1]{\refstepcounter{condition}
\par\medskip
\noindent \textbf{#1:~} \it}{\medskip \newline}

\newcounter{externalTheorem}
\newenvironment{externalTheorem}[1]{\refstepcounter{externalTheorem}
\par\medskip
\noindent \textbf{#1:~} \it}{\medskip \newline}

\newcommand{\norm}[1]{\left\lVert #1 \right\rVert}

\begin{document}
\title{Multi-Action Restless Bandits with Weakly Coupled Constraints: Simultaneous Learning and Controlling}
\author{Jing~Fu,
Bill~Moran~and~Jos\'e Ni\~no-Mora
\thanks{
Jing Fu is with Department of Electrical and Electronic Engineering, School of Engineering, STEM College, RMIT University, Australia (e-mail: jing.fu@rmit.edu.au).}
\thanks{Bill Moran is with Department of Electrical and Electronic Engineering, the University of Melbourne, VIC 3010, Australia (e-mail:wmoran@unimelb.edu.au).}
\thanks{Jos\'e Ni\~no-Mora is with Department of Statistics, Carlos III University of Madrid, Spain (E-mail: jose.nino@uc3m.es)}
}

\maketitle

\begin{abstract}
We study a system with finitely many groups of multi-action bandit processes, each of which is a Markov decision process (MDP) with finite state and action spaces and potentially different transition matrices when taking different actions. The bandit processes of the same group share the same state and action spaces and, given the same action that is taken, the same transition matrix. All the bandit processes across various groups are subject to multiple weakly coupled constraints over their state and action variables. Unlike the past studies that focused on the offline case, we consider the online case without assuming full knowledge of transition matrices and reward functions a priori and propose an effective scheme that enables simultaneous learning and control. We prove the convergence of the relevant processes in both the timeline and the number of the bandit processes, referred to as the convergence in the time and the magnitude dimensions. Moreover, we prove that the relevant processes converge exponentially fast in the magnitude dimension, leading to exponentially diminishing performance deviation between the proposed online algorithms and offline optimality. 
\end{abstract}


\section{Introduction}\label{sec:introduction}

\subsection{Overview and Motivation}

Restless multi-armed bandits (RMABs) provide a computationally feasible approach to the solution of many problems that can be couched as Markov Decision Processes (MDPs); see, for instance, \cite{krishnamurthy2007structured,wang2020whittle,villar2016indexability,fu2016asymptotic,fu2020energy,avrachenkov2016whittle}, and the recent survey \cite{ninomora2023survey}. As can be seen,  they have been  applied in a wide range of application scenarios. Originally  proposed by \cite{whittle1988restless}, structurally they comprise multiple parallel individual and otherwise independent \emph{bandit} MDP processes with binary actions \cite{gittins2011multiarmed}, which are linked by linear resource constraints.

Typically, the complexity resides in the multiplicity of such processes.
These overall RMABs are themselves  MDP with state space exponentially increasing in the number of  bandit processes, and in that regard are PSPACE-hard  \cite{papadimitriou1999complexity}.
Where the number of such bandits is  large, conventional solution methods for  MDP, such as value or policy iteration, cannot be applied due to the \emph{curse of dimensionality}.   Whittle's approach \cite{whittle1988restless} provides  a scalable (though not optimal) policy that assigns to  each state of each of the bandit processes a  real-valued \emph{Whittle index} by solving a relaxed version of the original problem. At each decision epoch, the Whittle policy gives higher priority to bandits with larger indices.
Such a policy is  subsequently referred to as  the \emph{Whittle index policy}.
Whittle index computations are linear in  the number of bandit processes; such  algorithms are discussed in \cite{nino2006restless,nino2007dynamic}.
Whittle \cite{whittle1988restless} conjectured that the Whittle index policy approaches optimality as the scale of the problem, measured by the number of the bandit processes, tends to infinity:  \emph{ asymptotic optimality}. This conjecture was proved by Weber and Weiss,  \cite{weber1990index},  under a non-trivial condition related to the existence of a global attractor of its associated averaging process.
The  global attractor condition is crucial, although it may exist  naturally in special cases.  We refer the reader to  discussions in \cite{weber1990index,verloop2016asymptotically,fu2018restless,fu2020energy}.

The conventional RMAB formulation is limited by the assumptions on binary actions and a single linear constraint in the  action variables.
\emph{Multi-action} bandit processes with more complex, and more  practical, constraints on both the action and state variables have been considered in, e.g.,  \cite{nino2008index,nino2022multigear,fu2018restless,brown2023fluid,fu2024patrolling,gast2024reoptimization,fu2025restless}.
In this paper, a  multi-action bandit process (MAB) is, simply,  an MDP with finite state and action spaces and bounded reward rates. 
We continue to call it  a bandit process in the tradition of RMAB community and to distinguish it from more general MDPs.
In \cite{brown2023fluid,fu2024patrolling,gast2024reoptimization}, these  MABs  are coupled through multiple constraints that fall in the scope of \emph{weakly coupled} MDPs (see \cite{adelman2008relaxations}). They have been  analyzed through  linear-programming approximation techniques with achieved asymptotic optimality in the finite-time-horizon case.
In \cite{fu2024patrolling}, for a special case of the weakly coupled constraints, the performance deviation was proved to achieve $O(e^{-N})$, where $N$ is the number of  bandit processes.
In \cite{brown2023fluid,gast2024reoptimization}, for general weakly-coupled bandit processes, under certain \emph{non-degeneracy} assumptions for a relaxed version of the problem,
the performance deviation between optimality and the designed scalable policies were also proved to achieve $O(e^{-N})$.
As demonstrated in \cite[Figure 2]{gast2024reoptimization}, the non-degeneracy assumptions are non-trivial and do not hold in general.

A key issue in the applicability of RMABs (and, more widely of MDPs) is lack of knowledge of the parameters; specifically of the transition kernels or the reward/cost functions.
In real-world applications  the transition kernels of the underlying stochastic process are usually unknown and, in situations where the MDPs are solvable, effective reinforcement learning techniques are used to explore necessary statistics and estimate these parameters,  while implementing the  control policies/algorithms to provide an optimal or close to optimal solution.
Unfortunately, in our context, the  large state and action spaces  prevents (efficient) convergence of conventional reinforcement learning techniques.

In this paper, our aim is to combine reinforcement learning methods  (specifically, $Q$-learning) \cite{bertsekas2015parallel} with a more general formulation of, and solution to,  the RMAB problem,  to provide a methodology for solving real-world decision and control  problems  where the parameters of the system are initially  unknown.  We will provide theoretical results to demonstrate convergence of our methods as well as implementable  algorithms and performance analysis.

We first generalize the notion of RMAB. In this paper, the basic object consists of a finite number of different groups of multiple (multi-action) bandit processes. These groups are called \emph{gangs} and each gang comprises bandit processes that have identical state and action spaces and state-action-dependent transition probabilities; that is, bandit processes in the same gang are stochastically identical.
Bandit processes in different gangs  can be entirely different.

We consider the general case of weakly-coupled constraints and refer to such a problem as  Weakly Coupled Gangs (WCGs). 
A detailed definition and explanation on the relationship between WCGs and conventional RMAB are provided in Section~\ref{subsec:model}.
WCGs encompass conventional RMABs and are at least as difficult to solve.

This paper focuses on simultaneous learning and control for  general WCG problems.
We propose a stream of methodologies that enable simultaneous, coordinated learning and control, and provide quantified estimations on how the scale of the problem, measured by the number of bandit processes $N$, can positively contribute to the overall performance.
We refer to the scale of the problem; that is $N$, as the \emph{magnitude dimension },  in contradistinction to the time dimension (timeline), which provides the underlying dimension for  the convergence of the learning process and, due to the features of the WCG problem, also leads to asymptotic optimality of the employed control algorithms/policies.
We refer to  algorithms/policies that do not presume  knowledge of the  transition kernels and/or  reward functions \emph{a priori}  as  \emph{online} algorithms; where the transition kernels and reward functions  are all assumed known we say that the algorithms are \emph{offline}.

We prove that the performance deviation between the online algorithms proposed in this paper and the offline optimality diminishes in $O(e^{-N})$ in both the infinite  and finite time horizon cases. 
Unlike the past work analyzing only the learning accuracy as $T\rightarrow \infty$, our results theoretically demonstrate that, even if $T$ is small,  as long as the WCG problem is realistically large, the $O(e^{-N})$ convergence of the underlying process(es) with coordinated learning and control ensures near-optimality of the proposed online algorithms.

\subsection{Results}
\label{sec:Results}

Our main contributions consist of  three parts for the general \emph{online} WCG problem; specifically,
\begin{enumerate}[label=(\arabic*)]
\item A scheme to enable  simultaneous control and learning with theoretically guaranteed convergence. \label{contribution:1}
\item Theoretical quantification of  the effect of the magnitude  dimension on  the convergence of the stochastic process(es). \label{contribution:2}
\item A class of policies with theoretically bounded performance degradation,    decreasing  exponentially in the magnitude  dimension.\label{contribution:3}
\end{enumerate}

\subsubsection{A scheme for simultaneous learning and control}
Contribution~\ref{contribution:1} is discussed in detail
in Section~\ref{subsec:q-index}.
We propose a scheme driven by a deployed \emph{primary} policy and $K$ \emph{artificial secondary} processes, in parallel with the primary  process. The primary process represents the unknown underlying  real-world process with a given policy,  whereas the $K$ secondary processes are  used to learn the  $Q$ factors and are guided by user-defined \emph{secondary} policies and reward functions.
The primary  process is the underlying stochastic process of the WCG  system, referred to as the  \emph{WCG process}.

We prove that the $K$ learning processes converge to the real values of the corresponding $Q$ factors (for different, user-defined secondary policies and reward functions) in the time dimension.
The learned $Q$ factors, even in the early stages,   can be used to construct $Q$-factor-based index policies, corresponding to  the Whittle index policy for the  RMAB case~\cite{nino2007dynamic,nino2022multigear}. 
We provide a detailed decision rule  for constructing an  online version of the Whittle index policy in Section~\ref{subsec:q-factor-based}.
This  online policy converges to the Whittle index policy since the estimated $Q$ factors of the $K$ learning processes converge to those of the primary process, and these  have been proved to be asymptotically optimal in a range of scenarios~\cite{weber1990index,fu2016asymptotic,fu2020energy,fu2018restless,fu2025restless}.

\subsubsection{Fast convergence in the magnitude dimension}

Contribution~\ref{contribution:2} is described in detail in Section~\ref{subsec:asym_regime}.
Under a mild assumption on  the probability distributions of the reward rates, we prove that the deviation between the estimated $Q$ factors of the $K$ learning processes and the real ones decreases as  $O(e^{-N})$, where $N$ is the magnitude  dimension.

The assumption only requests the reward rates should not vary too much from their true means so that the convergence of the learning processes is reasonably fast.
For example, if the reward rates are all normally distributed with bounded variance, then the assumption is satisfied, and the learning deviation diminishes in $O(e^{-N})$.
In this context, based on the convergence of the learning processes proved in the time dimension,   convergence is further improved  by  the large problem size (large $N$).
We believe that this  is the first attempt to  quantify the relationship between the convergence of the learning processes and the size of the WCG problem.

Observe that this  fast convergence happens in the magnitude dimension, and require relatively few time steps  for initial exploration, so that  it does not take a long time to approximate the $Q$ factors, thereby deleteriously affecting the policy.
Additionally large problem size also improves convergence between the proposed control policy, such as the online version of the Whittle index policy, and optimality. 

\subsubsection{Asymptotically optimal online policies}
Contribution~\ref{contribution:3}:  We describe  online policies for the WCG problem for the  infinite and finite time horizons cases in n Sections~\ref{sec:long-run} and \ref{sec:LP-algo}.
Fast convergence of the learned $Q$ factors, as well as the learned transition kernels and reward functions, are used to ensure that the proposed online policies are asymptotically optimal and that the performance degradations diminish in $O(e^{-N})$.
For the infinite time horizon case, we implement  $Q$-factor-based online policies, that generalize the online version of the Whittle index policy. For the  finite time horizon case,  we proposed a class of policies, for both the  online and offline cases, that achieve $O(e^{-N})$ performance degradation.
Unlike previous  work in \cite{brown2023fluid,gast2024reoptimization}, our proof for the $O(e^{-N})$ decrease in  performance degradation in the offline case does not rely on the non-degeneracy assumption.
A detailed literature survey is provided in Section~\ref{subsec:rWork}.


The remainder of the paper is organized as follows.
In Section~\ref{sec:model}, we give a detailed definition of the WCG problem. 
In Section~\ref{sec:multiple-dimensions}, we present our main results: a scheme with simultaneous learning and control for the WCG system with proved convergence in timeline, and theorems assure fast convergence of the learning and control processes in the magnitude dimension. 
Based on the fast convergence in the magnitude dimension,
in Sections~\ref{sec:long-run} and \ref{sec:LP-algo}, we propose online policies/algorithms for both infinite and finite time horizon cases, respectively, that achieve exponentially diminishing performance degradation for large-scale problems. 
In Section~\ref{sec:conclusions}, we present conclusions.

\subsection{Relations to the Literature} \label{subsec:rWork}

Previous work  has discussed the trade-offs between estimation of  unknown parameters and implementation of good policies/algorithms for  RMABs through the incorporation of  classical learning techniques with the  Whittle index policy and/or linear-programming-based approximations.
For instance, several work focused on approximating Whittle indices through UCB-based learning \cite{wang2023optimistic}, Bayes method \cite{jiang2023online}, two-time-scale learning \cite{avrachenkov2022whittle}, and Thompson sampling \cite{akbarzadeh2023learning}.
An optimistic-linear-programming-based learning technique that periodically explores necessary statistics and exploits the estimated Whittle indices was reported in \cite{wang2023optimistic}. .
Under a mild condition, \cite{wang2023optimistic} considered $N$ stochastically identical bandit processes and proved the convergence of the cumulative discounted total cost of the $N$ bandit processes for the relaxed RMAB problem as $T\rightarrow \infty$, achieving a regret  of $O(SN\sqrt{T\log T})$, where $S$ is the size of the state space for each bandit process and $T$ is the time horizon. A method periodically switching  the system between  control  and  Whittle index estimation using Thompson sampling is described in  \cite{akbarzadeh2023learning} . The performance regret between that proposed policy and the Whittle index policy is shown to be  $O(NS\sqrt{T \log T})$.

Bayesian estimation is used in \cite{jiang2023online}  for the online  case where  learning and control are done  simultaneously.
This method  exhibited relatively large computational complexity and, with some assumptions on the form of the transition kernels, achieved  a regret of $O((N^2+S^2)NS\sqrt{T\log (NT)})$ for the original RMAB problem.
A two-time-scale $Q$-learning technique for estimating Whittle indices is given in \cite{avrachenkov2022whittle}.
It combined learning the Lagrange-multiplier-based Q factors and the associated indices and requests the indices to be learned in a significantly slower time scale that that of the Q factors, and proved the convergence of the estimated Whittle indices to the real ones.


\section{The Problem}\label{sec:model}

The sets of positive and non-negative integers are denoted by $\mathbb{N}_{+}$ and $\mathbb{N}_{0}$ , respectively, and, for any $N\in\mathbb{N}_{+}$,  $[N]$ represents the set $\{1,2,\ldots,N\}$ with $[0]=\emptyset$.
We use $[N]_0$ to represent the set $\{0\}\cup[N]$.
Similarly,  $\mathbb{R}$, $\mathbb{R}_{+}$ and $\mathbb{R}_{0}$ denote the  sets of all, positive and non-negative reals, respectively.
Define $\db{\bR}$ as the set of all real-valued random variables.
For a finite set $\bS$, let $\bF(\bS)$ and $\db{\bF}(\bS)$ represent the sets of all real valued functions $\bS \rightarrow \mathbb{R}$ and all functions $\bS\rightarrow \db{\bR}$, respectively.
For any $f\in\bF(\bS)$, we define $\lVert f \rVert \coloneqq \max_{s\in\bS}\lvert f(s) \rvert$.

\subsection{System Model}\label{subsec:model}
We consider a system consisting of $I$ gangs  of restless multi-action bandit (RB) processes, where the $i$th gang  ($i\in [I]$) has $N_i$ stochastically identical   RB processes; that is,  with identical state and action spaces and transition kernels.
An RB process of class $i\in[I]$ is a discrete-time Markov decision process (MDP) with \emph{finite} state and action spaces $\bS_i$ and $\bA_i$, respectively.
Let $s_{i,n}(t)\in \bs_{i}$ and $\actionphi_{i,n}(t)\in \bA_{i}$ ($i\in[I],n\in[N_i]$) represent the state and action variables, respectively, of the $n$th RB process in gang    $i$ at time $t\in\mathbb{N}_0$. We write $\sS\coloneqq \prod_{i\in [I]}\bS_{i}^{{N_{i}}} $ and $\sA\coloneqq \prod_{i\in I}\bA_{i}^{{N_{i}}} $ for the full state and action space, respectively. These are given labels: $\sS_{t}$ and $\sA_{t}$  ($t\in \mathbb{N_{0}}$), so as to indicate time,  though,  as sets, they are identical; so that
$\bs(t)\coloneqq (s_{i,n}(t): i\in[I], n\in[N_i])\in \sS_{t}$ and $\actionphi(t)\coloneqq (\actionphi_{i,n}(t): i\in[I], n\in[N_i])\in \sA_{t}$.
We will write the history variable as
\[\bh(t)\coloneqq (\bs(0),\actionphi(0),\bs(1),\actionphi(1),\ldots,\bs(t))\in \sS_{0}\times \sA_{0}\times\sS_{1}\times\sA_{1}\times\cdots\times \sA_{t-1}\times \sS_{t}\eqqcolon \sH_{t}.\]
At each time $t\in \mathbb{N}_0$, the system controller chooses  the value of the action vector $\actionphi(t)$ as a function of the preceding history  $\bh(t)\in \sH_{t}$.

For a system of this kind, a \emph{policy} $\bm \phi$ is a sequence of maps $\actionphi_{t}:\sH_{t}\to \sA_{t}$ ($t\in \bm{N}_{0}$). The set of all policies is denoted by $\Phi$.  Evolution of the system under such a policy requires that the action $\actionphi(t)$ at time $t$ is
\begin{equation*}
  \actionphi_{t}=\actionphi_{t}(\bh(t))\in \sA_{t}.
\end{equation*}
Usually, but not always, we will focus on \emph{stationary policies}; that is, ones where
\[\actionphi_{t}=\actionphi:\sS_{t}=\sS\to \sA_{t}=\sA\qquad (t\in {\bm N}_{0})\] is just a function on  the current state and is unchanged in time. The set of all stationary policies is denoted by $\Phi^{\text{stat}}$.

Choice of the action $\actionphi(t)=(\actionphi_{i,n})$ provides the transition matrix  $\mathcal{P}_i(a)=\bigl[p_i(s,a,s')\bigr]_{|\bS_i|\times|\bS_i|}\in[0,1]^{|\bS_i|\times|\bS_{i}|}$ ($a\in\bA_i$),
so that  the RB process $\bigl\{s_{i,n}(t),t\in\mathbb{N}_0\bigr\}$ evolves according to  the transition probability from $s_{i,n}(t)$ to $s_{i,n}(t+1)$ is $p_i\bigl(s_{i,n}(t),\actionphi_{i,n}(t),s_{i,n}(t+1)\bigr)$. This state transition generates a real-valued, non-negative and bounded random reward $R_{i,n}(s_{i,n}(t),\phi_{i,n}(t))$,  with expectation $r_i\bigl(s_{i,n}(t),\actionphi_{i,n}(t)\bigr)\in\mathbb{R}$, where  $R_{i,n}\in\db{\bF}(\bS_i\times\bA_i)$ and $r_i\in\bF(\bS_i\times\bA_i)$.
Such $r_i$ is called the \emph{reward function.}
The boundedness assumption on the random reward  is uniform:  $R_{i,n}(t) \leq R_{\max} < \infty$ for all $i\in[I]$, $n\in[N_i]$, and $t\in\mathbb{N}_0$. Typically, policies are chosen to maximize long-term reward, as will be described below.

To highlight the dependence of the various objects studied on the choice of  policy $\phi$,  we write $s_{i,n}(t)$ and $R_{i,n}(t)$,  as $s^{\phi}_{i,n}(t)$ and $R^{\phi}_{i,n}(t)$, respectively. We will adopt a similar notation, as required, for other derived entities in relation to the system.

The original  RMAB formulation~\cite{whittle1988restless} imposes a simple form of constraint:
\begin{equation}\label{eqn:constraint:rmab}
\sum_{i\in[I],n\in[N_i]}\actionphi_{i,n}(t) = M,~\forall t=0,1,\ldots,T,
\end{equation}
where $T\in\mathbb{N}_+\cup\{+\infty\}$ is the time horizon, $M\in[N]$, and $\bA_i$ is specified to be $\{0,1\}$ for all $i\in[I]$.
It couples the RB processes $\bigl\{s_{i,n}(t),t\in\mathbb{N}_0\bigr\}$, enforcing dependencies among $s_{i,n}(t)$ and $\actionphi_{i,n}(t)$, for all $i\in[I]$ and $n\in[N_i]$.
For $T\in\mathbb{N}_+$ and $\beta \in (0,1]$, let 
\begin{equation}\label{eqn:obj_func}
\Gamma^\beta(T,\bs_0) \coloneqq \mathbb{E}\Bigl[\sum_{t\in[T]}\sum_{i\in[I],n\in[N_i]}\beta^t R_{i,n}(t)\Bigr|~\bs(0)=\bs_0\Bigr],
\end{equation}
where the initial state $\bs_0$ can be a random variable and, of course, the value is dependent on the action choices.
Define $\Gamma(T,\bs_0) \coloneqq \Gamma^1(T,\bs_0)$.
In \cite{whittle1988restless,weber1990index,nino2007dynamic,verloop2016asymptotically,ninomora2020verification,fu2018restless,brown2020index}, RMAB techniques  deal with maximization of  the long-run average expected reward, $\lim_{T\rightarrow \infty}\Gamma(T,\bs_0)/T$,  the long-run discounted expected cumulative reward $\lim_{T\rightarrow\infty} \Gamma^{\beta}(T,\bs_0)$ with $\beta<1$, and/or the expected (discounted) cumulative reward of the entire system $\Gamma^{\beta}(T,\bs_0)$, where the state and action variables are fully observed and the transition matrices and the reward functions are known a priori  (the offline case).

In \cite{weber1990index,verloop2016asymptotically,fu2018restless,brown2020index,fu2024patrolling}, scalable  algorithms were proposed for large-scale problems and were proved (under certain conditions) to approach optimality as the problem sizes, usually measured by the number of  RB processes, tend to infinity. In particular, the deviations between the proposed algorithms and optimality decrease exponentially or polynomially as the problem size increases.

In this paper, we consider a set of constraints that generalize \eqref{eqn:constraint:rmab}, 
\begin{equation}\label{eqn:constraint:linear}
    \sum_{i\in[I]}\sum_{n\in[N_i]}f_{i,\ell}\bigl(s_{i,n}(t),\actionphi_{i,n}(t)\bigr)=0,~\forall \ell\in[L], t\in[T]_0,
\end{equation}
where $L\in\mathbb{N}_0$, $T\in\mathbb{N}_+\cup\{\infty\}$,  and $f_{i,\ell}:\bS_i\times\bA_i \rightarrow \mathbb{R}$ is a bounded function.
We consider equality in \eqref{eqn:constraint:linear} for the sake of simplicity. It is straightforward  to modify the results presented here to the case where we replace some of the  equalities  \eqref{eqn:constraint:linear} by  inequalities; that is, having a mix of equalities and inequalities. Inequalities involve a positivity constraint on the Lagrange multipliers in a standard way.

We refer to such a process $\bigl\{\bs^{\phi}(t), t\in[T]_0\bigr\}$, consisting of the $I$ classes of RB processes coupled through \eqref{eqn:constraint:linear}, as a
\emph{Weakly Coupled Gang (WCG)} process.
A WCG is a Markov decision process with state space $\sS=\prod_{i\in[I]}\bS_i^{N_i}$ and, of course,  is at least as difficult as a general RMAB problem.

Unlike most of the past studies, we generalize the conventional constraints \eqref{eqn:constraint:rmab} to \eqref{eqn:constraint:linear} and do not assume any knowledge of the transition matrices $\mathcal{P}_i(a)$ or the reward functions $r_i$ for $i\in[I]$ and $a\in\bA_i$ \emph{a priori}.
Here, we develop methodologies, for both long-term and short-term objectives, that enable simultaneous learning and control of WCGs, for which the newly proposed algorithms are proved to converge to optimality, with exponential decrease in deviation.

\section{Convergence in Multiple Dimensions}\label{sec:multiple-dimensions}

\subsection{The Offline Problem}
\label{subsec:offline-problem}

Here we consider the \emph{offline} case: we assume that all parameter values are known to the decision maker.  WCGs are somewhat more general than RMAB generalizations discussed in the literature, both in terms of the structure of the state space and in the set of constraints,  so the results here are new. However, the key feature of this paper is that we can achieve close to optimal policies with online learning. This section serves as a stepping stone in that direction. We will talk about online learning in the following section.

We focus on a subset $\sPhi_{0}\subset \sPhi$ of stationary policies $\phi$, for which the following condition holds.
\begin{condition}{Ergodic Condition}\label{condition:ergodic}
There is a state $\bs_{0}=(s^{0}_{i})\in\bS$ and $T<\infty$ such that
\begin{equation}\label{eqn:proper-condition}
    \min_{t=0,1,\ldots,T}\mathbb{P}\Bigl\{s^{\phi}_{i,n}(t) \neq s^{0}_{i}\Bigl| s^{\phi}_{i,n}(0) = s_{i}\Bigr\} < 1,~\forall i\in[I],s_{i}\in\bS_i,
\end{equation}
and the process $\{s^{\phi}_{i,n}(t),t\in\mathbb{N}_0\}$ is aperiodic. 
\end{condition}
Such a state $\bs_{0}$ is called \emph{ergodic} and policies in $\sPhi_{0}$ are called \emph{ergodic policies}.

We are interested in the case with $\sPhi_{0}\neq \emptyset$.
For an ergodic policy, $\phi$, the  underlying Markov chain for $\{s^{\phi}_{i,n}(t),t\in\mathbb{N}_0\}$ is aperiodic and includes at most one recurrent class. If $\bS^{REC}_i$ is  the set of all recurrent states in $\bS_i$,
then $s_0\in\bS^{REC}_i$ and all states $s\in\bS^{REC}_i$ communicate with each other;  the corresponding Markov chain is irreducible.
Note that such an $\bs_0=\bs_{0}^{\phi}$ is policy-dependent.
From \cite[Theorem 1.8.3, Theorem 1.10.2]{norris1998markov},
there is a unique steady state distribution $\bm{\pi}^{\phi}_i$ on $\bS_i$ with support equal to $\bS^{REC}_{i}$. This also leads, via the map $s\to (s, \actionphi(s))$ to a ``steady state'' distribution on the set $\bS_i\times \bA_i$ that, with abuse of notation, we also denote by $\bm{\pi}^{\phi}_{i}(s,a)$. In this context, the Ergodic Theorem states that
\begin{equation}
  \label{eq:4}
  \lim_{T\to\infty } \frac1{T}\sum_{t\in [T]} g(s^{\phi}(t), \actionphi(t)) = \sum_{s\in \bS_{i}\\a\in \bA_{i}} g(s, a) \pi^{\phi}_{i}(s,a)
\end{equation}
for any function $g:\bS_{i}\times \bA_{i}\to \bbR$ and for $i\in [I]$.

We are interested in finding  a policy $\phi\in \sPhi_{0}$ that, at least approximately, achieves  the long-run objective
\begin{equation}\label{eqn:obj:long-run-average}
\max_{\phi\in\sPhi_{0}} \limsup_{T\rightarrow \infty} \frac{1}{T}\Gamma^{\phi}(T,\bs_0),
\end{equation}
subject to \eqref{eqn:constraint:linear}.

We follow the general methodology  of \cite{whittle1988restless}: we randomize the policy $\phi$.  We will indicate  randomized policies by putting a double bar over the symbol, thus $\dbphi$ and use $\db{\sA}$ to indicate the space of randomized actions; that is, the space of random variables with values in $\db(\sA)$.

Given a state $\bs^{\db{\phi}}(t)=\bs$, $\actiondbphi(\bs)\coloneqq(\actiondbphi_{i,n}(\bs):i\in[I],n\in[N_i])$ is a random action in $\db{\sA}$, with distribution
\[\alpha^{\db{\phi}}_{{i,a}}(\bs)\coloneqq
 \bbP\Bigl[\actiondbphi_{i,n}(\bs)= a\Bigr], \qquad (n\in [N_{i}], i\in [I],a\in\bA_i).\]
Naturally,  the action probabilities satisfy $\sum_{a\in\bA_i}\alpha^{\dbphi}_{i,a}(\bs)=1$.
For $i\in[I]$, define $\balpha^{\dbphi}_i(\bs)\coloneqq (\alpha^{\dbphi}_{i,a}(\bs):a\in\bA_i)$, and a  simplex  $\Delta_{\bA_i}\coloneqq \bigl\{\balpha: \bA_{i}\to \bbR_{+}|\sum_{a\in\bA_i}\alpha_a = 1\bigr\}$.
The random action $\actiondbphi(\bs)$ and the action probability $\balpha^{\dbphi}_i(\bs)$ are functions $\sS\rightarrow \db{\sA}$ and $\sS \rightarrow \Delta_{\bA_i}$, respectively, where $\db{\sA}$ is the space of the random actions.
Define $\dbsPhiz$ as the set of all the policies $\dbphi$ determined by $\actiondbphi(\bs)$, or equivalently $\balpha^{\dbphi}(\bs)$, for all $\bs\in\sS$, such that the \partialref{condition:ergodic}{ergodic condition} is satisfied with  $\phi$ replaced by  $\dbphi$. Note that, because of \eqref{eq:4},
\begin{equation}
  \label{eq:4ddd}
  \lim_{T\to\infty } \frac1{T}\sum_{t\in [T]} \bbE\bigl[g(s^{\dbphi}(t), \actiondbphi(t))\bigr] = \sum_{s\in \bS_{i}\\a\in \bA_{i}} g(s, a) \pi^{\dbphi}_{i}(s,a),
\end{equation}
for $\dbphi\in\dbsPhiz$, where  $\pi^{\dbphi}_{i}(s,a)$ is the distribution averaged over the random policy.

For policies $\dbphi\in\dbsPhiz$, relax \eqref{eqn:constraint:linear} to
\begin{equation}\label{eqn:constraint:long-run:relax}
        \limsup_{T\rightarrow \infty}\frac{1}{T}\sum_{t\in[T]}\sum_{i\in[I]}\sum_{n\in[N_i]}\bbE\Bigl[f_{i,\ell}\bigl(s^{\dbphi}_{i,n}(t),\actiondbphi_{i,n}(t)\bigr)\Bigr]=0,~\forall \ell\in[L_0],
      \end{equation}
with given initial state $\bs^{\dbphi}(0) = \bs_0$,  where $\bs_0$ is a random variable and the expectation is taken over the randomized policy.

Consider a modified objective 
\begin{equation}\label{eqn:obj:long-run-average:relax}
    \max_{\dbphi\in\dbsPhiz} \limsup_{T\rightarrow \infty} \frac{1}{T}\Gamma^{\dbphi}(T,\bs_0).
\end{equation}
Following the convention of the RMAB community, we refer to the problem described in \eqref{eqn:obj:long-run-average:relax} and \eqref{eqn:constraint:long-run:relax} as the \emph{relaxed} version of the original WCG described in \eqref{eqn:obj:long-run-average} and \eqref{eqn:constraint:linear}.
In much  previous  work \cite{whittle1988restless,weber1990index,nino2007dynamic,fu2016asymptotic,fu2018restless,fu2020energy,ninomora2020verification,brown2020index}, such a relaxed problem is an intermediate formulation used to quantify marginal rewards for each of the states $s\in\bS_i$ ($i\in[I]$) that lead to the \emph{(Whittle) index policy} with proved asymptotic optimality in a range of cases (see \cite{weber1990index,fu2020energy,fu2018restless,fu2025restless}).

Taking into account the constraint \eqref{eqn:constraint:long-run:relax}, we perform a further relaxation with  Lagrange multipliers $\bm{\gamma}=(\gamma_{\ell}:\ell\in[L])\in\mathbb{R}^L$ for the $L$ constraints in \eqref{eqn:constraint:long-run:relax},
 to produce  the dual function
\begin{equation}\label{eqn:long-run:dual-func}
    D(\bm{\gamma}) = \max_{\phi\in\dbsPhiz}\sum_{i\in[I]}\sum_{n\in[N_i]}\limsup_{T\rightarrow \infty}\frac{1}{T} \sum_{t\in[T]}\mathbb{E}\biggl[ R^{\dbphi}_{i,n}(t)- \sum_{\ell\in[L]}\gamma_{\ell} f_{i,\ell}\bigl(s^{\dbphi}_{i,n}(t),\actiondbphi_{i,n}(t)\bigr)\biggr],
\end{equation}
Now, Whittle's trick allows the  maximization on the right hand side of \eqref{eqn:long-run:dual-func} to  be decomposed into $I$ sub-problems.
For $\dbphi\in\dbsPhiz$, $i\in[I]$, $a\in\bA_i$, and $s\in\bS_i$, let
  \begin{align}\label{eq:1}
    \oalpha^{\dbphi}_{i}(a,s) &\coloneqq \lim_{t\to \infty}\bbE\bigl[\alpha^{\dbphi}_{i,a}(\bs^{\dbphi}(t)) \bigl| s^{\dbphi}_{i,n}(t) = s\bigr],\qquad  ( n\in [N_i]) \\
    \intertext{which limit exists by ergodicity and aperiodicity, and let}
 \obmalpha^{\dbphi}_i &\coloneqq (\oalpha^{\dbphi}_{i}(a,s): a\in\bA_i,s\in\bS_i).
  \end{align}
For $i\in[I]$, the value of $\obmalpha^{\dbphi}_i$ is sufficient to determine the transition matrix of the underlying (stationary)  process $\{s^{\dbphi}_{i,n}(t), t\in[T]_0\}$ (for any $n\in[N_i]$).

For $i\in[I]$, recall the simplex $\Delta_{\bA_i}$ with extreme points the set of functions $\bS_i\rightarrow \Delta_{\bA_i}$ that taking values only in $\{0,1\}$.  We write $\ext(\Delta_{\bA_{i}})$ for the set of extreme points.
Evidently, for $\dbphi\in \dbsPhiz$, $\oalpha^{\dbphi}_{i}\in\Delta_{\bA_i}$.
We write $\oPhiLocal$ for the set of those policies taking values only in $\ext(\Delta_{\bA_{i}})$.
Such a policy $\dbphi\in\oPhiLocal$ is almost sure to take a specific  action $a\in\bA_i$; that is, all the action variables of $\dbphi\in\dbsPhiz$ reduce to the deterministic version. 
For $\phi\in\oPhiLocal$, with some abuse of notation, we denote such a deterministic action as $\actionphi_i(s^{\phi}_{i,n}(t))=a$ with a function $\actionphi_i: \bS_i\rightarrow \bA_i$.

We define the dual function for the $i$th subproblem ($i\in[I]$),  $n\in[N_i]$, $\bs^{\phi}(0)=\bs_0$, and $\bm{\gamma}\in \bbR^{L}$,
\begin{equation}\label{eqn:sub-problem}
    D_i(\bm{\gamma})\coloneqq \max_{\dbphi\in\dbsPhiz}\limsup_{T\rightarrow \infty}\frac{1}{T} \sum_{t\in[T]}\mathbb{E}\biggl[ R^{\dbphi}_{i,n}(t)- \sum_{\ell\in[L]}\gamma_{\ell} f_{i,\ell}\bigl((s^{\dbphi}_{i,n}(t),\actiondbphi_{i,n}(t)\bigr)\biggr]
\end{equation}
Now we have  the following proposition.
\begin{proposition}\label{prop:decomposition} Fix $\bm{\gamma}\in \bbR^{L}$. Then,
  \begin{enumerate}
    \item  for $i\in [I]$,
  \begin{equation}
    \label{eq:3}
   D_{i}(\bm{\gamma})     = \max_{\phi\in\oPhiLocal}\sum_{s\in\bS_i}\sum_{a\in\bA_i}\pi^{\phi}_i(s,a)\Bigl(r_i(s,a)-\sum_{\ell\in[L]}\gamma_{\ell}f_{i,\ell}(s,a)\Bigr),
  \end{equation} where  $\bm{\pi}^{\phi}_i = (\pi^{\phi}_i(s,a):s\in\bS_i,a\in\bA_i)$ is the steady state distribution of the process $\bigl\{s^{\phi}_{i,n}(t), t\in\mathbb{N}_0\bigr\}$, under the policy $\phi\in\oPhiLocal$; and
\item    \begin{equation}\label{eqn:decomposition}
    D(\bm{\gamma}) = \sum_{i\in[I]}N_i D_i(\bm{\gamma}).
\end{equation}

  \end{enumerate}
\end{proposition}
The proof of Proposition~\ref{prop:decomposition} is given  in Appendix~\ref{app:prop:decomposition}.

The optimal solutions to the $I$ sub-problems in general do not satisfy the original constraints \eqref{eqn:constraint:linear} and hence are not feasible for the original WCG.
For the RMAB problem, \cite{whittle1988restless} proposed the well-known \emph{Whittle index} policy, which is based on the Whittle indices assigned to each state-action (SA) pair $(s,a)\in\bS_i\times\bA_i$ ($i\in[I]$) obtained by solving the $I$ sub-problems.
More precisely, for $i\in[I]$ and SA pair $(s,a)\in\bS_i\times\bA_i$, let $Q^{\dbphi,\bm{\gamma}}_i(s,a)$ represent the $Q$-factor for the associated process $\bigl\{s^{\phi}_{i,n}(t),t\in\mathbb{N}_0\bigr\}$ for any $n\in[N_i]$ under policy $\dbphi\in\dbsPhiz$, satisfying the following Bellman equation
\begin{equation}\label{eqn:q-factor}
    Q^{\dbphi,\bm{\gamma}}_i(s,a) = r_i(s,a) - \sum_{\ell\in[L]}\gamma_{\ell}f_{i,\ell}(s,a) +\sum_{s'\in\bS_i\backslash\{s_0\}}p_i(s,a,s')\sum_{a'\in\bA_i}\alpha^{\dbphi}_{i}(a',s')Q^{\dbphi,\bm{\gamma}}_i(s',a'),
  \end{equation}
where $s_0$ is an ergodic state.   
Let $Q^{\bm{\gamma}}_i(s,a)$ represents the  $Q$-factor corresponding to the optimal long-run average reward, so that,  for all $(s,a)\in\bS_i\times\bA_i$,
\begin{equation}\label{eqn:bellman}
    D_i(\bm{\gamma}) + Q^{\bm{\gamma}}_i(s,a) = r_i(s,a) - \sum_{\ell\in[L]}\gamma_{\ell}f_{i,\ell}(s,a)+\sum_{s'\in\bS_i\backslash\{s_0\}}p_i(s,a,s')\max_{a'\in\bA_i}Q^{\bm{\gamma}}_i(s',a'),\qquad \bigl((s,a)\in\bS_i\times\bA_i\bigr),
\end{equation}
where $D_i(\bm{\gamma})$ is defined in \eqref{eqn:sub-problem}.
In the  original  RMAB problem with $\bA_i = \{0,1\}$ and $L=1$, the vector $\bm{\gamma}$ is just  a scalar $\gamma$, and the Whittle index of SA pair $(s,1)$ ($s\in\bS_i,i\in[I]$) is
\begin{equation}\label{eqn:whittle-index}
    \gamma_i(s)\coloneqq \min \bigl\{\gamma \in\mathbb{R}| Q^{\gamma}_i(s,0)=Q^{\gamma}_i(s,1)\bigr\}.
\end{equation}

The Whittle index policy greedily prioritizes the RB processes with the highest Whittle indices and, when the transition matrices and reward functions are known, \emph{a priori}, has been proved to approach optimality as $N_i\to \infty$ ($i\in[I]$),  proportionately, in a wide range of scenarios \cite{weber1990index,verloop2016asymptotically,fu2020energy,fu2018restless,fu2025restless}.
We will further discuss the $Q$-factor based index policies in Section~\ref{sec:long-run}, where we will also introduce more recent results for multi-action RMAB problems (generalized $|\bA_i|\geq 2$ from the classic binary-action RMAB).

\subsection{Q-Factor-Based Learning and Control: Convergence in Timeline}\label{subsec:q-index}

Here, we describe how  learning and control can be conducted simultaneously with learned $Q$-factors for the general MCG process.
The $Q$-factors can further be used for approximating $Q$-factor-based algorithms, such as the Whittle indices through the Marginal Productivity (MP) indices, presented and analyzed in \cite{nino2001restless,nino2002dynamic,nino2006restless,nino2007dynamic,nino2022multigear}, for the conventional (multi-action) RMAB.

For $\phi\in\oPhiLocal$ and $i\in[I]$, and $\bar{r}_i\in\bF(\bS_i\times\bA_i)$, we 
define an affine  operator $\mathcal{T}_i^{\phi}:\bF(\bS_i\times\bA_i)\to\bF(\bS_i\times\bA_i)$ by
\begin{equation}\label{eqn:define_H_phi}
\mathcal{T}^{\phi}_i Q(s,a) = \bar{r}_i(s,a) + \sum_{s'\in\bS_i\backslash\{s_0\}}p_i(s,a,s')Q(s',\actionphi_i(s')), \qquad \bigl(Q\in \bF(\bS_i\times \bA_i)\bigr). 
\end{equation}
 We define $\mathcal{T^{\phi}}$ by
$\mathcal{T}^{\phi}\bm{Q}\coloneqq (\mathcal{T}^{\phi}_i Q_i)_{i\in [I]}$, where $\bm{Q}=(Q_i)_{i\in [I]}$ and $Q_i\in \bF(\bS_i\times \bA_i)$.

For $\phi\in\Phi$, $t\in[T]$ and $i\in[I]$, let
\begin{equation}
\label{eq:subprocesses}
    \begin{aligned}
        \hat{\mathscr{N}}^{\phi}_{i,t}(s,a)&=\{n\in  [N_i]: (s^{\phi}_{i,n}(t),\actionphi_{i,n}(t))  = (s,a)\}\quad (s\in \bS_i,\, a\in \bA_i)\\
        \hat{\mathscr{N}}^{\phi}_{i,t}(s,a,s')&=\{n\in  [N_i]: (s^{\phi}_{i,n}(t),\actionphi_{i,n}(t))  = (s,a), \, s^{\phi}_{i,n}(t+1)  =s'\}\quad (s, s'\in \bS_i,\, a\in \bA_i)
    \end{aligned}
\end{equation}
so that $ \hat{\mathscr{N}}^{\phi}_{i,t}(s,a,s')\subset \hat{\mathscr{N}}^{\phi}_{i,t}(s,a)$, for all $s'\in \bS_i$. 

For convenience, we fix a sequence of  real numbers, $\eta_t>0$  (such as $\eta_t =1/(t+1)$), satisfying  
\begin{equation}
    \label{eqn:assumption:stepsize}
    \sum_{t=0}^\infty \eta_t =\infty,\quad    \sum_{t=0}^\infty \eta_t^2 <\infty. 
\end{equation}For the  WCG process $\bigl\{\bs^{\phi}(t),t\in [T]_0\bigr\}$,  we define  \emph{step sizes}, $\bm{\eta}_{t}\coloneqq(\eta_{i,t}(s,a):i\in[I],(s,a)\in\bS_i\times\bA_i)$ ($t\in\mathbb{N}_0$), by 
, 
Now we assume a map, to be specified later,  $t\mapsto\phi_t: [T]_0 \rightarrow \oPhiLocal$, and  define, recursive, two objects: the estimated $Q$-factor, $\hat{\functionQ}_{i,t}$, and the estimation deviation $w_{i,t}$.  This is done as follow:  given $\hat{\functionQ}_{i,t}$, and observed $\bs^{\phi}(t+1)$, we update
\begin{align}
 \label{eqn:q-iteration}
    \hat{\functionQ}_{i,t+1} (s,a) & \coloneqq  \begin{cases}  (1-\eta_{t})\hat{\functionQ}_{i,t}(s,a) + \eta_{t}\Bigl(\bigl(\mathcal{T}^{\phi_t}_i\hat{\functionQ}_{i,t}\bigr)(s,a)+w_{i,t}(s,a)\Bigr)&   \text{ if } \hat{\mathscr{N}}^{\phi}_{i,t}(s,a)\neq \emptyset,\\
    0&\text{ otherwise}
    \end{cases},\\\label{eqn:define_w}
    w_{i,t}(s,a) &\coloneqq \begin{cases}
        \frac{1}{\bigl|\hat{\mathscr{N}}^{\phi}_{i,t}(s,a)\bigr|}\Bigl(\sum\limits_{n\in\hat{\mathscr{N}}^{\phi}_{i,t}(s,a)}\bar{R}_{i,n}(s,a) +\sum\limits_{s'\in\bS_i\backslash\{s_0\}}\bigl|\hat{\mathscr{N}}^{\phi}_{i,t}(s,a,s')\bigr|\hat{\functionQ}_{i,t}(s',\actionphi_{t,i}(s'))\Bigr)\\ 
        \hfil - \bar{r}_i(s,a) - \sum_{s'\in\bS_i\backslash\{s_0\}}p_i(s,a,s')\hat{Q}(s',\actionphi_{t,i}(s')),
    & \text{ if } \hat{\mathscr{N}}^{\phi}_{i,t}(s,a)\neq \emptyset,\\
    0, & \text{ otherwise},
    \end{cases}
\end{align}
where the initial value $\hat{\functionQ}_{i,0}$ is prescribed to be any given function in $\bF(\bS_i\times\bA_i)$ for all $i\in[I]$.

Here the  policy map  $\phi_t:[T]_0\rightarrow\oPhiLocal$ is a  user-defined function of time $t$ and can be different from the policy $\phi$  in the WCG process $\bigl\{\bs^{\phi}(t),t\in[T]_0\bigr\}$.
The variable $\bar{R}_{i,n}\in\db{\bF}(\bS_i\times\bA_i)$ (same  distribution for any $n\in[N_i]$) is the instantaneous reward rate with $\bar{R}_{i,n}(s,a)$ bounded for all $(s,a)\in\bS_i\times \bA_i$, and $\bar{r}_i(s,a)=\mathbb{E}[\bar{R}_{i,n}(s,a)]$ is the expectation given $(s,a)\in\bS_i\times\bA_i$ for any $n\in[N_i]$.
Let $\hat{\bm{Q}}_t\coloneqq (\hat{\functionQ}_{i,t})_{ i\in [I]}$, and $\bm{w}_t\coloneqq (w_{i,t}(s,a):i\in[I],s\in\bS_i, a\in\bA_i)$.

We refer to the process $\bigl\{\hat{\bm{Q}}_t,t\in [T]_0\bigr\}$ as the \emph{learning process} accompanying $\big\{\bs^{\phi}(t), t\in[T]_0\bigr\}$ (for $\phi\in\lPhi$), and refer to the combined process $\big\{\bs^{\phi}(t),\hat{\bm{Q}}_t, t\in[T]_0\bigr\}$ as the \emph{WCG-Learning} process.
Note that the learning process  $\bigl\{\hat{\bm{Q}}_t,t\in [T]_0\bigr\}$ is an artificial process utilizing the real process $\big\{\bs^{\phi}(t), t\in[T]_0\bigr\}$ as an event generator, and is dependent on both $\phi\in\lPhi$ for the real process and $\phi_t\in\oPhiLocal$  ($t\in [T]_0$) for translating the collected data.
We refer to such a $\phi\in\lPhi$ as the \emph{primary policy} and $\phi_t\in\oPhiLocal$ as the \emph{secondary policy} at time $t$.
Recall that the secondary policy $\phi_t\in\oPhiLocal$ is user-defined and can be different from the primary policy $\phi\in\lPhi$, for which the action variable $\phi(t)$ is a function of the history $\bh_t$. 

\begin{lemma}\label{lemma:exist_Q_star}
For any $\phi\in\oPhiLocal$, the operator $\mathcal{T}^{\phi}$  has a unique fixed point $\hat{\bm{Q}}^{\phi} \coloneqq (\hat{Q}^{\phi}_i)_{i\in [I]}$ with $\hat{Q}^{\phi}_i \in \bF(\bS_i\times\bA_i)$: 
\begin{equation}\label{eqn:lemma:exist_Q_star}
    \mathcal{T}^{\phi}\hat{\bm{Q}}^{\phi} = \hat{\bm{Q}}^{\phi}.
\end{equation}
\end{lemma}
The proof of Lemma~\ref{lemma:exist_Q_star} is provided in Appendix~\ref{app:lemma:exist_Q_star}.

\begin{proposition}\label{prop:q-convergence}
For any primary policy $\psi\in\lPhi$, secondary policy map $t\mapsto \phi_t=\phi\in\oPhiLocal$ ( $t\in\mathbb{N}_0$), $\epsilon, \delta >0$, 
there exists $T<\infty$ such that, for all $t>T$, 
\begin{equation}\label{eqn:q-convergence}
    \mathbb{P}\Bigl\{\bigl\lVert \hat{\bm{Q}}_t - \hat{\bm{Q}}^*\bigr\rVert > \epsilon\Bigr\} < \delta,
\end{equation}
where $\hat{\bm{Q}}^* \coloneqq (\hat{Q}^*_i)_{i\in[I]}$ satisfying \eqref{eqn:lemma:exist_Q_star} with $\hat{\bm{Q}}^{\phi}=\hat{\bm{Q}}^*$, and $\bigl\{\hat{\bm{Q}}_t,t\in\mathbb{N}_0\bigr\}$ is the learning process accompanying the real process $\bigl\{\bs^{\psi}(t), t\in\mathbb{N}_0\bigr\}$.
\end{proposition}
The proof of Proposition~\ref{prop:q-convergence} is provided  in Appendix~\ref{app:prop:q-convergence}.

Proposition~\ref{prop:q-convergence} indicates that we can control the real process $\bigl\{\bs^{\psi}(t),t\in\mathbb{N}_0\bigr\}$ through $\actionpsi(t)$ while simultaneously learning the Q factors $\hat{\bm{Q}}^*$ for those MDPs with the same state and action spaces $\bS_i$ and $\bA_i$ ($i\in[I]$), respectively, transition matrices $\mathcal{P}_i(a)$ ($i\in[I],a\in\bA_i$) and reward functions $\bar{r}_i$ ($i\in[I]$) under a different policy $\phi\in\oPhiLocal$.
Also, for the same real process $\bigl\{\bs^{\psi}(t), t\in\mathbb{N}_0\bigr\}$, we can have $K$ different learning processes $\bigl\{\hat{\bm{Q}}^k_t,t\in\mathbb{N}_0\bigr\}$ with $\phi_t = \bar{\phi}_k\in\oPhiLocal$ and different reward functions, for which Proposition~\ref{prop:q-convergence} still holds.
More precisely, 
we can specify $\hat{\bm{Q}}_t = (\hat{\bm{Q}}^1_t;\hat{\bm{Q}}^2_t;\ldots\hat{\bm{Q}}^K_t)$.
For each $k\in[K]$, $\hat{\bm{Q}}^k_t=(\hat{Q}^k_{i,t})_{i\in[I]}$ with $\hat{Q}^k_{i,t}\in\bF(\bS_i\times\bA_i)$ updated through \eqref{eqn:q-iteration}  with specified $\hat{\functionQ}_{i,t} = \hat{Q}^k_{i,t}$, $\phi_t = \bar{\phi}_k$,  $\bar{r}_i=\bar{r}_i^k$, and the associated random rewards $\bar{R}_{i,n} = \bar{R}^{k}_{i,n}\in\bF(\bS_i\times\bA_i)$. 
In this context, we rewrite $\mathcal{T}^{\phi}_i$ defined in \eqref{eqn:define_H_phi} as $\mathcal{T}^{k,\phi}_i$ to indicate its dependency on $\bar{\bm{r}}_i^k$ and $\bar{R}^{k}_{i,n}$.
From Proposition~\ref{prop:q-convergence}, for any primary policy $\psi\in\lPhi$, $\epsilon,\delta>0$, and the step sizes $\bm{\eta}_t$ satisfying \eqref{eqn:assumption:stepsize}, there exist $T<\infty$ such that for all $k\in[K]$ and $t>T$, 
\begin{equation}
\mathbb{P}\Bigl\{\bigl\lVert \hat{\bm{Q}}_t^k - \hat{\bm{Q}}^{*,k}\bigr\rVert > \epsilon \Bigr\} < \delta,
\end{equation}
where, for $k\in[K]$, $\hat{\bm{Q}}^{*,k}\coloneqq (\hat{Q}^{*,k}_i)_{i\in[I]}$ satisfying $\bigl(\mathcal{T}_i^{k,\bar{\phi}_k}\hat{Q}^{*,k}_i\bigr)=\hat{Q}^{*,k}_i$ for all $i\in[I]$.
The existence of such $\hat{\bm{Q}}^{*,k}$ is ensured by Lemma~\ref{lemma:exist_Q_star}.
It enables approximation and implementation of those Q-factor-based policies in the scope of the WCG problem without assuming full knowledge of the transition matrices or reward functions. 

In the rest of this paper, we generalize the WCG-learning scheme to the one with $K\in\mathbb{N}_+$ learning processes, for which the estimated Q factor $\hat{\bm{Q}}_t = (\hat{\bm{Q}}^1_t;\hat{\bm{Q}}^2_t;\ldots;\hat{\bm{Q}}^K_t)$ and Lemma~\ref{lemma:exist_Q_star} and Proposition~\ref{prop:q-convergence} are applicable straightforwardly. 

We consider in total $K+1$ stochastic processes in parallel: the WCG process and the $K$ learning processes for estimating Q factors under user-defined secondary policies and reward functions.
The $K$ learning processes treat the WCG process as an event generator for ``simulating" real-world statistics, but evolve/iterate under the secondary policies which can be any policy defined by the human operator/administrator/user and not dependent on the employed primary policy for the WCG process.
Meanwhile, upon a decision making epoch, the primary policy that manages the WCG process makes decisions based on the current states/learned results of the $K$ learning processes.
In this way, the learning and control are being conducted simultaneously.

Recall that numbers of classic policies in the case of (multi-action) RMAB are computed based on Q factors under certain policies, such as the Whittle index policy~\cite{nino2007dynamic,nino2022multigear},  and a range of online algorithms, such as \cite{wang2023optimistic,jiang2023online,avrachenkov2022whittle,akbarzadeh2023learning}, are basically approximating the Whittle indices.
In Section~\ref{subsec:q-index-policy}, based on the past results on computing Whittle indices in~\cite{nino2001restless,nino2002dynamic,nino2006restless,nino2007dynamic,nino2022multigear}, we will provide detailed steps about how to approximate Whittle indices in an online manner.
Such an online algorithm will coincide with the offline Whittle index policy based on the convergence of the Q factors as $T\rightarrow \infty$.

Apart from the classic Whittle index policies, the WCG-learning scheme can also be used to approximate any Q-factor-based policies, since the $K$ learning processes are sufficiently flexible for any secondary policy and/or reward function.

\subsection{Fast Convergence Towards Asymptotic Optimality: Convergence in Magnitude Dimension}

\label{subsec:asym_regime}
Consider a parameter $h\in\mathbb{N}_+\cup\{\infty\}$ such that $N_i = h N_i^0$ for some $N_i^0 \in\mathbb{N}_+$ ($i\in[I]$). We refer to $h$ as the \emph{scaling parameter} of the WCG system. 
In this context, for $\phi\in\lPhi$, the constraints in \eqref{eqn:constraint:linear} can be rewritten as
\begin{equation}\label{eqn:constraint:linear:h}
\frac{1}{h}\sum_{i\in[I]}\sum_{n\in[N_i]}f_{i,\ell}(s^{\phi,h}_{i,n}(t),\actionphi^{h}_{i,n}(t)) =0,~\forall \ell\in[L],t\in[T]_0,
\end{equation}
where recall $L\in\mathbb{N}_+$, $T\in\mathbb{N}_+\cup\{\infty\}$, and $f_{i,\ell}
\in\bF(\bS_i\times\bA_i)$.
Recall that $f_{i,\ell}(s,a)$ is bounded.
Here, the $1/h$ are used to bound both sides of \eqref{eqn:constraint:linear:h} for all $h\in\mathbb{N}_+\cup \{\infty\}$.
We start with a finite-time horizon objective and later discuss  the long-run average case in Section~\ref{subsec:q-factor-based}. 
For $T <\infty$, discounting parameter $\beta\in(0,1]$ and initial state $\bs_0$, we aim to maximize
\begin{equation}\label{eqn:obj:h}
    \max_{\phi\in\lPhi} \frac{1}{h}\Gamma^{\beta,\phi,h}(T,\bs_0),
\end{equation}
subject to \eqref{eqn:constraint:linear:h}, where $\frac{1}{h}$ ensures finiteness of the maximum for all $h\in\mathbb{N}_+\cup\{\infty\}$, and the superscript $\phi$ and $h$ are attached to indicate their influence on $\Gamma^{\beta}(T,\bs_0)$.
We re-write the state and action vectors $\bs^{\phi}(t)$ and $\actionphi(t)$ as $\bs^{\phi,h}(t)$ and $\actionphi^{h}(t)$, respectively, to indicate effects of $h$.

Define $\mathcal{I}\coloneqq \sum_{i\in[I]}\bigl|\bS_i\times \bA_i\bigr|$, and we label all the SA pairs $(i,s,a)$ by $\iota \in[\mathcal{I}]$.
Let $(i_{\iota},s_{\iota},\actionphi_{\iota})$ represent the SA pair labeled by $\iota\in[\mathcal{I}]$, and we alternatively refer to SA pair $(i_{\iota},s_{\iota},\actionphi_{\iota})$ as SA pair $\iota$.
For $h\in\mathbb{N}_+$, $\phi\in\lPhi$, $\iota\in[\mathcal{I}]$, and $t\in[T]_0$,
define $Z^{\phi,h}_{\iota}(t)$ as the proportion of processes $\bigl\{s^{\phi,h}_{i,n}(r), t\in[T]_0\bigr\}$ ($i\in[I],n\in[N_i]$) that are in SA pair $\iota$ at time $t$; that is, 
\begin{equation}\label{eqn:define_Z}
Z^{\phi,h}_{\iota}(t) \coloneqq \frac{1}{\sum_{i\in[I]}N_i}\biggl|\Bigl\{n\in[N_{i_{\iota}}]~|~s^{\phi,h}_{i_{\iota},n}(t) = s_{\iota}, \actionphi^{h}_{i_{\iota},n}(t)=\actionphi_{\iota}\Bigr\}\biggr|,
\end{equation}
for which define $\bm{Z}^{\phi,h}(t) \coloneqq \bigl(Z^{\phi,h}_{\iota}(t):\iota\in[\mathcal{I}]\bigr)$.
For the maximization in \eqref{eqn:obj:h}, since the bandit processes associated with $(i,n)$ gain the same expected reward $r_i(s,a)$ and transition rates $\mathcal{P}_i(a)$ when they are in the same SA pair $(i,s,a)$ ($i\in[I]$, $(s,a)\in\bS_i\times\bA_i$), we say that such bandit processes in the same SA pair are \emph{identical}.
In this context, for each time $t$,  the expected total reward $\sum_{i\in[I]}\sum_{n\in[N_i]}\mathbb{E}\bigl[R^{\phi,h}_{i,n}(t)\bigr] = \sum_{\iota\in[\mathcal{I}]}\sum_{i\in[I]}N_i Z^{\phi,h}_{\iota}(t)r_{i_{\iota}}(s_{\iota},\actionphi_{\iota})$, where $R^{\phi,h}_{i,n}(t)$ is the instantaneous reward generated by process $\bigl\{s^{\phi,h}_{i,n}(t), t\in[T]_0\bigr\}$ with scaling parameter $h$.
The problem described in \eqref{eqn:obj:h} and \eqref{eqn:constraint:linear:h} is equivalent to
\begin{equation}\label{eqn:obj:Z}
\max_{\phi\in\lPhi}\sum_{t=0}^T\beta^t \sum_{\iota\in[\mathcal{I}]}Z^{\phi,h}_{\iota}(t)r_{i_{\iota}}(s_{\iota},\actionphi_{\iota}),
\end{equation}
subject to 
\begin{equation}\label{eqn:constraint:linear:Z}
 \sum_{\iota\in[\mathcal{I}]}Z^{\phi,h}_{\iota}(t)f_{i_{\iota},\ell}(s_{\iota},\actionphi_{\iota})=0,~\forall \ell\in[L],t\in[T]_0.
\end{equation}
We say that the process $\bigl\{\bs^{\phi,h}(t),t\in[T]_0\bigr\}$ can be interpreted by $\bigl\{\bm{Z}^{\phi,h}(t),t\in[T]_0\bigr\}$.

Define a \emph{ranking} (permutation) $\bm{\rank}\coloneqq(\rank_1,\rank_2,\ldots,\rank_{\mathcal{I}})$ of all the SA pairs $\iota\in[\mathcal{I}]$, where $\rank_{\iota}\in[\mathcal{I}]$ represent the ranking of SA pair $\iota$ with respect to $\bm{\rank}$.
Let $\mathscr{R}$ be the set of all such rankings.
Consider a subset $\Phi^1\subset \lPhi$ of policies, where, for any $\varphi\in\Phi^1$, the action variable $\actionvarphi^{h}_{i,n}(t)$ is dependent on the history of the WCG-learning process through only the current state $\bs^{\varphi,h}(t)$, time stamp $t$, and an SA-pair ranking $\bm{\mathcal{R}}_t=\bm{\rank}_K(\hat{\bm{Q}}^1_t,\hat{\bm{Q}}^2_t,\ldots,\hat{\bm{Q}}^K_t)$, where $\bm{\rank}_K:\bigl(\mathbb{R}^{\sum_{i\in[I]}|\bS_i||\bA_i|}\bigr)^K\rightarrow \mathscr{R}$ is a function of $\hat{\bm{Q}}_t = (\hat{\bm{Q}}^1_t;\hat{\bm{Q}}^2_t;\ldots\hat{\bm{Q}}^K_t)$.
Let $\bm{z}^{\varphi,h}(t)\coloneqq\mathbb{E}[\bm{Z}^{\varphi,h}(t)]$. We have the following theorem.
\begin{theorem}\label{theorem:convergence-Z}
For any primary policy $\varphi\in\Phi^1$, secondary policies $\bar{\phi}^k_t\in\oPhiLocal$ ($k\in[K]$), $T<\infty$, $\epsilon>0$, and given initial state $\bm{Z}^{\varphi,h}(0)=\bm{z}_0$ and $\hat{\bm{Q}}_0=\bm{Q}_0$, the limit $\lim_{h\rightarrow \infty}\bm{z}^{\varphi,h}(t)$ exists, and
\begin{equation}\label{eqn:theorem:convergence-Z}
    \lim_{h\rightarrow\infty} \mathbb{P}\Bigl\{\max_{t=0,1,\ldots,T}\bigl\lVert\bm{Z}^{\varphi,h}(t) - \lim_{h\rightarrow\infty}\bm{z}^{\varphi,h}(t)]\bigr\rVert > \epsilon\Bigr\} = 0.
\end{equation}
\end{theorem} 
The proof of Theorem~\ref{theorem:convergence-Z} is provided in Appendix~\ref{app:theorem:convergence-Z}.

Theorem~\ref{theorem:convergence-Z} indicates that, for any primary policy $\varphi\in\Phi^1$, increasing the size of the WCG problem, measured by $h$, the process $\bigl\{\bm{Z}^{\varphi,h}(t), t\in[T]_0\bigr\}$ converges to a deterministic averaging trajectory $\bm{z}^{\varphi}(t)\coloneqq \lim_{h\rightarrow \infty}\bm{z}^{\varphi,h}(t)$.
More importantly,  based on Freidlin's theorem (\cite{freidlin2012random}), we prove that the convergence speed is exponential in $h$. 
\begin{theorem}\label{theorem:convergence_Z_exp}
For any primary policy $\varphi\in\Phi^1$, secondary policies $\bar{\phi}^k_t\in\oPhiLocal$ ($k\in[K]$), $T<\infty$, $\epsilon>0$, and given initial state $\bm{Z}^{\varphi,h}(0)=\bm{z}_0$ and $\hat{\bm{Q}}_0=\bm{Q}_0$, there exist positive constants $C,H<\infty$ such that, for all $h>H$,
\begin{equation}\label{eqn:theorem:convergence_Z_exp}
    \mathbb{P}\Bigl\{\max_{t=0,1,\ldots,T}\bigl\lVert\bm{Z}^{\varphi,h}(t) - \bm{z}^{\varphi}(t)\bigr\rVert > \epsilon\Bigr\} \leq e^{-Ch}.
\end{equation}
\end{theorem}
The proof of Theorem~\ref{theorem:convergence_Z_exp} is provided in Appendix~\ref{app:theorem:convergence_Z_exp}.

Theorem~\ref{theorem:convergence_Z_exp} strengthens Theorem~\ref{theorem:convergence-Z} and ensures a fast convergence speed to the asymptotic regime with respect to $\bm{Z}^{\varphi,h}(t)$.
In particular, 
when $\bm{Z}^{\varphi,h}(t)$ is sufficiently close to the averaging trajectory $\bm{z}^{\varphi,h}(t)$ for all $t$, we can conclude the Corollary~\ref{eqn:convergence_transition_prob:1} of Theorem~\ref{theorem:convergence_Z_exp}.

For $k\in[K]$, $i\in[I]$, $(s,a)\in\bS_i\times\bA_i$, $t\in\mathbb{N}_0$, $\bm{Q}=(Q_i)_{i\in[I]}\in \prod_{i\in[I]}\bF(\bS_i\times\bA_i)$, a primary policy $\varphi\in\Phi^1$, and a secondary policy $\bar{\phi}_k\in\oPhiLocal$, we define
\begin{multline}\label{eqn:define_general_w}
w^{k,\bar{\phi}_k}_{i,t}(s,a,\bm{Q}) \coloneqq \frac{1}{|\hat{\mathscr{N}}^{\varphi,h}_{i,t}(s,a)|}\Bigl(\sum_{n\in\hat{\mathscr{N}}_{i,t}^{\varphi,h}(s,a))}\bar{R}^{k}_{i,n}(s,a)+\sum_{s'\in\bS_i\backslash\{s_0\}}|\hat{\mathscr{N}}^{\varphi,h}_{i,t}(s,a,s')|Q_i(s',\actionbarphi_{k,i}(s'))\Bigr)\\
- \bar{r}^k_i(s,a) - \sum_{s'\in\bS_i\backslash\{s_0\}} p_i(s,a,s')Q_i(s',\actionbarphi_{k,i}(s')),
\end{multline}
for which $w_{i,t}(s,a)$ defined in \eqref{eqn:define_w} is a special case of $w^k_{i,t}(s,a,\bm{Q}) $ with specified $\bm{Q} = \hat{\bm{Q}}^k_t$, $\phi=\varphi$, and $\bar{\phi}_t=\bar{\phi}_k$.
Let $\bm{w}^{k,\actionbarphi_k}_t(\bm{Q}) \coloneqq (w^{k,\actionbarphi_k}_{i,t}(s,a,\bm{Q}):i\in[I],s\in\bS_i,a\in\bA_i)$.
\begin{corollary}\label{coro:converge_transition_prob}
For any primary policy $\varphi\in\Phi^1$, secondary policies $\bar{\phi}^k_t\in\oPhiLocal$ ($k\in[K]$), $T<\infty$, $\epsilon>0$, and given initial state $\bm{Z}^{\varphi,h}(0)=\bm{z}_0$ and $\hat{\bm{Q}}_0=\bm{Q}_0$,   there exist $C>0$ and $H<\infty$ such that, for all $h>H$ and $\iota\in[\mathcal{I}]$,
\begin{equation}\label{eqn:convergence_transition_prob:1}
\mathbb{P}\Bigl\{\max_{t=0,1,\ldots,T}Z^{\varphi,h}_{\iota}(t)\sum_{s'\in\bS_{i_{\iota}}}\Bigl|\frac{|\hat{\mathscr{N}}^{\varphi,h}_{i_{\iota},t}(s_{\iota},\actionphi_{\iota},s')|}{|\hat{\mathscr{N}}^{\varphi,h}_{i_{\iota},t}(s_{\iota},\actionphi_{\iota})|} - p_{i_{\iota}}(s_{\iota},\actionphi_{\iota},s')\Bigr|>\epsilon\Bigr\}\leq e^{-Ch},
\end{equation}
and, for all $k\in[K]$ and any given $\bm{Q}\in \prod_{i\in[I]}\bF(\bS_i\times\bA_i)$,
\begin{equation}\label{eqn:convergence_transition_prob:2}
\lim_{h\rightarrow \infty}\mathbb{P}\Bigl\{\max_{t=0,1,\ldots,T}\bm{Z}^{\varphi,h}(t)\cdot|\bm{w}^{k,\bar{\phi}_k}_t(\bm{Q})| >\epsilon\Bigr\} = 0,
\end{equation}
where $|\bm{w}^{k,\bar{\phi}_k}_t(\bm{Q})| = (|w^{k,\bar{\phi}_k}_{i,t}(s,a,\bm{Q})|: i\in[I], s\in\bS_i, a\in\bA_i)$ with
$w^k_{i,t}(s,a,\bm{Q})$ defined in \eqref{eqn:define_w}, and $\bm{Z}^{\varphi,h}(t)\cdot|\bm{w}^{k,\bar{\phi}_k}_t(\bm{Q})|$ is the dot production of the two vectors.

If we also assume, for any $\epsilon>0$ and $N^0\in\mathbb{N}_+$, there exists constant $C>0$ such that, for any $i\in[I]$ and $(s,a)\in\bS_i\times\bA_i$, given $s^{\varphi,h}_{i,n}(t)=s$ and $\actionvarphi^{h}_{i,n}(t)=a$ for all $n\in[hN^0_i]$,
\begin{equation}\label{eqn:assumption:R_exp}
\lim_{h\rightarrow \infty} \frac{1}{h} \ln\mathbb{P}\Bigl\{\Bigl\lvert\frac{\sum_{n\in hN^0_i}\bar{R}^{k}_{i,n}(s,a)}{hN^0_i}-\bar{r}^k_{i}(s,a)\Bigr\rvert > \epsilon\Bigr\}   = -C,
\end{equation}
then for any $\epsilon>0$, there exist $C>0$ and $H < \infty$ such that, for all $h > H$ and $\iota\in[\mathcal{I}]$, 
\begin{equation}\label{eqn:convergence_transition_prob:3}
\mathbb{P}\Bigl\{\max_{t=0,1,\ldots,T}\bm{Z}^{\varphi,h}(t)\cdot|\bm{w}_t| >\epsilon\Bigr\}\leq e^{-Ch}.
\end{equation}

\end{corollary}
The proof of Corollary~\ref{coro:converge_transition_prob} is provided in Appendix~\ref{app:coro:convergence_transition_prob}.

From Corollary~\ref{coro:converge_transition_prob}, for any SA pair $\iota\in[\mathcal{I}]$ and time $t\in[T]_0$, if it has positive probability $\mathbb{P}\bigl\{s^{\varphi,h}_{i_{\iota},n}(t) = s_{\iota},\actionvarphi^{h}_{i_{\iota},n}(t)=\actionphi_{\iota}~\bigl|~\bm{Z}^{\varphi,h}(0)=\bm{z}_0, \hat{\bm{Q}}_0=\bm{Q}_0 \bigr\} =z^{\varphi,h}_{\iota}(t) > 0$, then, the estimated transition probabilities $|\hat{\mathscr{N}}^{\varphi,h}_{i_{\iota},t}(s_{\iota},\actionphi_{\iota},s')|/|\hat{\mathscr{N}}^{\varphi,h}_{i_{\iota},t}(s_{\iota},\actionphi_{\iota})|$ and the perturbations of the learning processes $w_{i_{\iota},t}(s_{\iota},\actionphi_{\iota})$ quickly converge to the real transition kernel $p_{i_{\iota}}(s_{\iota},\actionphi_{\iota},s')$ and zero, respectively, as $h\rightarrow \infty$. 
Increasing the size of the WCG-learning system, measured by the scaling parameter $h$, can lead to fast increasing accuracy of the estimated transition probabilities and the learned Q factors.

Equation~\eqref{eqn:assumption:R_exp} is a mild assumption that requests $\mathbb{E}\bigl[(\bar{R}^{k}_{i,n}(s,a))^2\bigr]$ to be somehow bounded.
For instance, for $(s,a)\in\bS_i\times\bA_i$, if $\bar{R}^{k}_{i,n}(s,a)$ follows a normal distribution with finite variance, then \eqref{eqn:assumption:R_exp} is satisfied.

Given the high accuracy of estimating unknown parameters, we are more interested in the overall performance of the algorithms potentially produced;
that is, how tightly that the convergence of the WCG process $\bm{Z}^{\varphi,h}(t)$ affects the overall performance of the potential algorithms?
In Sections~\ref{sec:long-run} and \ref{sec:LP-algo}, we will discuss such tightness for two main directions of restless-bandit-based studies: Q-factor-based index policies for long-run objectives, such as~\cite{whittle1988restless,nino2007dynamic,fu2018restless}, and linear-programming-based approximations for finite-time-horizon cases, such as~\cite{brown2020index, gast2024reoptimization,fu2024patrolling}.



In Sections~\ref{sec:long-run} and \ref{sec:LP-algo}, we theoretically demonstrate that, for the WCG-learning scheme, the performance degradation of the Q-factor-based algorithms and Linear-programming-based approximations can diminish exponentially in the scaling parameter $h$.

\section{Long-Run Optimization: Q-factor-based algorithms}\label{sec:long-run}

\subsection{Q-Factor-Based Index Policies}\label{subsec:q-index-policy}

Consider a special case of the WCG problem with $L=1$.
For $i\in[I]$, we label all the actions $a\in\bA_i$ as $\La = 1,2,\ldots,|\bA_i|$, and consider the following condition.
\begin{condition}{Single Constraint}\label{condition:single-constraint} 
For $i\in[I]$, $(s,a)\in\bS_i\times\bA_i$, and $\ell=L=1$,
$f_{i,\ell}(s,a) = f_i(\La)$ for a function $f_i(\La):\bigl[|\bA_i|\bigr] \rightarrow \mathbb{R}_0$ that satisfies
\begin{equation}
0\leq f_i(1)< f_i(2) < \ldots < f_i(|\bA_i|) <\infty.
\end{equation}\vspace{-1cm}
\end{condition}
The \partialref{condition:single-constraint}{Single Constraint} condition is the same as the \emph{multi-gear} assumption in \cite{nino2022multigear}.
If a WCG problem satisfies \partialref{condition:single-constraint}{Single Constraint} condition, then we say it reduces to a \emph{multi-action RMAB} problem. 
For a multi-action RMAB problem, $\phi\in\oPhiLocal$, $i\in[I]$, and $s\in\bS_i$, let $\La^{\phi}_i(s)$ represent the action label for $\actionphi_i(s)$.

Consider the conventional RMAB problem as a special case of the multi-action RMAB, where $\bA_i = \{0,1\}$.
In this case, the Whittle index policy (\cite{whittle1988restless}) that prioritizes bandit processes $(i,n)$ according to the descending order of $\gamma_i(s^{\phi}_{i,n}(t))$ (referred to as the Whittle indices) have demonstrated advantages in general and have been proved to approach optimality as $N_i\rightarrow \infty$ proportionally for all $i\in[I]$ under a range of scenarios (see \cite{weber1990index,verloop2016asymptotically,fu2020energy,fu2018restless,fu2024patrolling}).
The Whittle indices can be exactly computed, e.g., through the adaptive-greedy index algorithm under so-called PCL-indexability conditions, see \cite{ninomora2001restless,nino2002dynamic,nino2006restless,nino2007dynamic}. 
For the more general RMAB problem with multiple actions with $|\bA_i| \geq 2$, \cite{weber2007comments} extended the Whittle indices to the multi-action RMAB problem, and \cite{nino2008index} outlined the extension of PCLs to the multi-action RMAB case, which was subsequently analyzed thoroughly in \cite{nino2022multigear}.
When the RMAB is PCL indexable (see \cite{nino2002dynamic,nino2022multigear}), the MP indices coincide with the Whittle indices.


For the case without knowledge of the transition matrices and reward functions a priori, we propose in the following a scheme consisting of a primary policy and a sequence of coordinated secondary policies.
Based on Proposition~\ref{prop:q-convergence}, the proposed scheme will coincide with the offline MP index policy (see \cite{nino2007dynamic,nino2022multigear}) with respect to the long-run average objective defined in \eqref{eqn:obj:long-run-average}.
In Section~\ref{subsec:q-factor-based}, we will explain that, based on Theorem~\ref{theorem:convergence_Z_exp}, such a convergence can be significantly stimulated in the magnitude dimension.

\subsubsection{Offline MP Indices}\label{subsubsec:MP-indices}
Define the long-run average reward of the class-$i$ process $\bigl\{s^{\phi}_{i,n},t\in\mathbb{N}_0\bigr\}$ (for any $n\in[N_i]$) under a policy $\phi\in\oPhiLocal$ as 
\begin{equation}\label{eqn:obj_i}
\Gamma^{\phi}_i\coloneqq \limsup_{T\rightarrow \infty}\frac{1}{T}\sum_{t\in[T]} R^{\phi}_{i,n}(t) = \frac{V^{\phi}_i(s_0)}{L^{\phi}_i(s_0)},
\end{equation}
where $s_0\in\bS_i$ is the ergodic state, $V^{\phi}_i(s_0)$ and $L^{\phi}_i(s_0)$ are the expected cumulative reward and time, respectively, of the process $\bigl\{s^{\phi}_{i,n},t\in\mathbb{N}_0\bigr\}$ when it starts from state $s_0$ until re-enters $s_0$.
If the transition probability $p_i(s_0,\actionphi_i(s_0),s_0) = 1$ (that is, the underlying process is trapped in state $s_0$), then $V^{\phi}_i(s_0)/L^{\phi}_i(s_0) = r_i(s_0,\actionphi_i(s_0))$ with $V^{\phi}_i(s_0),L^{\phi}_i(s_0)\rightarrow \infty$.
The second equality in \eqref{eqn:obj_i} is based on \cite{ross1992applied}: given a criteria $g=\Gamma^{\phi}_i\in\mathbb{R}$, the value functions $V^{\phi,g}_i(s)$ (for $s\in\bS_i$) of the underlying process with expected reward rate $r_i(s,\actionphi_i(s))-g$ for state $s$ satisfy Bellman equation
\begin{equation}
V^{\phi,g}_i(s) = r_i(s,\actionphi_i(s)) - g + \sum_{s'\in\bS_i\backslash \{s_0\}} p_i(s,\actionphi_i(s),s') V^{\phi,g}_i(s'),
\end{equation}
for all $s\in\bS_i$, and 
\begin{equation}
    V^{\phi,g}_i(s_0) = \sum_{t=0}^{L^{\phi}_i(s_0)}\Bigl[R^{\phi}_{i,n}(t)-g\Bigr]=V^{\phi}_i(s_0)-gL^{\phi}_i(s_0)=0.
\end{equation}
Similarly, define the long-run average \emph{marginal cost} of the class-$i$ process under a policy $\phi\in\oPhiLocal$ as
\begin{equation}\label{eqn:marginal-cost_i}
\Omega^{\phi}_i\coloneqq \limsup_{T\rightarrow \infty}\frac{1}{T}\sum_{t\in[T]} f_i\bigl(\La^{\phi}_i(s^{\phi}_{i,n}(t))\bigr) = \frac{U^{\phi}_i(s_0)}{L^{\phi}_i(s_0)},
\end{equation}
where $U^{\phi}_i(s_0)$ is the expected cumulative cost of the process $\bigl\{s^{\phi}_{i,n}(t),t\in\mathbb{N}_0\bigr\}$ since it starts from state $s_0$ until re-enters $s_0$, where the cost rate for state $s\in\bS_i$ is $f_i\bigl(\La^{\phi}_i(s^{\phi}_{i,n}(t))\bigr)$.
Recall that for $\phi\in\oPhiLocal$, $\La^{\phi}_i(s)$ is the label for action $\actionphi_i(s)$.

\cite{nino2022multigear} generalized the adaptive-greedy algorithm for computing MP indices (Whittle indices) to the multi-action RMAB case and referred to it as the \emph{downshift adaptive-greedy (DS) algorithm}. 
We start with introducing the offline DS algorithm proposed in \cite{nino2022multigear} for the multi-action RMAB problem.
Given $i\in[I]$ and $\LaVec=(\La(s):s\in\bS_i)\in\bigl[|\bA_i|\bigr]^{|\bS_i|}$,
we denote a policy $\phi\in\oPhiLocal$ with $\actionphi_i(s)$ labeled by $\La(s)$ as $\phi(\LaVec)$, and, for $s\in\bS_i$, consider an iteration of the action label vector: for $s'\in\bS_i$,
\begin{equation}
    (\Ia_i^{s}\LaVec)_{s'} \coloneqq \begin{cases}
        \La_{s'}-1, & \text{if }s = s', \La_{s'}\geq 2,\\
        \La_{s'}, & \text{otherwise}.
    \end{cases}
\end{equation}
Let $\LenM_i \coloneqq |\bS_i||\bA_i|$.
For any $i\in[I]$, based on \cite{nino2022multigear}, if the bandit process $i$ is PCL indexable, then there exist a sequence of SA pairs $(s_m,\La_m)$, $m=1,2,\ldots, \LenM_i$, and a sequence of action label vectors $\LaVec_m$ with $\LaVec_{m+1} = \Ia_i^{s_m}\LaVec_m$, 
$\LaVec_1 = |\bA_i|\bm{1}^{|\bS_i|}$ and $\LaVec_{\LenM_i} = \bm{1}^{|\bS_i|}$, where $\bm{1}^n$ is a vector of $n$ elements equal to $1$, such that the MP indices are given by, for $m=1,2,\ldots,\LenM_i-1$, 
\begin{equation}\label{eqn:define_index_multi_action}
    \nu_i(s_m,\La_m)\coloneqq \frac{\Gamma^{\phi(\LaVec_m)}_i-\Gamma^{\phi(\LaVec_{m+1})}_i}{\Omega^{\phi(\LaVec_m)}_i-\Omega^{\phi(\LaVec_{m+1})}_i},
\end{equation}
In particular, under the PCL-indexability, 
\[\nu_i(s_1,\La_1) \geq \nu_i(s_2,\La_2)\geq \ldots \geq \nu_i(s_{\LenM_i-1},\La_{\LenM_i-1}).\]
We provide the pseudo-code for finding such a sequence through the offline DS algorithm in Algorithm~\ref{algo:DS-agi}.

Given the MP indices, at each decision epoch $t$, we rank the pair of the bandit processes $s^{\phi}_{i,n}(t)$ and their potential actions $\La$ according to the descending order of $\nu_i(s^{\phi}_{i,n}(t),\La)$ for all $\La=2,3,\ldots,|\bA_i|$, $i\in[I]$, and $n\in[N_i]$. 
In the event of a tie with the same index value, prioritize one at random.
We use the ranks, $\rank =1,2,\ldots,\sum_{i\in[I]}N_i(|\bA|-1)$, of these bandit processes to label them; that is, we write the $\rank$th bandit process and its potential action as $(i(\rank),n(\rank),\La(\rank))$. 
We use the action label $\La^{\phi}_{i,n}(t)$ to represent the action $\actionphi_{i,n}(t)$ under policy $\phi$ for the bandit process $s^{\phi}_{i,n}(t)$ at time $t$.
Initialize all action labels $\La^{\phi}_{i,n}(t)= 1$ (the passive action).
From $\rank = 1$ to $\sum_{i\in[I]}N_i(|\bA|-1)$, 
we check the single constraint for the action variables
\begin{equation}\label{eqn:DS_policy:singl_constraint}
\sum_{i\in[I]}\sum_{n\in[N_i]}f_i\bigl(\La^{\phi}_{i,n}(t)\bigr) + f_{i(\rank),n(\rank)}(\La(\rank)) - f_{i(\rank),n(\rank)}\bigl(\La^{\phi}_{i(\rank),n(\rank)}(t)\bigr)\leq 0.
\end{equation}
If \eqref{eqn:DS_policy:singl_constraint} holds, then set $\La^{\phi}_{i(\rank),n(\rank)}(t) = \La(\rank)$.
After we check all the pairs of the bandit process and potential action pairs from $\rank=1$ to $\sum_{i\in[I]}N_i(|\bA|-1)$,
we use the resulting values of $\La^{\phi}_{i,n}(t)$ for all $i\in[I]$ and $n\in[N_i]$ as the decided actions to be implemented at time $t$.
We write such actions being implemented when given $\bm{s}^{\phi}(t) = \bm{s}$ and the MP indices $\bm{\nu}\coloneqq \bigl(\nu_i(s,\La):i\in[I],s\in\bS_i,\La\in\bigl[|\bA_i|\bigr]\backslash\{1\}\bigr)$ as $\LaVec^{\text{MP}}(\bm{\nu},\bm{s})$, which is a function 
\[\LaVec^{\text{MP}}:~\mathbb{R}^{\sum_{i\in[I]}|\bS_i|(|\bA_i|-1)}\times\prod_{i\in[I]}\bS^{N_i}_i\rightarrow \prod_{i\in[I]}\bigl[|\bA_i|\bigr]^{N_i}.\]
Such a decision making policy is the offline MP index policy for the multi-action RMAB problem.
Recall that, based on \cite{nino2022multigear}, when the bandit processes are PCL-indexable, MP indices and the MP index policy coincide with the Whittle indices and the Whittle index policy, respectively.

\IncMargin{1em}
\begin{algorithm}
\small 
\linespread{0.4}\selectfont

\SetKwFunction{FDownshiftAdaptiveGreedyIndex}{DownshiftAdaptiveGeedyIndex}
\SetKwProg{Fn}{Function}{:}{\KwRet}
\SetKwInOut{Input}{Input}\SetKwInOut{Output}{Output}
\SetAlgoLined
\DontPrintSemicolon
\Input{Given $i\in[I]$}
\Output{$(s_m,\La_m,\LaVec_m)$ and $\nu_i(s_m,a_m)$ for all $m\in[\LenM_i-1]$.}
\Fn{\FDownshiftAdaptiveGreedyIndex{}}{
    $\LaVec\gets |\bA_i|\bm{1}^{|\bS_i|}$\;
    \For{$m=1,2,\ldots,\LenM_i-1$}{
        \For{$s\in\{s'\in\bS_i|\La(s')\geq 2\}$}{
            $\LaVec' \gets \Ia_i^s \LaVec$\;
            \eIf{$\Omega^{\phi(\LaVec)}_i\neq\Omega^{\phi(\LaVec')}_i$}{                              
                $\nu(s)\gets \frac{\Gamma^{\phi(\LaVec)}_i-\Gamma^{\phi(\LaVec')}_i}{\Omega^{\phi(\LaVec)}_i-\Omega^{\phi(\LaVec')}_i}$\;
            }{
                $\nu(s)\gets 0$\;
            }
        }
        $s_m\gets \arg\max_{\begin{subarray}~s\in\bS_i\\\La(s)\geq 2\end{subarray}}\nu(s)$\;
        $\La_m\gets \La(s_m)$\;
        $\nu_i(s_m,\La_m) \gets \nu(s_m)$\;
        $\LaVec \gets \Ia^{s_m}_i\LaVec$\;
    }
}
\caption{Downshidt (DS) Adaptive-greedy index algorithm}\label{algo:DS-agi}
\end{algorithm}
 \DecMargin{1em}

\subsubsection{Online MP Index Policy}
\label{subsubsec:online-MP}
In the online case without knowing $\mathcal{P}_i(a)$ and $r_i(s,a)$ for $s\in\bS_i$ and $a\in\bA_i$ a priori, we can approximate $\Gamma^{\phi}_i$ and $\Omega^{\phi}_i$ by estimating $V^{\phi}_i(s_0)$, $U^{\phi}_i(s_0)$, and $L^{\phi}_i(s_0)$ through the WCG-learning scheme and, simultaneously, implementing the MP index policy based on the up-to-date data.

Define $\hatnu_{i,t}(s,\La)\in\mathbb{R}$ as the estimated MP index for $i\in[I]$ and $(s,\La)\in\bS_i\times\bigl[|\bA_i|\bigr]$ at time $t\in[T]_0$.
Let $\hatnuVec_t \coloneqq (\hatnu_{i,t}(s,\La): i\in[I],s\in\bS_i,\La\in\bigl[|\bA_i|\bigr])$.
Consider a primary policy $\psi \in \lPhi$ and a sequence of parameters $\bar{p}_t\in[0,1)$ ($t\in[T]_0$).
With probability $\bar{p}_t$, the policy $\psi$ take actions $\LaVec^{\psi}(t) = \LaVec^{\text{MP}}\bigl(\hatnuVec_t,\bs^{\psi}(t)\bigr)$.
With probability $1-\bar{p}_t$, we initialize $\LaVec^{\psi}(t) = \bm{1}^{\sum_{i\in[I]}N_i}$ and a set of pairs of bandit process and its action $\bB = \prod_{i\in[I]}[N_i]\times \bigl([|\bA_i|]\backslash\{1\}\bigr)$.
We uniformly randomly select a pair of bandit process and its action $(i,n,\La)$ from $\bB$.
If \eqref{eqn:DS_policy:singl_constraint} is not violated by replacing $\La(\rank)$ with $\La$, then take action $\La^{\psi}_{i,n}(t) = \La$, and updating $\bB$ with $\bB\backslash\{(i,n,\La)\}$.
We repeat such a random selection process until \eqref{eqn:DS_policy:singl_constraint} is violated by replacing $\La(\rank)$ with $\La$. 
The resulting $\LaVec^{\psi}(t)$ is the decided vector of actions for time $t$. 
The pseudo-code of implementing $\psi$ with updated $\hatnuVec_t$ is described in Lines~\ref{line:psi:begin}-\ref{line:psi:end} of Algorithm~\ref{algo:ompi} in Appendix~\ref{app:ompi} as part of the online MP index policy.

Here, the parameter $\bar{p}_t$ is measurable with respect to the history of the WCG learning process $\sH_t$ and is determined right before taking actions at time $t$.
The parameter $\bar{p}_t\in[0,1)$ is used to ensure a positive probability of visiting every SA pair $(s,a)\in\bS_i\times\bA_i$ ($i\in[I]$) within a finite time horizon.

Consider the WCG-learning system consisting of the WCG process $\bigl\{\bs^{\phi}(t), t\in[T]_0\bigr\}$ and $K= 6\max_{i\in[I]}|\bS_i|(|\bA_i|-1)$ learning processes that are used to estimate the MP indices for all $in\in[I]$ and $(s,\La)\in\bS_i\times\bigl(\bigl[|\bA_i|\bigr]\backslash\{1\}\bigr)$ in parallel.
Similar to the offline DS algorithm in Algorithm~\ref{algo:DS-agi}, we update the estimated MP indices $\hatnuVec_t$ through the following steps.
\begin{enumerate}[label=\roman*)]
\item \label{step:1}Initialize $\LaVec_{i,1}=|\bA_i|\bm{1}^{|\bS_i|}$ for all $i\in[I]$, $m=1$, and $\hatnuVec_0 = \bm{\nu}_0$ with given $\bm{\nu}_0\in\mathbb{R}^{\prod_{i\in[I]}|\bS_i|(|\bA_i|-1)}$.
\item \label{step:period}Let $J = \max_{i\in[I]}\bigl|\bS_i\bigr|\bigl(|\bA_i|-1\bigr)$ and a matrix $\sS_m\coloneqq[\bm{s}_1;\bm{s}_2;\ldots,\bm{s}_I]^T\coloneqq[s_{i,j}]_{I\times J}$, where  $\bm{s}_1, \bm{s}_2,\ldots,\bm{s}_I$ are column vectors satisfying that, for any $i\in[I]$, $\bm{s}_i = (s_{i,1},s_{i,2},\ldots,s_{i,\bigl|\bS_i\bigr|}; \bm{0})$ with $s_{i,1},s_{i,2},\ldots,s_{i,\bigl|\bS_i\bigr|}$ a permutation of all elements in $\bS_i$ and $\bm{0}$ is a vector of $J-|\bS_i|$ zeros.
\item \label{step:3} For all $j\in[J]$, construct policies $\bar{\phi}^0_m,\bar{\phi}_m^1\in\oPhiLocal$ such that, for all $i\in[I]$ and $j=1,2,\ldots, |\bS_i|$,
$ \La^{\bar{\phi}^0_m}_i(s_{i,j}) = \La_{i,m}$ and $\La^{\bar{\phi}^1_m}_i(s_{i,j}) = \bigl(\Ia_i^{s_{i,j}}\La_{i,m}\bigr)_{i,j}$.
\item\label{step:4} \label{step:q-iteration}At time $t$, iterate the $6J$ learning processes  $\hat{\bm{Q}}_t^{V,j,\varsigma},\hat{\bm{Q}}_t^{U,j,\varsigma},\hat{\bm{Q}}_t^{L,j,\varsigma}$ with $j\in[J]$ and $\varsigma\in\{0,1\}$ according to \eqref{eqn:q-iteration} under the secondary policies $\bar{\phi}_m^\varsigma$.
The first superscript, $V,U,L$, in the estimated Q factors indicate different reward rates of the corresponding learning processes; that is, $\bar{R}^{\cdot}_{i,n}$ for the learning processes $\hat{\bm{Q}}_t^{V,j,\varsigma}$, $\hat{\bm{Q}}_t^{U,j,\varsigma}$, and $\hat{\bm{Q}}_t^{L,j,\varsigma}$ are set to be $R^{\bar{\phi}_m^\varsigma}_{i,n}$, $ f_i\bigl(\La^{\bar{\phi}_m^\varsigma}_{i,n}(t)\bigr)$, and constant $1$, respectively (with $\bar{r}_i(s,a)$ set to be $r_i(s,a)$, $f_i(\La)$, and $1$, respectively).
\item \label{step:5}For $i\in[I]$ and $j\in\Bigl[\min\bigl\{J,|\bS_i|\bigr\}\Bigr]$, update the estimated MP indices: if $\La_{i,m}(s_{i,j})\geq 2$, then
\begin{multline}
    \hat{\nu}_{i,t}(s_{i,j},\La_{i,m}(s_{i,j})) = \\\begin{cases}
        \frac{\hat{Q}^{V,j,0}_{i,t}(s_0,a_i(0))\Bigl/\hat{Q}^{L,j,0}_{i,t}(s_0,a_i(0))~-~\hat{Q}^{V,j,1}_{i,t}(s_0,a_i(1))\Bigl/\hat{Q}^{L,j,1}_{i,t}(s_0,a_i(1))}{\hat{Q}^{U,j,0}_{i,t}(s_0,a_i(0))\Bigl/\hat{Q}^{L,j,0}_{i,t}(s_0,a_i(0))-\hat{Q}^{U,j,1}_{i,t}(s_0,a_i(1))\Bigl/\hat{Q}^{L,j,1}_{i,t}(s_0,a_i(1))},~\text{if }\frac{\hat{Q}^{U,j,0}_{i,t}(s_0,a_i(0))}{\hat{Q}^{L,j,0}_{i,t}(s_0,a_i(0))}\neq\frac{\hat{Q}^{U,j,1}_{i,t}(s_0,a_i(1))}{\hat{Q}^{L,j,1}_{i,t}(s_0,a_i(1))},\\
        0,~\text{otherwise}.
    \end{cases}
\end{multline}
where $a_i(\varsigma) = \actionbarphi_{m,i}^\varsigma(s_0)$ with $\varsigma\in\{0,1\}$.
If $\La_{i,m}(s_{i,j}) =1$, then we do not need to update the index because $\La_{i,m}(s_{i,j}) =1$ is considered as the passive action with the least marginal cost and we do not assign any index for passive actions.
Implement the primary policy $\psi$ based on the updated $\hatnuVec_t$, and observe $\bs^{\psi}(t+1)$.

\item For a prescribed precision parameter $\epsilon >0$, if
\begin{equation}\label{eqn:stop}
\Bigl\lVert \hat{\bm{Q}}_{t}^{\Xi,j,\varsigma} -\hat{\bm{Q}}_{t-1}^{\Xi,j,\varsigma}\Bigr\rVert < \epsilon,
\end{equation}
for all $\Xi \in \{V,U,L\}$, $j\in[J]$ and $\varsigma\in\{0,1\}$, then take an element \[s^*_i\in\arg\max_{s\in\bS_i: \La_{i,m}(s)\geq 2}\hat\nu_{i,t}(s,\La_{i,m}(s)),\] update $\LaVec_{i,m+1}$ with $\Ia_i^{s^*_i}\LaVec_{i,m}$ for all $i\in[I]$, increment $m$ by $1$,  and go to Step~\ref{step:stop-condition}; otherwise, increment time stamp by one and go back to Step~\ref{step:q-iteration}.
\item \label{step:stop-condition}
If $m=\LenM_i$ for all $i\in[I]$, then
stop all the learning processes, set $\hatnuVec_{\tau} = \hatnuVec_t$ and $\bar{p}_{\tau} = 0$  for all $\tau \geq t$, and continue implementing the primary policy $\psi$; otherwise, increment time stamp by one and go back to Step~\ref{step:period}.\label{step:end}
\end{enumerate}

We refer to the primary policy $\psi$ described earlier in this subsection and the coordinated secondary policies $\bar{\phi}_m^{\varsigma}$ described in Steps~\ref{step:1}-\ref{step:end} as the \emph{Online MP Index (OMPI) Algorithm}.
We also provide the pseudo-codes for Steps~\ref{step:1}-\ref{step:end} in Appendix~\ref{app:ompi}.
This OMPI algorithm is a special case of the WCG-learning scheme with $K$ learning processes described in Section~\ref{subsec:q-index}.

Recall that, based on Proposition~\ref{prop:q-convergence}, for any precision parameter $\epsilon>0$, the estimated MP indices converge to the offline MP indices within a finite time horizon.
If all states have different MP indices (no tie case), then there exists sufficiently small $\epsilon >0$ such that, for any $i\in[I]$, the state ranking according to the descending order of the estimated $\hatnu_i(s,\La)$ coincide with that for the real MP indices within finite time horizon, and the long-run average reward also coincides with that of the offline MP index policy.
If the offline MP index policy is also asymptotically optimal (approaching optimality as $N_i\rightarrow \infty$ proportionately for all $i\in[I]$) , then the OMPI is also asymptotically optimal.
Moreover, from Theorem~\ref{theorem:convergence_Z_exp}, increasing the number of bandit processes (increasing $N_i$) will significantly stimulate the convergence of the estimated Q factors. We will discuss in details in Section~\ref{subsec:q-factor-based} about how fast the problem size, measured by the number of bandit processes $N_i$ ($i\in[I]$), can stimulate the convergence.

\subsection{Fast convergence in magnitude dimension}
\label{subsec:q-factor-based}


From Theorem~\ref{theorem:convergence_Z_exp} and Corollary~\ref{coro:converge_transition_prob}, when $h$ is large and an SA pair is visited at time $t$ with significant value of $Z^{\varphi,h}_{\iota}(t)$, it is likely that $z^{\varphi,h}_{\iota}(t) > 0$, and that the estimated transition probability $|\hat{\mathscr{N}}^{\varphi,h}_{i_{\iota},t}(s_{\iota},\actionphi_{\iota},s')|/|\hat{\mathscr{N}}^{\varphi,h}_{i_{\iota},t}(s_{\iota},\actionphi_{\iota})|$ is close to the real transition kernel $p_{i_{\iota}}(s_{\iota},\actionphi_{\iota},s')$ for all $s'\in\bS_{i_{\iota}}$.
We propose the following steps, which are added to the WCG-learning scheme with $K$ learning processes described in Section~\ref{subsec:q-index}, to stimulate the convergence of the $K$ Q learning processes. 
We prescribe a threshold parameter $\bar{\bm{\epsilon}}\in (0,1)^{\mathcal{I}}$.
\begin{enumerate}[label=\roman*)]\setcounter{enumi}{7}
\item\label{step:8} Before the learning processes start, for $i\in[I]$ and $a\in\bA_i$, initialize $\hat{\mathcal{P}}_i(a)=\bigl[\hat{p}_i(s,a,s')\bigr]_{|\bS_i|\times|\bS_i|}$ to be zero matrices and $\hat{r}^k_i(s,a) = 0$ ($i\in[I],(s,a)\in\bS_i\times \bA_i,k\in[K]$).
\item \label{step:9} For each time $t$, after observe $\hat{\mathscr{N}}^{\varphi,h}_{i_{\iota},t}(s_{\iota},\actionphi_{\iota},s_{\iota'})$ for all $\iota,\iota'\in[\mathcal{I}]$, if $Z^{\varphi,h}_{\iota}(t) > \bar{\epsilon}_{\iota}$ for some $\iota\in[\mathcal{I}]$, then update $\hat{p}_{i_{\iota}}(s_{\iota},\actionphi_{\iota},s_{\iota'})= |\hat{\mathscr{N}}^{\varphi,h}_{i_{\iota},t}(s_{\iota},\actionphi_{\iota},s')|/|\hat{\mathscr{N}}^{\varphi,h}_{i_{\iota},t}(s_{\iota},\actionphi_{\iota})|$ and $\hat{r}^k_{i_{\iota}}(s_{\iota},\actionphi_{\iota}) = \frac{1}{|\hat{\mathscr{N}}^{\varphi,h}_{i_{\iota},t}(s_{\iota},\actionphi_{\iota})|}\sum_{n\in\hat{\mathscr{N}}^{\varphi,h}_{i_{\iota},t}(s_{\iota},\actionphi_{\iota})}\bar{R}^{k}_{i_{\iota},n}(s_{\iota},\actionphi_{\iota})$, where $\bar{R}^{k}_{i_{\iota},n}(s_{\iota},\actionphi_{\iota})$ is the instantaneous reward rate for the $k$th learning process when it is in SA pair $\iota$ at time $t$.
Recall that we allow each of the $K$ learning processes to have different reward functions.
\item \label{step:10}Upon time $T$, if, for all $\iota\in[\mathcal{I}]$, there was at least a time slot $t\in[T]$ such that $Z^{\varphi,h}_{\iota} > \bar{\epsilon}_{\iota}$, then, denote such $T$ as $T^*$, stop updating $\hat{\mathcal{P}}_i(a)$ and $\hat{r}_i$, and go to Step~\ref{step:11};
otherwise, go back to Step~\ref{step:8}.
\item \label{step:11}
We update the estimated Q factors $\hat{\bm{Q}}_t=(\hat{\bm{Q}}^1_t;\hat{\bm{Q}}^2_t;\ldots;\hat{\bm{Q}}^K_t)$ through the value iteration method with respect to the operator $\mathcal{T}^{\phi}_i$ for all $i\in[T^*]$ defined in \eqref{eqn:define_H_phi} and the estimated transition kernels $\hat{\mathcal{P}}_i(a)$.
More precisely, for each $k\in[K]$, secondary policy $\bar{\phi}_k\in\oPhiLocal$, and the most updated estimated Q factors $\hat{Q}^k_{i,T^*}$ ($i\in[I]$), initialize $\tilde{Q}^k_{i,0} = \hat{Q}^k_{i,T^*}$.
For $\phi\in\oPhiLocal$, $i\in[I]$, $(s,a)\in\bS_i\times\bA_i$, and $Q\in\bF(\bS_i\times\bA_i)$, we define
\begin{equation}\label{eqn:define_hat_T}
(\hat{\mathcal{T}}^{k,\phi}_iQ)(s,a) = \hat{r}^k_i(s,a) + \sum_{s'\in\bS_i\backslash\{s_0\}}\hat{p}_i(s,a,s')Q(s',\actionphi_i(s')),
\end{equation} 
where recall $s_0$ is the ergodic state.
We keep updating  $\tilde{Q}^k_{i,m+1} = \hat{\mathcal{T}}^{k,\bar{\phi}_k}_i\tilde{Q}^k_{i,m}$ until $sp\bigl( \tilde{Q}^k_{i,m+1} - \tilde{Q}^k_{i,m}\bigr) < \epsilon$  for all $i\in[I]$, where $\epsilon>0$ is a hyper-parameter for precision purpose, $sp(\bm{v})\coloneqq\max_{n\in[N]}v_n - \min_{n\in[N]}v_n$ is the span seminorm of vector $\bm{v}=(v_n:n\in [N])\in\mathbb{R}^N$.
Let $\tilde{Q}^k_i$ represent the output of such a value iteration process.
We then update $\hat{Q}^k_{i,T^*} = \tilde{Q}^k_i$ for all $i\in[I]$ and $k\in[K]$ and stop the stimulate process.
\end{enumerate}

We refer to the above Steps~\ref{step:8}-\ref{step:11} as the \emph{stimulate process}. We can incorporate the stimulate process in to the WCG-learning scheme with the above-mentioned additional steps for all time slots $t\in[T^*]$.
We refer to the WCG-learning process with the incorporated stimulate process as the \emph{WCG-StimL} process.

For the $K$ secondary polices $\bar{\phi}_k\in\oPhiLocal$ ($k\in [K]$), let $\bar{\bm{\phi}}\coloneqq (\bar{\phi}_k: k\in[K])$.
For a probability simplex $\Delta_{[\mathcal{I}]}\coloneqq \{\bm{z}\in[0,1]^{\mathcal{I}}|\sum_{\iota\in[\mathcal{I}]}z_{\iota}=1\}$, $\bm{z}_0\in\Delta_{[\mathcal{I}]}$,  $\bm{Q}_0\in \prod_{i\in[I]}\bF(\bS_i\times\bA_i)$ and secondary polices $\bar{\bm{\phi}}\in(\oPhiLocal)^K$, define $\Phi^2(\bm{z}_0,\bm{Q}_0,\bar{\bm{\phi}})\subset \Phi^1$ as a subset of $\Phi^1$, where, for any $\varphi\in\Phi^2(\bm{z}_0,\bm{Q}_0,\bar{\bm{\phi}})$, given initial states $\bm{Z}^{\varphi,h}(0)=\bm{z}_0$ and $\hat{\bm{Q}}_0=\bm{Q}_0$, and employed secondary policies $\bar{\bm{\phi}}$, there exists $T<\infty$ such that, for all $h\in\mathbb{N}_+\cup\{\infty\}$ and SA pairs $\iota\in[\mathcal{I}]$,
\begin{equation}\label{eqn:Phi^1}
\max_{t=0,1,\ldots,T}\mathbb{P}\Bigl\{s^{\varphi,h}_{i_{\iota},n}(T) = s_{\iota},\actionvarphi^{h}_{i_{\iota},n}(T)=\actionphi_{\iota}~\Bigl|~\bm{Z}^{\varphi,h}(0)=\bm{z}_0,\hat{\bm{Q}}^{\varphi,h}_0=\bm{Q}_0\Bigr\}=
\max_{t=0,1,\ldots,T}z^{\varphi,h}_{\iota}(t) >0,
\end{equation}
where $n$ can be any element in $[N_{i_{\iota}}]$.
Equation~\eqref{eqn:Phi^1} indicates that, for such $\varphi\in\Phi^2(\bm{z}_0,\bm{Q}_0,\bar{\bm{\phi}})$, every SA pair has been visited and the associated Q factors $\hat{Q}^k_{i_{\iota},t}(s_{\iota},\actionphi_{\iota})$ ($k\in[K]$) have been updated by time $T^*$.
We have the following proposition based on \cite[Proposition 6.6.1]{puterman2005markov}.

\begin{proposition}\label{prop:stimulating_algo:convergence_span}
For $\bm{z}_0\in\Delta_{[\mathcal{I}]}$,  $\bm{Q}_0\in\prod_{i\in[I]}\bF(\bS_i\times\bA_i)$, secondary polices $\bar{\bm{\phi}}\in\bigl(\oPhiLocal\bigr)^K$, and primary policy $\varphi\in\Phi^2(\bm{z}_0,\bm{Q}_0,\bar{\bm{\phi}})$,
there exist $B,C>0$, and $H,M<\infty$ such that, for all $h>H$, $k\in[K]$, and $i\in[I]$, 
\begin{equation}\label{eqn:stimulating_algo:convergence_span}
\mathbb{P}\Bigl\{sp(\tilde{Q}^k_{i,M+1} - \tilde{Q}^k_{i,M}) > e^{-B} sp(\tilde{Q}^k_{i,1} - \tilde{Q}^k_{i,0}) \Bigr\} < e^{-Ch},
\end{equation}
where recall $span(\bm{v})$ is the span seminorm of vector $\bm{v}$, and the randomness of $\tilde{Q}^k_{i,m}$ ($m\in\mathbb{N}_+$) is led by the estimated $\hat{p}_i(s,a,s')$ and $\hat{r}^k_i(s,a)$ (both dependent on the scaling parameter $h$).
\end{proposition}
The proof of Proposition~\ref{prop:stimulating_algo:convergence_span} is provided in Appendix~\ref{app:stimulating_algo:convergence_span}.

Based on Proposition~\ref{prop:stimulating_algo:convergence_span}, there exists $H<\infty$ and $B,C>0$ such that, for all $h>H$, $k\in[K]$, $i\in[I]$ and $m\in\mathbb{N}_+$,
\begin{equation}\label{eqn:stimulating_algo:convergence_span:add}
\mathbb{P}\Bigl\{sp(\tilde{Q}^k_{i,m+1} - \tilde{Q}^k_{i,m}) > e^{-B m} sp(\tilde{Q}^k_{i,1} - \tilde{Q}^k_{i,0}) \Bigr\} < e^{-Ch}.    
\end{equation}
That is, when such $e^{-Ch}\downarrow 0$, the span seminorm $sp(\tilde{Q}^k_{i,m+1} - \tilde{Q}^k_{i,m})$  diminishes exponentially in $m$ almost surely, and
the value iteration described in Step~\ref{step:11} finishes quickly.

For $k\in[K]$, $i\in[I]$, and $\bar{\phi}_k\in\oPhiLocal$, let $Q^{k,\bar{\phi}_k}_i\coloneqq \lim_{m\rightarrow\infty} \bigl(\mathcal{T}^{k,\bar{\phi}_k}_i\bigr)^m Q^k_0$ for any initial value $Q^k_0\in\bF(\bS_i\times\bA_i)$.
Such $Q^{k,\bar{\phi}_k}_i$ is also the unique solution to $\mathcal{T}^{k,\bar{\phi}_k}_iQ = Q$ and the true Q factors of the real process under policy $\bar{\phi}_k$. The uniqueness is led by the finiteness of $\bS_i$ and the existence of the ergodic state $s_0$ (that is, with positive probability, $s_0$ is reachable from any other state in $\bS_i$ within a finite time period).


\begin{proposition}\label{prop:stimulating_algo:convergence_norm}
For $\bm{z}_0\in\Delta_{[\mathcal{I}]}$,  $\bm{Q}_0\in\prod_{i\in[I]}\bF(\bS_i\times\bA_i)$, secondary polices $\bar{\bm{\phi}}\in\bigl(\oPhiLocal\bigr)^K$, and primary policy $\varphi\in\Phi^2(\bm{z}_0,\bm{Q}_0,\bar{\bm{\phi}})$, there exist $C_1,C_2,H<\infty$ and $C_3 > 0$ such that, for all $h>H$, $k\in[K]$, and $i\in[I]$, and any $\epsilon>0$ used for the stopping condition of the value iteration in Step~\ref{step:11},
\begin{equation}\label{eqn:stimulating_algo:convergence_norm}
\mathbb{P}\Bigl\{\lVert \tilde{Q}^k_i - Q^{k,\bar{\phi}_k}_i\rVert > C_1\epsilon + C_2\lVert \bm{w}^{k,\bar{\phi}_k}_t(\tilde{Q}^k_i)\rVert\Bigr\} < e^{-C_3h},
\end{equation} 
where recall $\tilde{Q}^k_i$ is the output of the value iteration in \rm{Step~\ref{step:11}}, and $\bm{w}^{k,\bar{\phi}_k}_t(\cdot)$ is defined in \eqref{eqn:define_general_w}.
\end{proposition}
The proof of Proposition~\ref{prop:stimulating_algo:convergence_norm} is provided in Appendix~\ref{app:stimulating_algo:convergence_norm}.

Proposition~\ref{prop:stimulating_algo:convergence_norm} is proved based on the fast convergence,  as $h\rightarrow \infty$, of the estimated probabilities $\hat{p}_i(s,a,s')$ to the real ones $p_i(s,a,s')$. 
Together with \eqref{eqn:convergence_transition_prob:2} in Corollary~\ref{coro:converge_transition_prob}, we can further obtain, for the primary policy $\varphi\in\Phi^2(\bm{z}_0,\bm{Q}_0,\bar{\bm{\phi}})$, and any $\delta > 0$, 
\begin{equation}\label{eqn:stimulating_algo:convergence_Q}
\lim_{h\rightarrow \infty}\mathbb{P}\Bigl\{\bigl\lVert \tilde{Q}^k_i - Q^{k,\bar{\phi}_k}_i\bigr\rVert > \delta\Bigr\}  =0,
\end{equation}
by taking $0<\epsilon <\delta/2C_1$ used for the stopping condition in Step~\ref{step:11} and sufficiently large $h$ such that $\lVert \bm{w}^{k,\bar{\phi}_k}_t(\tilde{Q}^k_i)\rVert < \delta/2C_2$.
If we further assume \eqref{eqn:assumption:R_exp}, then, based on \eqref{eqn:convergence_transition_prob:3},
the convergence in \eqref{eqn:stimulating_algo:convergence_Q} happens exponentially in $h$; that is, there exist $C>0$ and $H <\infty$ such that for all $h>H$,
\begin{equation}\label{eqn:stimulating_algo:convergence_Q:exp}
\mathbb{P}\Bigl\{\lVert\tilde{Q}^k_i - Q^{k,\bar{\phi}_k}_i\rVert > \delta\} \leq e^{-Ch}.
\end{equation}

In an ideal case with sufficiently large $h$, within a finite time horizon $T<\infty$ (for which all SA pairs have been visited at least once), the learned Q factors $\hat{\bm{Q}}_T$ have been sufficiently close to the real ones so that the SA-pair ranking $\bm{\rank}_K(\hat{\bm{Q}}^1_T,\hat{\bm{Q}}^2_T,\ldots,\hat{\bm{Q}}^K_T)$ already coincides with that led by the real Q-factor-based indices, such as the Whittle/MP indices for the (multi-action) RMAB case.
Note that, for coinciding with the SA-pair ranking led by the real indices, the estimated Q factors $\hat{\bm{Q}}_T$ do not need to be exactly the same as the real ones.
In the ideal case, by such a time $T<\infty$, we can stop all the learning processes and exploit the Q-factor-based index policies (for example, the OMPI for the (multi-action) RMAB case), based on $\hat{\bm{Q}}_T$, which coincide with the offline index policies and will achieve the same long-run average performance.
If, further, the offline index policies approach optimality of the WCG problem as $h\rightarrow \infty$ (asymptotically optimal), such as the Whittle/MP index policy \cite{whittle1988restless,weber1990index,nino2007dynamic} and other index policies in non-RMAB case \cite{fu2020energy,fu2018restless,fu2025restless,fu2024patrolling}, then the Q-factor-based primary policy also approaches optimality asymptotically.

Proposition~\ref{prop:stimulating_algo:convergence_norm} and \eqref{eqn:stimulating_algo:convergence_Q} (or \eqref{eqn:stimulating_algo:convergence_Q:exp}) clarifies how  $h$ contribute to the improved accuracy of the learned Q factors, which is a different dimension than the timeline $t\in\mathbb{N}_0$. 
It is remarkable to reveal this new dimension, which we refer to as the \emph{magnitude dimension}.
It can lead to efficient convergence of the learning process(es) without disturbing exploitation of the primary algorithm in the timeline.
It hence follows with simultaneous control and learning.


\section{Optimality for Finite Time Horizon: Linear-programming-based approximations}
\label{sec:LP-algo}

We will focus on finite time horizon objectives in this subsection.
\cite{gast2024reoptimization,brown2023fluid} proposed the linear-programming-based approximation method that casts the same action variables as those for an optimal solution of the relaxed problem, solved through linear programming, and then adapts the actions a little bit to fit the constraints of the original problem.
Such action variables are applicable to the original problem; meanwhile, the resulting stochastic process is ensured to converge to optimality of the relaxed problem, which is an upper/lower bound of the maximum/minimum of the original problem, in the asymptotic regime.
It leads to asymptotic optimality of the proposed policy.

\subsection{Asymptotic Optimality with Fully Known Transition Kernels}
Here, we start with analyzing the asymptotic behaviors of the WCG problem with full knowledge of the transition kernels $\mathcal{P}_i(a)$ and reward functions $r_i$ ($i\in[I],a\in\bA_i$) in the finite time horizon case.


For $\phi\in\lPhi$, recall that the action vector $\actionphi^{h}(t)\in \sA$ is a function of the history $\bh_t$.
Similar to the discussion in Section~\ref{subsec:offline-problem}, we consider a random vector $\actiondbphi^{h}(t)\coloneqq (\actiondbphi^h_{i,n}(t):i\in[T],n\in[N_i])\in\db{\sA}$,  where recall $\db{\sA}$ is the space of the random actions.
The probability of taking $\actiondbphi^{h}(t)=\actionphi$ for $\actionphi\in\sA$ is a function of the history $\bh_t$.
We define $\dblPhi$ as the set of all the policies $\dbphi$ determined by such random actions $\actiondbphi^{h}(t)$ for all $t\in[T]_0$.
Note that, when we have full knowledge of the transition kernels and reward rates and do not rely on any estimated Q factors $\hat{\bm{Q}}_t$, we can consider $\hat{\bm{Q}}_t\equiv \bm{0}$, and the policies in $\dblPhi$ are equivalent to those dependent on $\bh_t$ for $t\in[T]_0$.

We consider the problem, for $T<\infty$,
\begin{equation}\label{eqn:objective:finite-time:h}
    \max_{\dbphi\in\dblPhi} \frac{1}{h}\Gamma^{\dbphi,h}(T,\bs_0),
\end{equation}
subject to \eqref{eqn:constraint:linear:h} (with $\actionphi^{h}_{i,n}(t)$ replaced by $\actiondbphi^{h}_{i,n}(t)$).
Along the same lines as the relaxation from \eqref{eqn:constraint:linear} to \eqref{eqn:constraint:long-run:relax} described in Section~\ref{subsec:offline-problem}, 
we relax \eqref{eqn:constraint:linear:h} to
\begin{equation}\label{eqn:constraint:relax:h}
\sum_{i\in[I]}\sum_{n\in[N_i]} \mathbb{E}\Bigl[f_{i,\ell}\bigl(s^{\dbphi,h}_{i,n}(t),\actiondbphi^{h}_{i,n}(t)\bigr)\Bigr]=0,~\forall \ell\in[L], t\in[T]_0.
\end{equation}

The problem \eqref{eqn:objective:finite-time:h}
subject to \eqref{eqn:constraint:relax:h} is a relaxed version of the problem described in \eqref{eqn:objective:finite-time:h} and \eqref{eqn:constraint:linear:h}.
Unlike the Lagrange relaxation used in Section~\ref{subsec:offline-problem}, we rewrite the relaxed problem described in \eqref{eqn:objective:finite-time:h} and \eqref{eqn:constraint:relax:h} as a linear programming form, 
\begin{equation}\label{eqn:obj:linear-programming}
\sup_{\bm{x}\in[0,1]^{\mathcal{I}T}}\sum_{t\in[T]_0}\sum_{\iota\in[\mathcal{I}]}N^0_{i_{\iota}} r_{i_{\iota}}(s_{\iota},\actionphi_{\iota})x_{\iota,t},
\end{equation}
subject to
\begin{eqnarray}
\label{eqn:constraint:linear-programming:1}&
\bm{x}^T_t\mathcal{P} = \bm{x}_{t+1}^T\tilde{\mathbb{I}},&\forall t\in[T-1]_0\\
\label{eqn:constraint:linear-programming:2}&\sum_{\iota\in [\mathcal{I}]:i_{\iota}=i}x_{\iota,t} = 1, &~\forall i\in[I],t\in[T]_0,\\
&\sum_{\iota\in[\mathcal{I}]}N_{i_{\iota}}^0 x_{\iota,t}f_{i_{\iota},\ell}(s_{\iota},\actionphi_{\iota}) =0,&~\forall \ell\in[L],t\in[T]_0,
\label{eqn:constraint:linear-programming:3}
\end{eqnarray}
where $\bm{x}^T$ is the transpose of $\bm{x}$, $\bm{x}_t\coloneqq (x_{\iota,t}:\iota\in[\mathcal{I}])$, initial state $\bm{x}_0$ is given, $\mathcal{P}\coloneqq \bigl[p(\iota,i',s')\bigr]_{\mathcal{I}\times\sum_{i\in[I]}|\bS_i|}$ with
\begin{equation}\label{eqn:linear-programming:1}
    p(\iota,i',s') \coloneqq \begin{cases}
        p_{i_{\iota}}(s_{\iota},\actionphi_{\iota},s'),&\text{if }i' = i_{\iota},\\
        0,&\text{otherwise},
    \end{cases}
\end{equation}
and $\tilde{\mathbb{I}}\coloneqq \bigl[\tilde{\scri}(\iota,i',s')\bigr]_{\mathcal{I}\times\sum_{i\in[I]}|\bS_i|}$ with
\begin{equation}\label{eqn:linear-programming:2}
\tilde{\scri}(\iota,i',s') \coloneqq \begin{cases}
    1, &\text{if }i'=i_{\iota},s'=s_{\iota},\\
    0, &\text{otherwise}.
\end{cases}
\end{equation}
In \eqref{eqn:constraint:linear-programming:1}, $x_{\iota,t}$ ($\iota\in[\mathcal{I}]$) represents the probability that the process $\{s^{\dbphi,h}_{i,n}(t),t\in\mathbb{N}_0\}$, for any $n\in[N_{i_{\iota}}]$, stays in SA pair $\iota$.
Recall that we are interested in the case where solutions to \eqref{eqn:obj:linear-programming}-\eqref{eqn:linear-programming:3} exist.
Given $\bm{x}^*$ optimal to the linear programming problem \eqref{eqn:obj:linear-programming}-\eqref{eqn:constraint:linear-programming:3}, we can construct a corresponding policy $\dbphi^*\in\dblPhi$ by setting 
\begin{equation}\label{eqn:linear-programming:3}
    \alpha^{\dbphi^*}_{\iota}(t)\coloneqq\mathbb{P}\Bigl\{\actiondbphi^{*,h}_{i_{\iota},n}(t) = \actionphi_{\iota}\Bigl|s^{\dbphi^*,h}_{i_{\iota},n}(t) = s_{\iota}\Bigr\}=\frac{x^*_{\iota,t}}{\sum_{\iota'\in[\mathcal{I}]:i_{\iota}=i_{\iota'},s_{\iota}=s_{\iota'}}x^*_{\iota',t}},~\forall \iota\in[\mathcal{I}],n\in[N_{i_{\iota}}],t\in[T]_0.
\end{equation}
Such $x^*_{\iota,t}$ is the probability of a bandit process $\{s^{\dbphi^*,h}_{i_{\iota},n}(t),t\in[T]_0\}$ (for any $n\in[N_{i_{\iota}}]$) staying in SA pair $\iota$ at time $t$ under the policy $\dbphi^*$ that is optimal to the relaxed problem \eqref{eqn:objective:finite-time:h} and \eqref{eqn:constraint:relax:h}.
We also have $z^{\dbphi^*,h}_{\iota}(t) =\frac{N^0_{i_{\iota}}}{\sum_{i\in[I]}N^0_i}x^*_{\iota,t}$ (for all $\iota\in[\mathcal{I}]$, $t\in[T]_0$, and $h\in\mathbb{N}_+$) and $\bm{x}^*_t$ 
as $x^*_{\iota,t}  = \frac{\sum_{i\in[I]}N_i^0 }{N_{i_{\iota}}^0} z^{\dbphi^*,h}_{\iota}(t)$.
Define conversion matrices $\mathcal{Z},\mathcal{X}\in\mathbb{R}^{\mathcal{I}\times\mathcal{I}}$ for which $\bm{z}^{\dbphi^*,h}(t) = \mathcal{Z}\bm{x}_0^*$ and $\bm{x}^*_t = \mathcal{X}\bm{z}^{\dbphi^*,h}(t)$.
Apparently, such $\dbphi^*$ is not necessarily applicable to the original problem.
We define a set $\dblPhiz\subset \dblPhi$ of policies $\dbphi$ determined by action probabilities 
$\alpha^{\dbphi}_{\iota}(t)\coloneqq \mathbb{P}\bigl\{\actiondbphi^{h}_{i_{\iota},n}(t) = \actionphi_{\iota}\bigl|s^{\dbphi,h}_{i_{\iota},n}=s_{\iota}\bigr\}$ (any $n\in[N_{i_{\iota}}]$)
for all $\iota\in[\mathcal{I}]$ and $t\in[T]_0$, for which $\bm{\alpha}^{\dbphi}_{i}(s,t) \coloneqq (\alpha^{\dbphi}_{\iota}(t): \iota\in[\mathcal{I}],i_{\iota}=i,s_{\iota}=s)$ is a function $\bS_i\times \bigl([T]_0\bigr)\rightarrow \Delta_{\bA_i}$.
Let $\bm{\alpha}^{\dbphi}(t) \coloneqq (\alpha^{\dbphi}_{\iota}(t) : \iota\in[\mathcal{I}])$.


For $\dbphi\in\dblPhi$, define $\bm{\Upsilon}^{\dbphi,h}(t)\coloneqq (\Upsilon^{\dbphi,h}_{i,s}(t):i\in[I],s\in\bS_i)$, where
\begin{equation}
    \Upsilon^{\dbphi,h}_{i,s}(t) \coloneqq \sum_{\begin{subarray}~\iota\in[\mathcal{I}]:\\i_{\iota}=i,\\s_{\iota}=s\end{subarray}}Z^{\dbphi,h}_{\iota}(t)=\mathcal{U}\bm{Z}^{\dbphi,h}(t)),
\end{equation}
for a matrix $\mathcal{U}\in \mathbb{R}^{(\sum_{i\in[I]}|\bS_i|)\times\mathcal{I}}$.
Let $\bm{\upsilon}^{\dbphi,h}(t) \coloneqq \mathbb{E}\bm{\Upsilon}^{\dbphi,h}(t)$.
Also, for a policy $\dbphi\in\dblPhi$,
define the proportion of bandit processes that are in SA pair $\iota\in[\mathcal{I}]$ at time $t$ under $\phi$ as
\begin{equation}
\alpha^{\dbphi,h}_{\iota}(t) \coloneqq 
\begin{cases}
    \frac{Z^{\dbphi,h}_{\iota}(t)}{\Upsilon^{\dbphi,h}_{i_{\iota},s_{\iota}}(t)}, &\text{if } \Upsilon^{\dbphi,h}_{i_{\iota},s_{\iota}}(t) > 0,\\
    \frac{1}{|\bA_{i_{\iota}}|},&\text{otherwise},
\end{cases}
\end{equation}
with $\bm{\alpha}^{\dbphi,h}(t)\coloneqq(\alpha^{\dbphi,h}_{\iota}(t):\iota\in[\mathcal{I}])$.
In general, for $\dbphi\in\dblPhi$, such $\bm{\alpha}^{\dbphi,h}(t)$ is a random vector dependent on the history $\bh(t)$.
In the special case with $\dbphi\in\dblPhiz$,  $\bm{\alpha}^{\dbphi}(t)\coloneqq\bbE\bm{\alpha}^{\dbphi,h}(t)$, which is independent from $h$. 
We have the following proposition.

\begin{proposition}\label{prop:asym_opt:LP}
For any $\dbpsi\in\dblPhi$ and $\dbphi\in\dblPhiz$, $\bm{Z}^{\dbpsi,h}(0)=\bm{Z}^{\dbphi,h}(0)$, if 
there exists constant $K < \infty$ such that 
\begin{equation}\label{eqn:assumption:action}
\lim_{h\rightarrow \infty} \bbP\Bigl\{\bigl\lVert \bm{\alpha}^{\dbpsi,h}(t)- \bm{\alpha}^{\dbphi}(t)\bigr\rVert > K\bigl\lVert \bm{\Upsilon}^{\dbpsi,h}(t)-\bm{\upsilon}^{\dbphi,h}(t)\bigr\rVert\Bigr\} =0,~\forall t\in[T]_0,
\end{equation}
then, for any $\epsilon>0$ and $T<\infty$, there exist $C_1,H<\infty$ and $C_2>0$ such that, for all $h>H$,
\begin{equation}\label{eqn:prop:asym_opt:LP:1}
 \max_{t\in[T]_0} \bigl\lVert \bm{z}^{\dbpsi,h}(t)-\bm{z}^{\dbphi,h}(t)\rVert \leq C_1e^{-C_2h} + \epsilon,
\end{equation}
where recall $\bm{z}^{\dbpsi,h}(t) = \mathbb{E}[\bm{Z}^{\dbpsi,h}(t)]$ and $\bm{z}^{\dbphi,h}(t) = \mathbb{E}[\bm{Z}^{\dbphi,h}(t)]$. 
\end{proposition}
The proof of Proposition~\ref{prop:asym_opt:LP} is based on Theorem~\ref{theorem:convergence_Z_exp} and is provided in Appendix~\ref{app:asym_opt:LP}. 
For any policy $\dbphi\in\dblPhi$, the expected total reward $\frac{1}{h}\Gamma^{\dbphi,h}(T,\bs_0) = \sum_{t\in[T]_0}\bm{r}\cdot \bm{z}^{\dbphi,h}(t)$, where $\bs^{\phi,h}(0)=\bs_0$ is the initial state, and
recall $\bm{r} = (r_{i_{\iota}}(s_{\iota},\actionphi_{\iota}): \iota\in[\mathcal{I}])$.
Based on Proposition~\ref{prop:asym_opt:LP}, for $\dbphi=\dbphi^*$, if we propose a policy $\dbphi\in\dblPhi$ such that \eqref{eqn:assumption:action} and the original constraints~\eqref{eqn:constraint:linear:h} are both satisfied ($\dbpsi$ is applicable to the original WCG problem), then $\dbpsi$ is asymptotically optimal with performance deviation exponentially diminishing in $h$.
Because $\dbphi^*$ is optimal to the relaxed problem, which achieves an upper bound of the maximum of the original WCG problem described in \eqref{eqn:obj:h} and \eqref{eqn:objective:finite-time:h}.


For such $\dbpsi$, condition \eqref{eqn:assumption:action} requests convergence between $\bm{\alpha}^{\dbpsi,h}(t)$ and $\bm{\alpha}^{\dbphi^*,h}(t)$ and (Lipschitz) continuity of $\bm{\alpha}^{\dbpsi,h}(t)$ in $\bm{\Upsilon}^{\dbpsi,h}(t)$ around the point $ \bm{\upsilon}^{\dbphi^*,h}(t) = \mathcal{U}\mathcal{Z}\bm{x}^*$ in the asymptotic regime $h\rightarrow \infty$. 
Examples for such $\dbpsi$ that satisfy both the original constraints \eqref{eqn:constraint:linear:h} and \eqref{eqn:assumption:action} with $\dbphi=\dbphi^*$ include the \emph{Lagrangian index policy}~\cite{brown2020index} and the \emph{water filling policy}~\cite{gast2023linear} for the RMAB case, and the \emph{fluid-budget balancing policy}~\cite{brown2023fluid} which was defined for the special case assuming at most one saturate constraint in \eqref{eqn:constraint:linear:h}. 
We refer to those $\dbpsi$ satisfying \eqref{eqn:constraint:linear:h} and \eqref{eqn:assumption:action} as the \emph{adapted-linear-programming} (ALP) algorithms.

Proposition~\ref{prop:asym_opt:LP} proves that ALP approaches optimality in the speed of $O(e^{-h})$.
When the WCG reduces to a standard RMAB problem, assuming $\bm{x}^*$ is a non-degenerate solution to the linear programming problem, \cite{gast2023linear} provided similar results with exponential convergence between a linear-programming-based policy and optimality in the asymptotic regime.
For the general case, under a slightly different non-degenerate assumption that requests at most one constraint in \eqref{eqn:constraint:linear:h} is saturate (achieves equality), 
\cite{brown2023fluid} proved that their policy is asymptotically optimal with performance sub-optimality $O(e^{-h})$. 
Here, Proposition~\ref{prop:asym_opt:LP} has no request on the form of $\bm{x}^*$ (not necessary to be non-degenerate) and is applicable for general WCG and every ALP satisfying \eqref{eqn:assumption:action}.

 

\subsection{Asymptotically Optimal Control with Unknown Transition Kernels}
We then consider the WCG-StimL process where $\mathcal{P}$ and $r_i$ are not known a priori.

Define $\hat{\mathcal{P}}\coloneqq [\hat{p}_{_{\iota}}(s_{\iota},\actionphi_{\iota},s')]_{\mathcal{I}\times\sum_{i\in[I]}|\bS_i|}$ where $\hat{p}_{i_{\iota}}(s_{\iota},\actionphi_{\iota},s')$ is the estimated transition probabilities updated in Step~\ref{step:9}.
Similar to the Q-factor-based algorithms in Section~\ref{subsec:q-factor-based}, 
we propose in the following the steps of approximating the optimum of the linear programming problem described in \eqref{eqn:linear-programming:1}-\eqref{eqn:linear-programming:3}, which, more importantly, leads to asymptotically optimal ALP algorithms with suboptimality diminishing in $O(e^{-h})$.

Recall the stimulate process described in Steps~\ref{step:8}-\ref{step:11} in Section~\ref{subsec:q-factor-based}, where the transition matrix and reward function for the $k$th ($k\in[K]$) learning process are estimated through $\hat{\mathcal{P}}$ and $\hat{\bm{r}}^k$, respectively.
Here, to approximate the ALP algorithm, we need only $K=1$ learning process for the unknown $\mathcal{P}$ and $\bm{r}$.
In this context, since $k=K=1$, we omit the superscript $k$ for all the learning processes in this subsection.

Define a set of probability matrices 
\begin{equation}
    \mathscr{P} \coloneqq \Bigl\{\mathcal{P}\in [0,1]^{\mathcal{I}\times (\sum_{i\in[I]}|\bS_i|)}~\Bigr|~\forall \iota\in[\mathcal{I}],~\sum_{s'\in\bS_{i'}}(\mathcal{P})_{\iota,s'} = 1\Bigr\}.
\end{equation}
For prescribed $\varepsilon>0$, define another linear programming problem as follows.
\begin{equation}\label{eqn:obj:linear-programming:optimistic}
\sup_{\begin{subarray}~\bm{x}\in[0,1]^{\mathcal{I}T},\\\bar{\mathcal{P}}\in\mathscr{P},\\\bar{\bm{r}}\in\mathbb{R}^{\mathcal{I}}\end{subarray}}\sum_{t\in[T]_0}\sum_{\iota\in[\mathcal{I}]}N^0_{i_{\iota}} \bar{r}_{i_{\iota}}(s_{\iota},\actionphi_{\iota})x_{\iota,t},
\end{equation}
subject to
\begin{eqnarray}
\label{eqn:constraint:linear-programming:optimistic:1}&\bm{x}^T_t\bar{\mathcal{P}} = \bm{x}_{t+1}^T\tilde{\mathbb{I}},&\forall t\in[T-1]_0\\
\label{eqn:constraint:linear-programming:optimistic:2}&\sum_{\iota\in [\mathcal{I}]:i_{\iota}=i}x_{\iota,t} = 1, &~\forall i\in[I],t\in[T]_0,\\
&\sum_{\iota\in[\mathcal{I}]}N_{i_{\iota}}^0 x_{\iota,t}f_{i_{\iota},\ell}(s_{\iota},\actionphi_{\iota}) =0,&~\forall \ell\in[L],t\in[T]_0,
\label{eqn:constraint:linear-programming:optimistic:3}
\\
 &\hat{p}_{i_{\iota}}(s_{\iota},\actionphi_{\iota},s')- \frac{\varepsilon}{\mathcal{I}}\leq \bar{p}_{i_{\iota}}(s_{\iota},\actionphi_{\iota},s')\leq \hat{p}_{i_{\iota}}(s_{\iota},\actionphi_{\iota},s') + \frac{\varepsilon}{\mathcal{I}},& ~\forall \iota\in[\mathcal{I}],s'\in\bS_{i_{\iota}},\\
 &\hat{r}_{i_{\iota}}(s_{\iota},\actionphi_{\iota}) - \frac{\varepsilon}{\mathcal{I}}\leq \bar{r}_{i_{\iota}}(s_{\iota},\actionphi_{\iota}) \leq \hat{r}_{i_{\iota}}(s_{\iota},\actionphi_{\iota}) + \frac{\varepsilon}{\mathcal{I}}, &~\forall \iota\in [\mathcal{I}],\label{eqn:constraint:linear-programming:optimistic:5}
\end{eqnarray}
where $\bm{x}_0\in[0,1]^{\mathcal{I}}$ is given and satisfies $\sum_{i\in[\mathcal{I}]:i_{\iota}=i}x_{i,0} = 1$ for all $i\in[I]$.
We refer to the problem described in \eqref{eqn:obj:linear-programming:optimistic}-\eqref{eqn:constraint:linear-programming:optimistic:5} with prescribed  $\varepsilon>0$ and given $\bm{x}_0\in[0,1]^{\mathcal{I}}$ as the $(\varepsilon,\bm{x}_0)$-linear programming ($(\varepsilon,\bm{x}_0)$-LP) problem.
Let $\bm{x}^*(\varepsilon,\bm{x}_0)$ represent an optimal solution of the $(\varepsilon,\bm{x}_0)$-LP problem, and $\dbphi^*(\varepsilon,\bm{x}_0)\in\dblPhiz$ is the optimal policy led by $\bm{x}^*(\varepsilon,\bm{x}_0)$; that is,
\begin{equation}
   \alpha^{\dbphi^*(\varepsilon,\bm{x}_0)}_{\iota}(t) = \frac{x_{\iota,t}^*(\varepsilon,\bm{x}_0)}{\sum_{\iota'\in[\mathcal{I}]:i_{\iota}=i_{\iota'},s_{\iota}=s_{\iota'}}x^*_{\iota',t}(\varepsilon,\bm{x}_0)}.
\end{equation}
Now we consider a policy $\dbpsi(\varepsilon,\bm{x}_0)\in\dblPhi$ that is applicable to the original problem described in \eqref{eqn:objective:finite-time:h} and \eqref{eqn:constraint:linear:h} and satisfies \eqref{eqn:assumption:action} for $\dbpsi = \dbpsi(\varepsilon,\bm{x}_0)$ and $\dbphi = \dbphi^*(\varepsilon,\bm{x}_0)$. We have the following proposition shows that, for the WCG-StimL process, the convergence in transition kernels leads to tight convergence in the optimum of the linear programming problem described in \eqref{eqn:linear-programming:1}-\eqref{eqn:linear-programming:3} and $(\varepsilon,\bm{x}_0)$-LP.

For $\bm{z}_0\in\Delta_{[\mathcal{I}]}$, $\bm{Q}_0\in\prod_{i\in[I]}\bF(\bS_i\times\bA_i)$, and any secondary policy $\bar{\phi}\in\oPhiLocal$, consider a primary policy $\varphi\in\Phi^2(\bm{z}_0,\bm{Q}_0,\bar{\phi})$ that satisfies \eqref{eqn:Phi^1}.
Here, $\bm{z}_0$ and $\bm{Q}_0$ are the initial value of $\bm{Z}^{\varphi,h}(0)$ and $\hat{\bm{Q}}_0$ for any primary policy $\varphi$.
In this subsection, unlike the Q-factor-based algorithms, we do not need to specify the secondary policy $\bar{\phi}$ or the estimated Q factor $\hat{\bm{Q}}_t$, but only request the estimated transition matrix $\hat{\mathcal{P}}$ and reward function $\hat{\bm{r}}$ that are updated in Step~\ref{step:9} of the stimulate process.

We construct a policy $\dbvtheta\in\dblPhi$ as follows.
Starting from time $t=0$, we employ the primary policy $\varphi\in\Phi^2(\bm{z}_0,\bm{Q}_0,\bar{\phi})$ until time $T^*<\infty$ for which the stopping condition in Step~\ref{step:10} of the stimulate process has been reached. 
For time $t=T^*,T^*+1,\ldots$, we employ the policy $\dbpsi(\varepsilon,\mathcal{X}\bm{z}^{\dbvtheta,h}(T^*))$ that is applicable to the original problem and satisfies \eqref{eqn:assumption:action} for $\dbpsi=\dbpsi(\varepsilon,\mathcal{X}\bm{z}^{\dbvtheta,h}(T^*))$ and $\dbphi = \dbphi^*(\varepsilon,\mathcal{X}\bm{z}^{\dbvtheta,h}(T^*))$. 
We refer to such a policy $\dbvtheta$, comprising of $\varphi$ for $t\in[T^*-1]_0$ and $\dbpsi(\varepsilon,\mathcal{X}\bm{z}^{\dbvtheta,h}(T^*))$ on and after $T^*$, as the \emph{online adapted linear-programming} (OALP) algorithm.
The expected total reward of policy $\dbvtheta$ is then given by $\Gamma^{\dbvtheta}(T^*+T,\bs_0) = \Gamma^{\varphi}(T^*,\bs_0)+\Gamma^{\dbpsi(\varepsilon,\bm{x})}(T,\bs^{\dbvtheta,h}(T^*))$, where $\bm{x}=\mathcal{X}\bm{z}^{\dbvtheta,h}(T^*)$ is the vector of proportions of bandit processes in all the SA pairs that can be interpreted from $\bs^{\dbvtheta,h}(T^*)$ through the relationship defined in \eqref{eqn:define_Z}.


\begin{proposition}\label{prop:converge_LP_opt}
For any $\epsilon > 0$ and $T < \infty$, 
there exist $C_1, 
 H <\infty$, $0<\varepsilon < \epsilon$, and $C_2>0$ such that, for all $h>H$, 
\begin{equation}\label{eqn:prop:convergence_LP_opt:1}
    \max_{t=T^*,T^*+1,\ldots,T^*+T}\bigl\lVert \bm{z}^{\dbvtheta,h}(t) - \bm{x}^*_{t-T^*}(\varepsilon,\bm{x}_0)\bigr\rVert\leq C_1e^{-C_2 h} +\epsilon,
\end{equation}
where $\bm{x}_0 = \mathcal{X}\bm{z}^{\dbvtheta,h}(T^*)$, 
\end{proposition}
The proof of Proposition~\ref{prop:converge_LP_opt} is provided in Appendix~\ref{app:prop:convergence_LP_opt}.

From Proposition~\ref{prop:converge_LP_opt},  $\dbvtheta$ approaches $\dbphi^*(\varepsilon,\mathcal{X}\bm{z}^{\dbvtheta,h}(T^*))$ as $h\rightarrow\infty$ at a speed of $O(e^{-h})$.
The stopping time $T^*$ of the stimulate process is achieved when all the SA pairs $\iota\in[\mathcal{I}]$ has been explored (see Step~\ref{step:10}), which is usually small when the scaling parameter $h$ is sufficiently large. 
For example, if the primary policy $\varphi$ for $t=0,1,\ldots,T^*$ leads to a fully connected Markov chain for all the SA pairs, then $T^*=1$ for sufficiently large $h$.
Unlike past work on learning-based algorithms, such $T^*$ of OALP is independent from the exploitation horizon $T$,
because the convergence of the estimated transition kernels to real ones can be achieved through only the magnitude dimension without requesting further moving on in the timeline.
When selecting $T^* \ll T$, the negative effects of learning the unknown parameters are mitigated.

Recall that the WCG process $\{\bs^{\dbvtheta,h}(t), t\in[T+T^*]_0\}$ can be interpreted to $\{\bm{Z}^{\dbvtheta,h}(t),t\in[T+T^*]_0\}$.
The relationship between $\bs^{\phi,h}(t)$ and $\bm{Z}^{\phi,h}(t)$ is defined in~\eqref{eqn:define_Z}.
Based on Corollary~\ref{coro:converge_transition_prob} and Proposition~\ref{prop:converge_LP_opt},
we have the following corollary.
\begin{corollary}\label{coro:asym_opt:OALP}
For any $\epsilon > 0$ and $T<\infty$, there exists $0<\varepsilon < \epsilon$, $C_1,H<\infty$, $0<\varepsilon<\epsilon$, and $C_2>0$ such that, for all $h>H$,
\begin{equation}\label{eqn:coro:asym_opt:OALP}
\frac{1}{h}\Bigl\lvert \Gamma^{\dbphi^*,h}(T,\bs^{\dbvtheta,h}(T^*))-\Gamma^{\dbpsi(\varepsilon,\bm{x}),h}(T,\bs^{\dbvtheta,h}(T^*))\Bigr\lvert \leq C_1e^{-C_2 h}+\epsilon,
\end{equation}
where $\dbphi^*$ is optimal to the linear programming problem \eqref{eqn:obj:linear-programming}-\eqref{eqn:constraint:linear-programming:3} with real transition matrix $\mathcal{P}$, reward function $\bm{r}$ and the same initial condition $\bm{x}_0=\mathcal{X}\bm{z}^{\dbvtheta,h}(T^*)$.
\end{corollary}
The proof of Corollary~\ref{coro:asym_opt:OALP} is provided in Appendix~\ref{app:coro:asym_opt:OALP}.

Corollary~\ref{coro:asym_opt:OALP} indicates the convergence between the optimum of the offline relaxed problem, led by $\dbphi^*$, and the $\dbvtheta$ policy after $t=T^*$ and onward, when $h\rightarrow \infty$.
It leads to quickly diminishing sub-optimality, at the speed of $O(e^{-h})$, of $\dbvtheta$ to the original WCG problem, because the optimum of the relaxed problem is an upper bound for that of the original problem;
that is, OALP is asymptotically optimal. 
Recall that the expected total reward of policy $\dbvtheta$ is $\Gamma^{\dbvtheta}(T^*+T,\bs_0) = \Gamma^{\varphi}(T^*,\bs_0)+\Gamma^{\dbpsi(\varepsilon,\bm{x})}(T,\bs^{\dbvtheta,h}(T^*))$, where $\bm{x}=\mathcal{X}\bm{z}^{\dbvtheta,h}(T^*)$, and the stopping time of the exploration period $T^*$ is independent from the time length of the exploitation period $T$.
Such $T^*$ is usually small, since it requests to explore each SA pair for only one time.

The asymptotic optimality in the magnitude dimension is marvelous for large-scale WCG-StimL systems where the probem size, measured by the scaling parameter $h$, is large. 
If the system is small, then conventional learning and optimization techniques, such Q learning and neural network, can achieve solutions without requesting excessive large amount of computational and storage resources. 
When the system becomes large, conventional techniques exhibit the curse of dimensionality with respect to learning, control, and, more importantly, simultaneous learning and control.
Optimal solutions become intractable and we resort good approximations of optimality.
The WCG-StimL model and the (O)MPI and (O)ALP are applicable to a wide range of application scenarios, including the standard RMAB case~\cite{whittle1988restless,weber1990index,nino2001restless,verloop2016asymptotically,fu2019towards,fu2020energy,brown2020index} and \cite{gast2023linear} and extended restless-bandit-based problems~\cite{fu2020resource,fu2025restless} and \cite{fu2024patrolling}.

\section{Conclusions}\label{sec:conclusions}

We have proposed the WCG-StimL scheme that enables simultaneous learning and control.
Such a WCG-StimL process can simultaneously estimate multiple Q factors, for which the implemented policies are user-defined and can be different from the primary policy employed by the real WCG process.
We have proved that the estimated Q factors converge to the real ones in both the time and the magnitude dimensions. 
Impressively, we have proved that the estimated Q factors, as well as the estimated transition kernels (and the reward functions), converge to the real ones exponentially fast in the magnitude dimension.
In such a scheme, we can propose online version of policies for the WCG problem that are based on some Q factors and/or fully knowledge of transition kernels and reward functions. 

For the long-run-objective case, we have proposed the OMPI algorithm for the (multi-gear) RMAB problem, which converges to the offline Whittle index policy, when the system is PCL-indexable, exponentially fast in the magnitude dimension. 
If the Whittle index policy further approaches offline optimality in the magnitude dimension (that is, asymptotically optimal as discussed in many cases \cite{weber1990index,fu2016asymptotic,fu2020energy,fu2018restless}), then the OMPI also approaches offline optimality as the magnitude dimension $h\to\infty$.

For the finite time horizon case, we have proposed a class of LP-based offline policies, referred to as the ALP policies, that approach optimality with performance degradation diminishes exponentially in the magnitude domain $h$. 
Unlike the past work assuming non-degenerate conditions~\cite{gast2023linear,brown2023fluid}, our proofs for the exponentially diminishing sub-optimality do not request any specific form of the optimal solutions for the associated LP problems.
More importantly, we have proposed an online version of the ALP policies, referred to as OALP, and have proved that their performance degradations (against offline optimality) also diminish exponentially in the magnitude dimension.

\appendices
\input{proof_propn1}

\section{Proof of Lemma~\ref{lemma:exist_Q_star}}\label{app:lemma:exist_Q_star}
\proof{Proof of Lemma~\ref{lemma:exist_Q_star}.}
First we observe that it is enough to prove this result for each $i\in[I]$ separately, accordingly we fix $i$. Since $\mathcal{T}^{\phi}_i$ is affine, 
 it is enough to prove that the linear part, that is, 
 \begin{equation}\label{eqn:linearpart}
 \mathcal{T}^{\phi}_{i,0} Q(s,a) =  \sum_{s'\in\bS_i\backslash\{s_0\}}p_i(s,a,s')Q(s',\actionphi_i(s'))
 \end{equation}
 is a contraction for an appropriate norm on $\bF(\bS_i\times \bA_i)$.

 For  the MDP with transition probabilities 
  $p_i(s,a,s')$,  $(s,a,s')\in \bS_i\times \bA_i\times \bS_i$, 
let $\hat{J}(s)$ be the reward-to-go (value function) when such an MDP (under the policy $\phi$) starts from state $s$ until it goes into the ergodic state $s_0\in\bS_i$.  This satisfies 
\begin{equation}
\hat{J}(s) = -1 +\sum_{s'\in\bS_i\backslash\{s_0\}}p_i(s,\actionphi_i(s),s')\hat{J}(s'), \quad (s\in \bS_i). 
\end{equation} 
Define $\xi(s) \coloneqq -\hat{J}(s)$ ($s\in \bS_i$). Observe that $\xi(s)> 0$, and that 
\begin{equation}
    \sum_{s'\in\bS_i\backslash\{s_0\}}p_i(s,\actionphi_i(s),s') \xi(s') = -1 - \hat{J}(s)  = \xi(s)-1 = \frac{\xi(s)-1}{\xi(s)}\xi(s)\leq \beta \xi(s).
\end{equation}
where  $\beta = \max_{s'\in\bS_i}\frac{\xi_i(s)-1}{\xi_i(s)} < 1$.
Now, we define the appropriate norm on $\bF(\bS_i\times \bA_i)$ by
\begin{equation}\label{eqn:maximum-norm}
    \norm{\functionQ}_{\bm{\xi}} \coloneqq \max_{(s,a)\in\bS_i\times\bA_i}\frac{\bigl\lvert Q(s,a)\bigr\rvert}{\xi(s)}, \quad (\functionQ\in \bF(\bS_i\times \bA_i)),  
\end{equation}
so that, for $\functionQ\in \bF(\bS_i\times \bA_i)$
\begin{multline}\label{eqn:exist_Q_star:1}
\bigl\lvert\bigl(\mathcal{T}_{i,0}^{\phi}\functionQ\bigr)(s,a)\bigr\rvert
= \bigl\lvert \sum_{s'\in\bS_i\backslash\{s_0\}}p_i(s,a,s')Q(s',\actionphi_i(s'))\bigr\rvert
\leq \sum_{s'\in\bS_i\backslash\{s_0\}}p_i(s,a,s')\bigl\lvert Q(s',\actionphi_i(s'))\bigr\rvert\\ \leq \sum_{s'\in\bS_i\backslash\{s_0\}} p_i(s,a,s') \norm{\functionQ}_{\bm{\xi}}\xi(s')
\leq \beta \xi(s) \norm{\functionQ}_{\bm{\xi}}. 
\end{multline}
This proves that $\mathcal{T}^\phi$ is a contraction and so has a unique fixed point by the Contraction Mapping Theorem.  
\endproof

\section{Proof of Proposition~\ref{prop:q-convergence}}
\label{app:prop:q-convergence}

\begin{lemma}\label{lemma:w_zero}
Given the history $\bh(t)$, $\bbE[w_{i,t}(s,a) |\bh(t)] =0$ for all $i\in[I]$ and $(s,a)\in\bS_i\times\bA_i$.    
\end{lemma}
\proof{Proof of Lemma~\ref{lemma:w_zero}.}
For $i\in[I]$ and $(s,a)\in\bS_i\times\bA_i$, if $\hat{\mathscr{N}}^{\phi}_{i,t}(s,a) = \emptyset$, then $w_{i,t}(s,a)=0$ based on the definition in \eqref{eqn:define_w}.
We then focus on the case where $\hat{\mathscr{N}}^{\phi}_{i,t}(s,a)\neq \emptyset$.
Based on the definition, for $i\in[I]$, $(s,a)\in\bS_i\times\bA_i$ and any $n\in[N_i]$, we have $\bbE[\bar{R}_{i,n}(s,a)|\bh(t)] = \bar{r}_i(s,a)$.
Then, we rewrite $\bbE[w_{i,t}(s,a)|\bh(t)]$ as follows.
\begin{equation}\label{eqn:lemma_w_zero:1}
    \bbE[w_{i,t}(s,a)|\bh(t)] = \bbE\biggl[\sum\limits_{s'\in\bS_i\backslash\{s_0\}}\Bigl(\frac{|\hat{\mathscr{N}}^{\phi}_{i,t}(s,a,s')|}{|\hat{\mathscr{N}}^{\phi}_{i,t}(s,a)|}-p_i(s,a,s')\Bigr)\hat{Q}(s',\phi_{t,i}(s')~\biggl|~\bh(t)\biggr].
\end{equation}
Recall the definition of $\hat{\mathscr{N}}^{\phi}_{i,t}(s,a)$ and $\hat{\mathscr{N}}^{\phi}_{i,t}(s,a,s')$ in \eqref{eq:subprocesses}.
Given $\bh(t)$, the only random part in \eqref{eqn:lemma_w_zero:1} is $|\hat{\mathscr{N}}^{\phi}_{i,t}(s,a,s')|$, which satisfies
\begin{equation}
    \bbE\Bigl[|\hat{\mathscr{N}}^{\phi}_{i,t}(s,a,s')|~\Bigl|~\bh(t)\Bigr] = \bbE\Bigl[\sum_{n\in\hat{\mathscr{N}}^{\phi}_{i,t}(s,a)}\indicator\bigl\{s^{\phi}_{i,n}(t+1) = s'\bigr\}~\Bigl|~\bh(t)\Bigr] = \bigl|\hat{\mathscr{N}}^{\phi}_{i,t}(s,a)\bigr|p_i(s,a,s'),
\end{equation}
where $\indicator\{\cdot\}$ is the indication function, and the second equality holds because, given $s^{\phi}_{i,n}(t) = s$ and $\phi_{i,n}(t)=a$, the event $s^{\phi}_{i,n}(t+1)=s'$ is independently and identically distributed for all $n\in \hat{\mathscr{N}}^{\phi}_{i,t}(s,a)$. It proves the lemma.

\endproof

\proof{Proof of Proposition~\ref{prop:q-convergence}.}
Based on Lemma~\ref{lemma:w_zero}, we have $\mathbb{E}\bigl[w_{i,t}(s,a) | \bh(t)\bigr] = 0$
for all $i\in[I]$, $t\in[T]_0$, and $(s,a)\in\bS_i\times \bA_i$, and, for the case with $\hat{\mathscr{N}}^{\psi}_{i,t}(s,a)\neq \emptyset$,
\begin{multline}
    \lvert w_{i,t}(s,a)\rvert \\= \Bigl\lvert \frac{1}{\bigl|\hat{\mathscr{N}}^{\phi}_{i,t}(s,a)\bigr|}\sum_{n\in\hat{\mathscr{N}}^{\phi}_{i,t}(s,a)}\bar{R}_{i,n}(s,a)-\bar{r}_i(s,a) +\sum_{s'\in\bS_i\backslash\{s_0\}}\hat{\functionQ}_{i,t}(s',\actionphi_{t,i}(s'))\Bigl(\frac{\bigl|\hat{\mathscr{N}}^{\phi}_{i,t}(s,a,s')\bigr|}{\bigl|\hat{\mathscr{N}}^{\phi}_{i,t}(s,a)\bigr|}-p_i(s,a,s')\Bigr)\Bigr\rvert\\
    \leq R_{max}+\max_{s'\in\bS_i}\bigl\lvert \hat{\functionQ}_{i,t}(s',\actionphi_{t,i}(s')\bigr\rvert
    \leq R_{max}+\max_{\begin{subarray}~s'\in\bS_i,\\a'\in\bA_i\end{subarray}} \bigl\lvert \hat{\functionQ}_{i,t}(s',a')\bigr\rvert,
\end{multline}
where $R_{max} \coloneqq \max_{t\in[T]_0,n\in\hat{\mathscr{N}}^{\psi}_i(s,a)}\bar{R}_{i,n}(s,a) - \bar{r}_i(s,a)< \infty$.
Hence, for all $i\in[I]$ and $(s,a)\in\bS_i\times \bA_i$, there exist $K,B<\infty$ such that
\begin{equation}\label{eqn:q-convergence:2}
\lvert w_{i,t}(s,a) \rvert^2 \leq K\max_{\begin{subarray}~s'\in\bS_i,\\a'\in\bA_i\end{subarray}}\hat{\functionQ}_{i,t}^2(s',a') + B.
\end{equation}
On the other hand, from Lemma~\ref{lemma:exist_Q_star} and \eqref{eqn:exist_Q_star:1}, for any $i\in[I]$, there exist $\beta\in[0,1)$ and $\bm{\xi}_i=(\xi_i(s):s\in\bS_i)\in \mathbb{R}_+^{|\bS_i|}$ such that, for all $t\in [T]_0$,
\begin{equation}\label{eqn:q-converge:3}
\Bigl\lVert \bigl(\mathcal{T}^{\bar{\phi}}_i \hat{\functionQ}_{i,t}\bigr) - \hat{\functionQ}^*\Bigr\rVert_{\bm{\xi}_i} \leq \beta \Bigl\lVert \hat{\functionQ}_{i,t} - \hat{\functionQ}^* \Bigr\rVert_{\bm{\xi}_i}.
\end{equation}

We can then prove the proposition by invoking \cite{bertsekas1996neuro}[Proposition 4.5].
\begin{externalTheorem}{Proposition 4.5 in \cite{bertsekas1996neuro}}
For the sequence $\bm{r}_t$ that takes values in $\mathbb{R}^N$ for some finite $N\in\mathbb{N}_+$ and is generated by
\begin{equation}
     r_{t+1}(n)= (1-\eta_t(n))r_t(n) + \eta_t(n)\Bigl((\mathcal{T}_t \bm{r}_t)(n) + w_t(n)+u_t(n)\Bigr), \forall n\in[N],
\end{equation}
we assume the following conditions.
\begin{enumerate}
    \item \label{condition:1} The stepsizes $\eta_u(n)$ are nonnegative and satisfy
\begin{equation}
    \sum_{t=0}^{\infty}\eta_t(n) = \infty,~\text{and}~\sum_{t=0}^{\infty}\eta^2_t(n) < \infty,~\forall n\in[N].
\end{equation}
\item \label{condition:2} The noise terms $w_t(n)$ satisfy (a) for every $n\in[N]$ and $t\in\mathbb{N}_0$, 
\[\mathbb{E}[w_t(n)|r_0(n),\ldots,r_t(n),\eta_0(n),\ldots,\eta_t(n)]=0;\] and (b) given any norm $\lVert \cdot \rVert $ on $\mathbb{R}^N$, there exist constants $A$ and $B$ such that 
\begin{equation}
    \mathbb{E}\Bigl[w^2_t(n)~\Bigl|~ r_0(n),\ldots,r_t(n),\eta_0(n),\ldots,\eta_t(n)\Bigr] \leq A + B \lVert \bm{r}_t \rVert^2, \forall n\in[N], t\in\mathbb{N}_0.
\end{equation}
\item \label{condition:3} There exists a vector $\bm{r}^*$, a positive vector $\bm{\xi}$, and a scaler $\beta\in[0,1)$, such that 
\begin{equation}
    \bigl\lVert \mathcal{T}_t\bm{r}_t - \bm{r}^*\bigr\rVert_{\bm{\xi}} \leq \beta \bigl\lVert \bm{r}_t - \bm{r}^*\bigr\rVert_{\bm{\xi}}, \forall t\in\mathbb{N}_0.
\end{equation}
\item \label{condition:4}
There exists a nonnegative random sequence $\theta_t$ that converges to zero with probability $1$, and is such that
\begin{equation}
    \lvert u_t(n)\rvert \leq \theta_t (\lVert \bm{r}_t\rVert_{\bm{\xi}} +1), \forall n\in[N],t\in\mathbb{N}_0.
\end{equation}
\end{enumerate}
Then, $\bm{r}_t$ converges to $\bm{r}^*$ with probability $1$.
\end{externalTheorem}

For the WCG-learning scheme in this paper, Conditions~\ref{condition:1}-\ref{condition:3} are  ensured by \eqref{eqn:assumption:stepsize}, \eqref{eqn:q-convergence:2}, and Lemma~\ref{lemma:exist_Q_star} and \eqref{eqn:q-converge:3}, respectively. Condition~\ref{condition:4} naturally holds in our case with $u_t(n)\equiv 0$.
It proves the proposition.

\endproof

\section{Pseudo-Code for OMPI}\label{app:ompi}

The pseudo-code for the OMPI algorithm at each time $t\in\mathbb{N}_0$ is provided in Algorithm~\ref{algo:ompi}.
\IncMargin{1em}
\begin{algorithm}
\small 
\linespread{0.5}\selectfont

\SetKwFunction{FOMPI}{OMPI}
\SetKwProg{Fn}{Function}{:}{\KwRet}
\SetKwInOut{Input}{Input}\SetKwInOut{Output}{Output}
\SetAlgoLined
\DontPrintSemicolon
\Input{Given $t\in[T]_0$, the history $\hat{\mathcal{\bm{H}}}_t$ upon time $t$, the iteration parameter $m$, the vectors of action labels $\LaVec_i$ for all $i\in[I]$.}
\Output{Updated estimated MP indices $\hatnuVec_t$, action variables for the primary policy $\psi$, the estimated Q factors $\hat{\bm{Q}}_t = \bigl(\hat{\bm{Q}}^{V,j,\varsigma}_t;\hat{\bm{Q}}^{U,j,\varsigma}_t;\hat{\bm{Q}}^{L,j,\varsigma}_t: j=1,2,\ldots,\max_{i\in[I]}|\bS_i|,\varsigma=0,1\bigr)$, and the  iteration parameter $m$.}
\Fn{\FOMPI{}}{
    \eIf{$\exists i\in[I]$, $m< \LenM_i$}{
        \tcc{The learning processes have not yet been stopped.}       
        $J\gets \max_{i\in[I]}\bigl|\bS_i\bigr|$ and construct $\bar{\phi}_m^{\varsigma}\in\oPhiLocal$ as described in Steps~\ref{step:period} and \ref{step:3}\;
        For all $j\in[J]$, $\Xi\in\{V,U,L\}$, and $\varsigma\in\{0,1\}$, update $\hat{\bm{Q}}^{\Xi,j,\varsigma}_t $ according to Step~\ref{step:4}.\; 
        Update the estimated MP indices $\hatnuVec_t$ according to Step~\ref{step:5}.\;
        \If{\eqref{eqn:stop} is satisfied for all $\Xi\in\{V,U,L\},j\in [J],\varsigma\in\{0,1\}$}{
            $s^*_i\gets$ a randomly selected element in $\arg\max_{s\in\bS_i: \La_{i,m}(s)\geq 2}\hat{\nu}_{i,t}(s,\La_{i,m}(s))$\;
            $\LaVec_{i,m+1}\gets \Ia_i^{s^*_i}\LaVec_{i,m}$\;
            $m\gets m+1$\;
        }
        $\bar{p}_t\gets$ a positive value in $(0,1)$\;
    }{
        $\bar{p}_t \gets 0$ and 
        $\hat{\bm{\nu}}_t\gets \hat{\bm{\nu}}_{t-1}$\;
    }
    \tcc{We then implement the primary policy $\psi$ with the updated $\hatnuVec_t$.}   
    \label{line:psi:begin}Uniformly randomly generate a number $\bar{p}$ from $[0,1]$\;
    \eIf{$\bar{p} > \bar{p}_t$}{
        $\LaVec^{\psi}(t) \gets \bm{1}^{\sum_{i\in[I]}N_i}$\;
        $\bB \gets \prod_{i\in[I]}[N_i]\times \bigl([|\bA_i|]\backslash\{1\}\bigr)$\;
        \While{$\bB\neq \emptyset$}{
            Uniformly randomly select an element $(i,n,\La)$ from $\bB$.\;
            \eIf{\eqref{eqn:DS_policy:singl_constraint} is not violated by replacing $\La(\rank)$ with $\La$}{
                $\La^{\psi}_{i,n}(t)\gets \La$\;
                $\bB\gets \bB\backslash\{(i,n,\La)\}$\;
            }{
                Break\;
            }
        }
    }{
        Take action $\LaVec^{\psi}(t) \gets \LaVec^{\text{MP}}(\hatnuVec_t,\bs^{\psi}(t))$\;
    }
    \label{line:psi:end}
}
\caption{The OMPI algorithm at each time $t$.}\label{algo:ompi}
\end{algorithm}
 \DecMargin{1em}

\section{Proof of Theorem~\ref{theorem:convergence-Z}}\label{app:theorem:convergence-Z}
Recall that Theorem~\ref{theorem:convergence-Z} focuses on policies $\varphi\in\Phi^1$, for which the action variables $\actionvarphi^{h}(t)$ are measurable through only current states $\bs^{\varphi}(t)$, time stamp $t$, and the ranking $\bm{\mathcal{R}}_t = \bm{\rank}_K(\hat{\bm{Q}}^1_t,\hat{\bm{Q}}^2_t,\ldots,\hat{\bm{Q}}^K_t)$.
For given $\varphi\in\Phi^1$, $\bar{\bm{\phi}}=(\bar{\phi}_t:t\in[T]_0)\in\bigl(\oPhiLocal\bigr)^{T+1}
$, $h\in\mathbb{N}_+$, and initial $\bs^{\varphi,h}(0)$ and $\hat{\bm{Q}}_0$, define $\varrho^{\varphi,\bar{\bm{\phi}},h}(T)\coloneqq (\bm{\mathcal{R}}_0,\bm{\mathcal{R}}_1,\ldots,\bm{\mathcal{R}}_T)$, which is a random variable taking values in $\mathscr{R}^{T+1}$ with certain probabilities.
For such $\varphi\in\Phi^1$ and a scenario conditioned on $\varrho^{\varphi,\bar{\bm{\phi}},h}(T)=\varrho$, we rewrite the state and action variables $s^{\varphi,h}_{i,n}(t)$, $Z^{\varphi,h}_{\iota}(t)$, and $\actionvarphi^{h}_{i,n}(t)$ as $s^{\varphi,h,\varrho}_{i,n}(t)$, $Z^{\varphi,h,\varrho}_{\iota}(t)$, and $\actionvarphi^{h,\varrho}_{i,n}(t)$, respectively, for the condition $\varrho$.

Define a state-action-time (SAT) triplet $(\iota,t)\in [\mathcal{I}]\times\bigl([T]_0\bigr)$, for which, since the element $\iota$ represents an SA pair, we use the term ``triplet" for the SAT.
Similar to the SA pairs, let $\mathcal{K}\coloneqq \Bigl|[\mathcal{I}]\times\bigl([T]_0\bigr)\Bigr|$, and label all the SAT triplets as $\kappa\in[\mathcal{K}]$, where, for any $t\in[T]_0$, we list all the SATs $(\iota,t+1)$ ($\iota\in[\mathcal{I}]$) in bulk preceding to those SATs $(\iota',t)$ ($\iota'\in[\mathcal{I}]$).
We denote the $\kappa$th SAT by $(\iota_{\kappa},t_{\kappa})$, and alternatively refer to such an SAT as SAT $\kappa\in[\mathcal{K}]$ or SAT $(i_{\kappa},s_{\kappa},\actionphi_{\kappa},t_{\kappa})$, where $(i_{\kappa},s_{\kappa})$ and $(i_{\kappa},\actionphi_{\kappa})$ specify the state and action, respectively, associated with SAT $\kappa$.
For $t\in[T]_0$, define the set of all the associated SATs as $\mathscr{K}_t\coloneqq \{\kappa\in[\mathcal{K}]~|~t_{\kappa} = t\}$.

For given $\varphi\in\Phi^1$, $\bar{\bm{\phi}}\in\bigl(\oPhiLocal\bigr)^{T+1}$, $h\in\mathbb{N}_+$,  $\varrho^{\varphi,\bar{\bm{\phi}},h}(T)=\varrho$, $t\in[T]_0$, and $\kappa\in[\mathcal{K}]$, define $Z^{\varphi,h,\varrho}_{\kappa}(t)$ as the proportion of bandit processes $\bigl\{s^{\varphi,h,\varrho}_{i,n}(t), t\in[T]_0\bigr\}$ ($i\in[I],n\in[N_i]$) that are in SA pair $\iota\in[\mathcal{I}]$ at time $t$; that is,
\begin{equation}\label{eqn:app:convergence_Z:1}
Z^{\varphi,h,\varrho}_{\iota}(t) = \frac{1}{\sum_{i\in[I]}N_i}\Biggl|\biggl\{(i,n)~\biggl|~
i\in[I],n\in[N_i],~s^{\varphi,h,\varrho}_{i,n}(t)=s_{\iota},\actionvarphi^{h,\varrho}_{i,n}(t)=\actionphi_{\iota}\biggr\}\Biggr|,
\end{equation}
where recall $N_i = hN_i^0$ for $i\in[I]$.
Let $\bm{Z}^{\varphi,h,\varrho}(t)\coloneqq (Z^{\varphi,h,\varrho}_{\iota}(t): \iota\in[\mathcal{I}])$.

Given policy $\varphi\in\Phi^1$ and $\varrho\in\mathscr{R}^{T+1}$, for $i\in[I]$, $t\in[T]_0$, and $\iota,\iota'\in[\mathcal{I}]$ with $i_{\iota}=i$,
define a random variable $T^{\varphi,\varrho}_{i,t}(\iota,\iota')$ that is exponentially distributed with rate 
\begin{equation}\label{eqn:app:convergence_Z:2}
\frac{1}{\mathbb{E}\bigl[T^{\varphi,\varrho}_{i,t}(\iota,\iota')\bigr]} = p^{\varphi,\varrho}_{i,t}(\iota,\iota')\coloneqq\mathbb{P}\Bigl\{s^{\varphi,h,\varrho}_{i,n}(t+1) = s_{\iota'},\actionvarphi^{h,\varrho}_{i,n}(t+1) = \actionphi_{\iota'}~\Bigl|~s^{\varphi,h,\varrho}_{i,n}(t) = s_{\iota},\actionvarphi^{h,\varrho}_{i,n}(t) = \actionphi_{\iota}\Bigr\},
\end{equation}
for any $n\in[N_i]$.
For $\kappa\in[\mathcal{K}]$ and $\iota'\in[\mathcal{I}]$, consider $N_{i_{\kappa}}$ homogeneous Poisson point (HPP) processes, denoted by $\{M^{\varphi,\varrho}_{\kappa,\iota'}(n,\tau), \tau\geq 0\}$ with $n\in[N_{i_{\kappa}}]$, on the positive half line, where the interval between two events is in distribution equivalent to $ T^{\varphi,\varrho}_{i_{\kappa},t_{\kappa}}(\iota_{\kappa},\iota')$.
For each time slot $t\in[T]_0$, $i\in[I]$, and $n\in[N_i]$, we associated the HPP processes $\{M^{\varphi,\varrho}_{\kappa,\iota'}(n,\tau), \tau\geq 0\}$, where $\kappa\in\mathscr{K}_t$ with $i_{\kappa} = i$ and $\iota'\in[\mathcal{I}]$, to the bandit process $s^{\varphi,h,\varrho}_{i,n}(t)$.
We assume that no simultaneous event occurrences among all the HPP processes almost surely.
For $t\in[T]_0$, define 
\begin{equation}
    M^{\varphi,\varrho}_t(\tau)\coloneqq \sum_{(\kappa,\iota')\in\mathscr{K}_t\times[\mathcal{I}]}\sum_{n\in[N_{i_{\kappa}}]}M^{\varphi,\varrho}_{\kappa,\iota'}(n,\tau),
\end{equation}
where let $\tau^{\varphi,\varrho}_t(m)\coloneqq \min \{\tau\geq 0~|~M^{\varphi,\varrho}_t(\tau)\geq m\}$, representing the occurrence time of the $m$th event in such a process, and $\tau^{\varphi,\varrho}_t(0) \equiv 0$.
In particular, when this $m$th event is contributed by the process $\bigl\{M^{\varphi,\varrho}_{\kappa,\iota'}(n,\tau),\tau\geq 0\bigr\}$, let variables $\bigl(\kappa^{\varphi,\varrho}_t(m),\iota^{\varphi,\varrho}_t(m),n^{\varphi,\varrho}_t(m)\bigr) = (\kappa,\iota',n)$.




For $\tau\in\mathbb{R}_0$, define a vector of random variables $\bm{\xi}^{\varphi,\varrho}_{\kappa,\tau}(n) = (\xi^{\varphi,\varrho}_{\kappa,\kappa',\tau}(n):\kappa'\in [\mathcal{K}])$ such that, if $(\kappa,\iota_{\kappa'},n) = (\kappa^{\varphi,\varrho}_{t_{\kappa}}(m),\iota^{\varphi,\varrho}_{t_{\kappa}}(m),n^{\varphi,\varrho}_{t_{\kappa}}(m))$ for $m=\min\{m'\in\mathbb{N}_+~|~\tau^{\varphi,\varrho}_t(m')> \tau\}$ and $t_{\kappa'} = t_{\kappa}+1$, then
\begin{equation}\label{eqn:app:convergence_Z:3}
    \xi^{\varphi,\varrho}_{\kappa,\kappa',\tau}(n)\coloneqq \frac{1}{\tau^{\varphi,\varrho}_t(m)-\tau^{\varphi,\varrho}_t(m-1)},
\end{equation}
otherwise, $\xi^{\varphi,\varrho}_{\kappa,\kappa',\tau}(n) = 0$.
Here, for all $\tau\geq 0$, there exists the unique $(\kappa,\kappa',n)$ such that $\xi^{\varphi,\varrho}_{\kappa,\kappa',\tau}(n) > 0$.
For $\tau<0$, define $\xi^{\varphi,\varrho}_{\kappa,\kappa',\tau}(n) =0$ almost surely for all $\kappa,\kappa'\in[\mathcal{K}]$ and $n\in[N_{i_{\kappa}}]$.
Note that the subscript $\tau$ for $\xi^{\varphi,\varrho}_{\kappa,\kappa',\tau}(n)$ is an artificial timeline for theoretical analysis (related to the artificial HPP processes) of the asymptotic regime, which is different from the realistic time points that are included in $\kappa$ and $\kappa'$.
In this context, the random vector $\bm{\xi}^{\varphi,\varrho}_{\kappa,\tau}(n)$ is identically distributed for all $\tau\geq 0$, and the trajectory $\xi^{\varphi,\varrho}_{\kappa,\kappa',\tau}(n)$ is almost continuous in $\tau\geq 0$ except finitely many discontinuous points of the first kind within every finite time interval.
For two appropriately selected first-kind discontinuity points $\tau_1,\tau_2\in[0,\infty)$, the integral $\int_{\tau_1}^{\tau_2}\xi^{\varphi,\varrho}_{\kappa,\kappa',\tau}(n) d \tau$ is an integer representing the potential number of transitions that the process $s^{\varphi,\varrho}_{i_{\kappa},n}(t)$ changes from SA pair $\iota_{\kappa}$ at time $t_{\kappa}$ to SA pair $\iota_{\kappa'}$ at time $t_{\kappa'}$.
Based on the definition in \eqref{eqn:app:convergence_Z:2}, this transition number is positive only if $t_{\kappa'}=t_{\kappa}+1$.
Such a $\xi^{\varphi,\varrho}_{\kappa,\kappa',\tau}(n)$ is considered as an event that may trigger a state transition of the corresponding process.
We use the word ``potential" because the process $\bigl\{s^{\varphi,\varrho}_{i_{\kappa},n}(t), t\in[T]_0\bigr\}$ may not be in 
SA pair $\iota_{\kappa}$ at time $t_{\kappa}$, causing less number of real transitions. 
Let $\bm{\xi}^{\varphi,h,\varrho}_\tau\coloneqq (\xi^{\varphi,\varrho}_{\kappa,\kappa',\tau}(n): \kappa,\kappa'\in[\mathcal{K}], n\in[N_{i_{\kappa}}])$ which takes values in $\mathbb{R}_0^{\mathcal{K}\sum_{\kappa\in[\mathcal{K}]}N_{i_{\kappa}}}$, where recall $N_i = hN_i^0$.

For $h\in\mathbb{N}_+$, $\bm{x}\in\mathbb{R}_0^{\mathcal{K}}$, $\bm{\xi}=(\xi_{\kappa,\kappa'}(n):\kappa,\kappa'\in[\mathcal{K}],n\in[N_{i_{\kappa}}])\in\mathbb{R}_0^{\mathcal{K}\sum_{\kappa'\in[\mathcal{K}]}N_{i_{\kappa'}}}$, and $\kappa,\kappa'\in[\mathcal{K}]$, define
\begin{equation}\label{eqn:app:convergence_Z:3}
Q^{h}(\kappa,\kappa',\bm{x},\bm{\xi}) \coloneqq 
\text{$\indicator$}\biggl\{\sum_{\kappa''\in [\mathcal{K}],t_{\kappa''} <t_{\kappa}}x_{\kappa''} =0\biggr\}\sum_{n=h\lceil x^-_{\kappa-1}/h\rceil+1}^{h\lceil x^-_{\kappa}/h\rceil }\xi_{\kappa,\kappa'}(n) + f^{h,\Lipschitza}_{\kappa,\kappa'}(\bm{x},\bm{\xi}),
\end{equation}
where $x^-_{\kappa} = \sum_{\kappa'=1,i_{\kappa'}=i_{\kappa}}^{\kappa}x_{\kappa'}$ for $\kappa\in[\mathcal{K}]$ and $x^-_0 = 0$, and $f^{h,\Lipschitza}_{\kappa,\kappa'}(\bm{x},\bm{\xi})$ are appropriate functions for which, when $h=1$, $Q^h(\kappa,\kappa',\bm{x},\bm{\xi})$ is Lipschitz continuous in $\bm{x}$ for any given $\bm{\xi}$ and $\Lipschitza\in(0,1)$. 
As in \cite{fu2018restless,fu2020energy,fu2024patrolling}, such functions $f^{h,\Lipschitza}_{\kappa,\kappa'}(\bm{x},\bm{\xi})$ can be constructed by incorporating the Dirac delta function.
In particular, $f^{h,\Lipschitza}_{\kappa,\kappa'}(\bm{x},\bm{\xi})$ is continuous in $\Lipschitza\in(0,1)$ and linear in $\bm{\xi}$, satisfying $\lim_{\Lipschitza\downarrow 0}\bigl(Q^h(\kappa,\kappa',\bm{x},\bm{\xi}) - Q^h(\kappa,\kappa',\lceil\bm{x}\rceil,\bm{\xi})\bigr) = 0$, where $\lceil \bm{x}\rceil = (\lceil x_{\kappa} \rceil: \kappa\in[\mathcal{K}])$, and $\lim_{\Lipschitza\downarrow0}d~ f^{h,\Lipschitza}_{\kappa,\kappa'}(\bm{x},\bm{\xi})/d\Lipschitza = 0$. 
We provide an example of such $f^{h,\Lipschitza}_{\kappa,\kappa'}$ in Appendix~\ref{app:dirac-delta}.
For $\bm{x}\in\mathbb{R}^{\mathcal{K}}$ with negative elements, define $Q^h(\kappa,\kappa',\bm{x},\bm{\xi}) \coloneqq Q^h(\kappa,\kappa',(\bm{x})^+,\bm{\xi})$, where $(\bm{x})^+\coloneqq \bigl(\max\{x_{\kappa},0\}:\kappa\in[\mathcal{K}]\bigr)$.
For $h\in\mathbb{N}_+$, $\kappa\in[\mathcal{K}]$, $\bm{x}\in\mathbb{R}_0^{\mathcal{K}}$, $\bm{\xi}\in\mathbb{R}_0^{\mathcal{K}\sum_{\kappa'\in[\mathcal{K}]}N_{i_{\kappa'}}}$, define a function
\begin{equation}\label{eqn:app:convergence_Z:4}
b^h_{\kappa}(\bm{x},\bm{\xi})\coloneqq \sum_{\kappa'\in[\mathcal{K}]}\Bigl(Q^h(\kappa',\kappa,\bm{x},\bm{\xi})-Q^h(\kappa,\kappa',\bm{x},\bm{\xi})\Bigr),
\end{equation}
which, for $h=1$, is Lipschitz continuous in $\bm{x}$ and $\bm{\xi}$.
Let $b^h(\bm{x},\bm{\xi})\coloneqq \bigl(b^h_{\kappa}(\bm{x},\bm{\xi}):\kappa\in[\mathcal{K}]\bigr)$, 
and, based on \eqref{eqn:app:convergence_Z:3} and \eqref{eqn:app:convergence_Z:4}, there exists matrix 
$\tilde{\mathcal{Q}}^h(\bm{x})\in\mathbb{R}^{\mathcal{K}\times\bigl(\mathcal{K}\sum_{\kappa'\in[\mathcal{K}]}N_{i_{\kappa'}}\bigr)}$ such that $b^h(\bm{x},\bm{\xi})=\tilde{\mathcal{Q}}^h(\bm{x})\bm{\xi}$.

For $U\in\mathbb{R}_0\cup\{\infty\}$, $\bm{x}\in\mathbb{R}_0^{\mathcal{K}}$ and $\bm{\xi}\in\mathbb{R}_0^{\mathcal{K}\sum_{\kappa\in[\mathcal{K}]}N_{i_{\kappa}}}$, we define an adapted version of $b^h(\bm{x},\bm{\xi})$ as
\begin{equation}\label{eqn:app:convergence_Z:5}
b^{h,U}(\bm{x},\bm{\xi})\coloneqq b^{h,U}\bigl(\bm{x},\min\bigl\{\bm{\xi},U\bigr\}\bigr),
\end{equation}
where $\min\bigl\{\bm{\xi},U\bigr\} \coloneqq \bigl(\min\{\xi_{\kappa,\kappa'}(n),U\}: \kappa,\kappa'\in[\mathcal{K}],n\in[N_{i_{\kappa}}]\bigr)$.
Let $b^U(\bm{x},\bm{\xi})$ represent the special case $b^{h,U}(\bm{x},\bm{\xi})$ with $h=1$.
It is obvious that $b^{h,\infty}(\bm{x},\bm{\xi})\coloneqq\lim_{U\rightarrow \infty}b^{h,U}(\bm{x},\bm{\xi}) = b^h(\bm{x},\bm{\xi})$.

For the special case with $h=1$ and any $U\in\mathbb{R}_0\cup\{\infty\}$, we construct a trajectory $\bm{X}^{\sigma,U}_{\tau}$ on $\tau\geq 0$ that satisfies
\begin{equation}\label{eqn:app:convergence_Z:6}
\dot{\bm{X}}^{\sigma,U}_\tau = b^U(\bm{X}^{\sigma,U}_{\tau},\bm{\xi}^{\varphi,1,\varrho}_{\tau/\sigma}),
\end{equation}
with given $\bm{X}^{\sigma,U}_0=\bm{x}_0\coloneqq (x_{\kappa,0}:\kappa\in[\mathcal{K}])$, where  $\bm{\xi}^{\varphi,1,\varrho}_{\tau}$ is the special case of $\bm{\xi}^{\varphi,h,\varrho}_{\tau}$ with $h=1$.
Let $\bm{X}^{\sigma}_\tau$ represent the special case $\bm{X}^{\sigma,\infty}_{\tau}\coloneqq \lim_{U\rightarrow \infty}\bm{X}^{\sigma,\infty}_{\tau}$.

Consider a case where $x_{\kappa,0} = Z^{\varphi,h,\varrho}_{\iota_{\kappa}}(0)\sum_{i\in[I]}N_i$ for all $\kappa\in\mathscr{K}_0$; and $x_{\kappa,0} = 0$ for other $\kappa\in[\mathcal{K}]\backslash\mathscr{K}_0$.
For $\sigma=h=1$, if we plug in $\bm{X}^{\sigma}_0 = \bm{x}_0$, then, for any time slot $t\in[T]$ of the real process and its associated artificial time point
\begin{equation}\label{eqn:app:convergence_Z:7}
    \tau^{\sigma}(t) \coloneqq \min\Bigl\{\tau\in\mathbb{R}_0~\Bigl|~\sum_{\kappa\in\bigcup_{t'<t}\mathscr{K}_{t'}}X^{\sigma}_{\kappa,\tau} = 0\Bigr\},
\end{equation}
it follows that $\lim_{\Lipschitza\downarrow 0}X^{\sigma}_{\kappa,\tau^{\sigma}(t)} \sim Z^{\varphi,h,\varrho}_{\iota_{\kappa}}(t)\sum_{i\in[I]}N_i$ for all $\kappa \in\mathscr{K}_t$, where $\sim$ is equivalence in distribution.
Recall that such an artificial process $\bm{X}^{\sigma}_{\tau}$ is used for the completeness of the theoretical analysis.
It coincides in distribution with the real process $\bm{Z}^{\varphi,h,\varrho}(t)$ for the special points $\tau = \tau^{\sigma}(t)$, but not equal in general.

\begin{lemma}\label{lemma:convergence_b}
For any $U\in\mathbb{R}_0\cup\{\infty\}$, $\delta>0$, $\tau\in\mathbb{R}_0$, and $\bm{x}\in\mathbb{R}_0^{\mathcal{K}}$, there exists $\bar{b}^U(\bm{x})=\mathbb{E}b^U(\bm{x},\bm{\xi}^{\varphi,1,\varrho}_{\tau})\in\mathbb{R}_0^{\mathcal{K}}$
such that
\begin{equation}\label{eqn:lemma:convergence_b}
\lim_{\bar{T}\rightarrow +\infty}\mathbb{P}\biggl\{\Bigl\lVert \frac{1}{\bar{T}}\int_{\tau}^{\tau+\bar{T}}b^U(\bm{x},\bm{\xi}^{\varphi,1,\varrho}_{\tau'}) d\tau' -\bar{b}^U(\bm{x})\Bigr\rVert > \delta\biggr\}=0,
\end{equation}
uniformly for all $\tau\geq 0$. We recall that, since $\bm{\xi}^{\varphi,1,\varrho}_{\tau}$ is identically distributed for all $\tau\geq 0$, $\mathbb{E}b^U(\bm{x},\bm{\xi}^{\varphi,1,\varrho}_{\tau})$ is independent for $\tau\geq 0$.
\end{lemma}

\proof{Proof of Lemma~\ref{lemma:convergence_b}.}
Let $\bar{b}^U(\bm{x})\coloneqq \tilde{\mathcal{Q}}^1(\bm{x})\bm{\lambda}^U$ for some $\bm{\lambda}^U\in\mathbb{R}_0^{\mathcal{K}}$.
For any $U<\infty$, $\delta>0$, $\tau,\bar{T}\in\mathbb{R}_0$, and $\bm{x}\in\mathbb{R}_0^{\mathcal{K}}$,
\begin{multline}\label{eqn:lemma:convergence_b:1}
\mathbb{P}\biggl\{\Bigl\lVert \frac{1}{\bar{T}}\int_{\tau}^{\tau+\bar{T}}b^U(\bm{x},\bm{\xi}^{\varphi,1,\varrho}_{\tau'}) d\tau' -\bar{b}^U(\bm{x})\Bigr\rVert > \delta\biggr\} \\
=\mathbb{P}\biggl\{\Bigl\lVert \tilde{\mathcal{Q}}^1(\bm{x})\Bigl(\frac{1}{\bar{T}}\int_{\tau}^{\tau+\bar{T}}\min\bigl\{\bm{\xi}^{\varphi,1,\varrho}_{\tau'},U\bigr\} d\tau' -\bm{\lambda}^U\Bigr)\Bigr\rVert> \delta\biggr\}.
\end{multline}
Based on the definition of $\bm{\xi}^{\varphi,h,\varrho}_{\tau}$, for any $\tau,\bar{T}\geq 0$, $\min\bigl\{\bm{\xi}^{\varphi,1,\varrho}_{\tau},U\bigr\}$ is always bounded, and there exists $\bm{\lambda}^{\varphi,\varrho}(U)\coloneqq \frac{1}{\bar{T}}\mathbb{E}\Bigl[\int_{\tau}^{\tau+\bar{T}}\min\bigl\{\bm{\xi}^{\varphi,1,\varrho}_{\tau'},U\bigr\}d\tau'\Bigr] < \infty$.
For $\bm{\lambda}^U = \bm{\lambda}^{\varphi,\varrho}(U)$, based on the Law of Large Numbers, we obtain \eqref{eqn:lemma:convergence_b} uniformly for all $\tau\geq0$.

For the case with $U\rightarrow\infty$, any $\delta>0$, $\tau,\bar{T}\in\mathbb{R}_0$, and $\bm{x}\in\mathbb{R}_0^{\mathcal{K}}$, we obtain
\begin{multline}\label{eqn:lemma:convergence_b:2}
\mathbb{P}\biggl\{\Bigl\lVert \frac{1}{\bar{T}}\int_{\tau}^{\tau+\bar{T}}b(\bm{x},\bm{\xi}^{\varphi,1,\varrho}_{\tau'}) d\tau' -\bar{b}^{\infty}(\bm{x})\Bigr\rVert > \delta\biggr\} \\
\leq \mathbb{P}\biggl\{\Bigl\lVert \tilde{\mathcal{Q}}^1(\bm{x})\frac{1}{\bar{T}}\Bigl\lfloor\int_0^{\bar{T}}\bm{\xi}^{\varphi,1,\varrho}_{\tau'} d\tau'\Bigr\rfloor -\bar{b}^{\infty}(\bm{x})\Bigr\rVert + \frac{o(\bar{T})}{\bar{T}}> \delta\biggr\}\\
=\mathbb{P}\biggl\{\Bigl\lVert \tilde{\mathcal{Q}}^1(\bm{x})\Bigl(\frac{1}{\bar{T}}\Bigl\lfloor\int_0^{\bar{T}}\bm{\xi}^{\varphi,1,\varrho}_{\tau'} d\tau'\Bigr\rfloor -\bm{\lambda}^{\infty}\Bigr)\Bigr\rVert + \frac{o(\bar{T})}{\bar{T}}> \delta\biggr\},
\end{multline}
where, for any vector $\bm{v}\in\mathbb{R}^M$, $\lfloor \bm{v}\rfloor \coloneqq \bigl(\lfloor v_m\rfloor: m\in[M]\bigr)$, $\bar{b}^{\infty} \coloneqq \lim_{U\rightarrow \infty}\bar{b}^U$, and $\bm{\lambda}^{\infty}\coloneqq \lim_{U\rightarrow \infty}\bm{\lambda}^U$. 
Based on the definition of $\bm{\xi}^{\varphi,h,\varrho}_{\tau}$, for any $\tau,\bar{T}\in\mathbb{R}_0$, each element of $\bigl\lfloor \int_0^{\bar{T}}\bm{\xi}^{\varphi,1,\varrho}_{\tau'}d\tau'\bigr\rfloor$ follows a Poisson distribution with expectation $\bar{T}\bm{\lambda}^{\varphi,\varrho}$, where $\bm{\lambda}^{\varphi,\varrho}\coloneqq \bigl(\lambda^{\varphi,\varrho}_{\kappa,\kappa'}(n): \kappa,\kappa'\in[\mathcal{K}],n\in[N_{i_{\kappa}}]\bigr)$ with
$\lambda^{\varphi,\varrho}_{\kappa,\kappa'}(n) =    p^{\varphi,\varrho}_{i_{\kappa},t_{\kappa}}(\iota_{\kappa},\iota_{\kappa'})$ for all $n\in[N_{i_{\kappa}}]$.

From the Law of Large Numbers, plugging  $\bm{\lambda}^{\infty}=\bm{\lambda}^{\varphi,\varrho}$ in \eqref{eqn:lemma:convergence_b:2}, we obtain 
\begin{multline}\label{eqn:lemma:convergence_b:3}
\lim_{\bar{T}\rightarrow +\infty}\sup_{\tau\geq 0}\mathbb{P}\biggl\{\Bigl\lVert \frac{1}{\bar{T}}\int_{\tau}^{\tau+\bar{T}}b(\bm{x},\bm{\xi}^{\varphi,1,\varrho}_{\tau'}) d\tau' -\bar{b}(\bm{x})\Bigr\rVert > \delta\biggr\}\\
\leq\lim_{\bar{T}\rightarrow +\infty}\mathbb{P}\biggl\{\Bigl\lVert\tilde{\mathcal{Q}}(\bm{x})\Bigl(\frac{1}{\bar{T}}\Bigl\lfloor\int_0^{\bar{T}}\bm{\xi}^{\varphi,1,\varrho}_{\tau'}d\tau'\Bigr\rfloor - \bm{\lambda}^{\varphi,\varrho}\Bigr)\Bigr\rVert+\frac{o(\bar{T})}{\bar{T}} > \delta\biggr\}=0.
\end{multline}
That is, for $U\rightarrow \infty$, any $\delta>0$ and $\bm{x}\in\mathbb{R}_0^{\mathcal{K}}$, there exists $\bar{b}^{\infty}(\bm{x}) = \tilde{\mathcal{Q}}^1\bm{\lambda}^{\varphi,\varrho}$ such that
\eqref{eqn:lemma:convergence_b} is satisfied uniformly for all $\tau\geq 0$. It proves the lemma.

\endproof

For $U\in\mathbb{R}_0\cup\{\infty\}$, define $\bar{\bm{x}}^U_{\tau}$ as the unique solution of $\dot{\bar{\bm{x}}}^U_{\tau} = \bar{b}^U(\bar{\bm{x}}_{\tau})$ with given $\bar{\bm{x}}_0\in\mathbb{R}^{\mathcal{K}}$, for which the unique existence of $\bar{\bm{x}}^U_{\tau}$ is guaranteed by Picard's Existence Theorem in \cite{coddington1955theory}.

From \cite[Theorem 2.1 in Chapter 7]{freidlin2012random}, if $\bar{\bm{x}}^U_t$ and $\bar{b}^U(\bm{x})$ exist and satisfy \eqref{eqn:lemma:convergence_b} uniformly for $\tau\geq 0$, and $\lVert b^U(\bm{x},\bm{\xi}^{\varphi,1,\varrho}_t)\rVert^2 < \infty$ for all $\bm{x}\in\mathbb{R}_0^{\mathcal{K}}$, then, for any $0<\bar{T}<\infty$ and $\delta >0$,
\begin{equation}\label{eqn:app:convergence_Z:8}
\lim_{\sigma\downarrow 0} \mathbb{P}\Bigl\{\sup_{0\leq \tau\leq \bar{T}}\bigl\lVert \bm{X}^{\sigma,U}_{\tau} - \bar{\bm{x}}^U_\tau\bigr\rVert > \delta\Bigr\}=0,
\end{equation}
where the initial point $\bm{X}^{\sigma,U}_0 = \bar{\bm{x}}^U_0$ is given.
From Lemma~\ref{lemma:convergence_b}, \eqref{eqn:app:convergence_Z:8} is achieved for all given $\Lipschitza\in(0,1)$ and $U\in\mathbb{R}_0$. 
Together with the continuity of $\bm{X}^{\sigma,U}_\tau$ and $\bar{\bm{x}}^U_\tau$ in $\Lipschitza\in(0,1)$ and $U\in(0,\infty)$ and the boundness of $\lim_{\Lipschitza\downarrow 0}\bm{X}^{\sigma,U}_\tau$, $\lim_{\Lipschitza\downarrow 0}\bar{\bm{x}}^U_\tau$, 
$\lim_{U\rightarrow \infty}\bm{X}^{\sigma,U}_\tau$, and $\lim_{U\rightarrow \infty}\bar{\bm{x}}^U_\tau$, 
\eqref{eqn:app:convergence_Z:8} also holds for $\Lipschitza\downarrow 0$ and $U\rightarrow \infty$.
That is, for any $0<\bar{T}<\infty$ and $\delta >0$,
\begin{equation}\label{eqn:app:convergence_Z:9}
    \lim_{\sigma\downarrow 0}\mathbb{P}\Bigl\{\sup_{0\leq \tau\leq \bar{T}}\bigl\lVert \bm{X}^{\sigma}_{\tau} - \bar{\bm{x}}_{\tau}\bigr\rVert > \delta\Bigr\} = 0,
\end{equation}
where $\bm{X}^{\sigma}_{\tau}$ and $\bar{\bm{x}}_{\tau}$ are those in the case with $\Lipschitza\downarrow 0$.

\proof{Proof of Theorem~\ref{theorem:convergence-Z}.}
In this proof, we consider the case with $\Lipschitza\downarrow 0$ and $U\rightarrow \infty$, and let $\sigma = 1/h$, where $h$ is the scaling parameter.
For any $\tau,\bar{T}\in\mathbb{R}_0$,
\begin{equation}\label{eqn:app:convergence_Z:10}
\int_\tau^{\tau+\bar{T}} b(\bm{X}^{\sigma}_{\tau},\bm{\xi}^{\varphi,1,\varrho}_{\tau/\sigma})d\tau = \sigma\int_{\tau/\sigma}^{(\tau+\bar{T})/\sigma}b(\bm{X}^{\sigma}_{\sigma\tau},\bm{\xi}^{\varphi,1,\varrho}_\tau)d\tau = \frac{1}{h}\int_{h\tau}^{h(\tau+\bar{T})}b(\bm{X}^{\sigma}_{\tau/h},\bm{\xi}^{\varphi,1,\varrho}_\tau)d\tau.
\end{equation}

Given $h\in\mathbb{N}_+$ and $\bm{X}^{\sigma}_0=\bar{\bm{x}}_0=\bm{x}_0\in\mathbb{N}_0^{\mathcal{K}}$, consider a trajectory $\bm{X}^h_\tau$ which satisfies $\dot{\bm{X}}^h_\tau = \frac{1}{h }b^h\bigl(h\bm{X}^h_{\tau},\bm{\xi}^{\varphi,h,\varrho}_{\tau}\bigr)$ and $\bm{X}^h_0 = \bm{x}_0$.
Define a set $\mathscr{X}^h \subset \mathbb{R}_0^{\mathcal{K}}$ where, for any $\bm{x}=(x_{\kappa}:\kappa\in[\mathcal{K}])\in\mathscr{X}^h$, each element $x_{\kappa}$ can be rewritten as $\frac{n_{\kappa}}{h}$ with $n_{\kappa}\in\mathbb{N}_0$. 
Based on the definitions of $\bm{\xi}^{\varphi,h,\varrho}_\tau$ and $Q^h$, for $t\in[T]_0$, and $\bar{T} = \tau^{\varphi,1,\varrho}_{t}(m)$ for some $m\in\mathbb{N}_0$, $\bm{X}^{\sigma}_{\bar{T}}$ and $\bm{X}^h_{\bar{T}}$ must take some values in $\mathscr{X}^h$.
In particular, for $\kappa\in[\mathcal{K}]$, given $\tau_m=\tau^{\varphi,1,\varrho}_{t_{\kappa}}(m)$ for $m\in\mathbb{N}_0$,
if $\bm{X}^{\sigma}_{\tau_m} = \bm{X}^h_{\tau_m}=\bm{x}$ with  $\bm{x}\in\mathscr{X}^h$ and $x_{\kappa}>0$, then, let $\bar{T}^{\sigma}$ and $\bar{T}^h$ represent the time such that $\lvert \int_{h\tau_m}^{h(\tau_m+\bar{T}^{\sigma})}b_{\kappa}(\bm{X}^{\sigma}_{\tau/h},\bm{\xi}^{\varphi,1,\varrho}_{\tau})d\tau\rvert= 1$ and $\lvert \int_{\tau_m}^{\tau_m+\bar{T}^h}b^h_{\kappa}(h\bm{X}^h_{\tau},\bm{\xi}^{\varphi,h,\varrho}_{\tau})d\tau\rvert= 1$, respectively, that are exponentially distributed with rate $h\lceil x_{\kappa}/h\rceil\sum_{\kappa'\in\mathscr{K}_{t_{\kappa}+1}}p^{\varphi,\varrho}_{i_{\kappa},t_{\kappa}}(\iota_{\kappa},\iota_{\kappa'})$, we have $\bar{T}^{\sigma}\sim \bar{T}^h$, where $\sim$ means equivalence in distribution.
Hence, for given $\bar{T}\in\mathbb{R}_0$, if $\bm{X}^{\sigma}_{\tau_m} = \bm{X}^h_{\tau_m}=\bm{x}$, then, for any $\kappa\in[\mathcal{K}]$,
\begin{equation}\label{eqn:app:convergence_Z:11}
    \frac{1}{h}\int_{h\tau_m}^{h(\tau_m+\bar{T})}b_{\kappa}(\bm{X}^{\sigma}_{\tau/h},\bm{\xi}^{\varphi,1,\varrho}_{\tau})d\tau \sim \frac{1}{h}\int_{\tau_m}^{\tau_m+\bar{T}}b^h_{\kappa}(h\bm{X}^h_{\tau},\bm{\xi}^{\varphi,h,\varrho}_{\tau})d \tau,
\end{equation}
leading to $\bm{X}^{\sigma}_\tau \sim \bm{X}^h_\tau$ for all $\tau\in\mathbb{R}_0$.

Given $\bm{X}^{\sigma}_0=\bm{X}^h_0=\bar{\bm{x}}_0= \sum_{i\in[I]}N_i^0\bm{\scrz}_0$ that takes values in $\mathscr{N}_0^{\mathcal{K}}$, consider a trajectory \[\bm{\mathcal{Z}}^h_\tau=\frac{\bm{X}^h_\tau}{\sum_{i\in[I]}N_i^0},\] 
which satisfies
\[\dot{\bm{\mathcal{Z}}}^h_\tau = \frac{1}{h\sum_{i\in[I]}N_i^0}b^h\bigl(h\sum_{i\in[I]}N_i^0\bm{\mathcal{Z}}^h_{\tau},\bm{\xi}^{\varphi,h,\varrho}_{\tau}\bigr),\] 
and $\bm{\mathcal{Z}}^h_0 = \bm{\scrz}_0$.
From \eqref{eqn:app:convergence_Z:11}, we obtain that, for $\tau \geq 0$,
\begin{equation}\label{eqn:app:convergence_Z:12}
    \frac{\bm{X}^{\sigma}_\tau}{\sum_{i\in[I]}N_i^0}\sim \bm{\mathcal{Z}}^h_\tau.
\end{equation}
Together with \eqref{eqn:app:convergence_Z:9}, for any $0<\bar{T}<\infty$ and $\delta>0$, 
\begin{equation}\label{eqn:app:convergence_Z:13}
    \lim_{h\rightarrow \infty}\mathbb{P}\biggl\{\sup_{0\leq \tau\leq \bar{T}}\Bigl\lVert \bm{\mathcal{Z}}^h_{\tau} - \frac{\bar{\bm{x}}_{\tau}}{\sum_{i\in[I]}N_i^0}\Bigr\rVert > \delta\biggr\} = 0,
\end{equation}
where recall that $\bar{\bm{x}}_{\tau}$ is a deterministic trajectory when given $\bar{\bm{x}}_0\in\mathbb{R}^{\mathcal{K}}$. 
Thus, given a random initial point $\bm{\mathcal{Z}}^h_0=\bar{\bm{x}}_0/\sum_{i\in[I]}N^0_i$, for any $0\leq \tau \leq \bar{T}$,
\begin{equation}\label{eqn:app:convergence_Z:14}
    \lim_{h\rightarrow \infty}\mathbb{E}\Bigl[\bm{\mathcal{Z}}^h_{\tau}\Bigr]= \frac{\mathbb{E}\bigl[\bar{\bm{x}}_{\tau}\bigr]}{\sum_{i\in[I]}N^0_i}.
\end{equation}
where $\mathbb{E}\bigl[\bar{\bm{x}}_{\tau}\bigr]$ takes expectation over the initial point $\bar{\bm{x}}_0$.

Similar to \eqref{eqn:app:convergence_Z:7}, for a time point of the real process $t\in[T]_0$, define
\begin{equation}\label{eqn:app:convergence_Z:15}
\tau^h(t)\coloneqq \min\Bigl\{\tau\in\mathbb{R}_0~\Bigl|~\sum_{\kappa\in\bigcup_{t'<t}\mathscr{K}_{t'}}\mathcal{Z}^h_{\kappa,\tau} = 0\Bigr\},
\end{equation}
where $\tau^h(0)=0$.
For $t\in[T]_0$ and such $\tau^h(t)$, based on the definitions of $b^h(\bm{x},\bm{\xi}^{\varphi,h,\varrho}_{\tau})$ and $\bm{\xi}^{\varphi,h,\varrho}_{\tau}$,  \begin{equation}\label{eqn:app:convergence_Z:16}
   \mathcal{Z}^h_{\kappa,\tau^h(t)}\sim Z^{\varphi,h,\varrho}_{\iota_{\kappa}}(t),~\forall \kappa\in[\mathcal{K}].
\end{equation}
Plugging \eqref{eqn:app:convergence_Z:16} and \eqref{eqn:app:convergence_Z:14} in \eqref{eqn:app:convergence_Z:13},  we obtain that \eqref{eqn:theorem:convergence-Z} holds for any $\varphi\in\Phi^1$ and $\varrho^{\varphi,\bar{\bm{\phi}},h}(T)=\varrho\in\mathscr{R}^{T+1}$. It proves the theorem.

\endproof

\section{An Example of  $f^{h,\Lipschitza}_{\kappa,\kappa'}$}\label{app:dirac-delta}

For $\Lipschitza\in(0,1)$ and $u\in[0,1]$, define
\begin{equation}\label{eqn:dirac-delta:1}
y_{\Lipschitza}(u)\coloneqq \begin{cases}
\int_{-\infty}^{\rho(u)}\frac{1}{\Lipschitza\sqrt{\pi}}e^{-(v-\frac{1}{\Lipschitza})^2/\Lipschitza^2} dv,&\text{if }u\in(0,1),\\
0, & \text{if }u=0,\\
1, & \text{if } u=1,
\end{cases}
\end{equation}
where $\rho(u)$ is a real-valued function for $u\in(0,1)$ satisfying $\rho(u)\rightarrow -\infty$ for $u\rightarrow 0$, $\rho(u)\rightarrow +\infty$ for $u\rightarrow 1$, and suitably smooth.
For example, we can specify it as $\rho(u)=-\cot{u\pi}$. 
In this context, $y_{\Lipschitza}(u)$ is continuous for all $u\in(0,1])$, semi-continuous for $u=0,1$, and, for $u\in(0,1)$, the derivative 
\begin{equation}\label{eqn:dirac-delta:2}
    \frac{d y_{\Lipschitza}}{d u} = \frac{1}{\Lipschitza\sqrt{\pi}}e^{-(\cot{u}-\frac{1}{\Lipschitza})^2/\Lipschitza^2} < \frac{1}{\Lipschitza\sqrt{\pi}} < \infty.
\end{equation}
For $h=1$, $\kappa,\kappa'\in[\mathcal{K}]$, $\bm{x}\in\mathbb{R}_0^{\mathcal{K}}$, and $\bm{\xi}\in\mathbb{R}_0^{\mathcal{K}\sum_{\kappa''\in[\mathcal{K}]}N_{i_{\kappa''}}}$, consider
\begin{multline}\label{eqn:dirac-delta:3}
f^{h,\Lipschitza}_{\kappa,\kappa'}(\bm{x},\bm{\xi}) = -\xi_{\kappa,\kappa'}\bigl(\lceil x^-_{\kappa-1}\rceil+1\bigr)y_{\Lipschitza}\bigl(x^-_{\kappa-1} - \lfloor x^-_{\kappa-1}\rfloor\bigr) + \xi_{\kappa,\kappa'}\bigl(\lceil x^-_{\kappa-1} \rceil\bigr)y_{\Lipschitza}\bigl(\lceil x^-_{\kappa-1} \rceil - x^-_{\kappa-1}\bigr)\\
-\xi_{\kappa,\kappa'}\bigl(\lceil x^-_{\kappa} \rceil\bigr)y_{\Lipschitza}\bigl(\lceil x^-_{\kappa}\rceil - x^-_{\kappa}\bigr) + \xi_{\kappa,\kappa'}\bigl(\lceil x^-_{\kappa}\rceil+1\bigr)y_{\Lipschitza}\bigl(x^-_{\kappa} - \lfloor x^-_{\kappa} \rfloor\bigr)\\
-\indicator\{0<x(t_{\kappa}-1)<1\}\sum_{n=\lceil x^-_{\kappa-1} \rceil+1}^{\lceil x^-_{\kappa} \rceil}\xi_{\kappa,\kappa'}(n)y_{\Lipschitza}(\lceil x(t_{\kappa}-1)\rceil - x(t_{\kappa}-1)),
\end{multline}
where  $x(t_{\kappa}-1) = \sum_{\kappa''\in \mathscr{K}_{t_{\kappa}-1}}x_{\kappa''}$ with $x(t)\equiv 0 $ for all $t <0$,
which ensures Lipschitz continuity of $Q^h(\kappa,\kappa',\bm{x},\bm{\xi})$ in $\bm{x}\in\mathbb{R}^{\mathcal{K}}$.
For $h >1$, let $f^{h,\Lipschitza}_{\kappa,\kappa'}(\bm{x},\bm{\xi})\equiv 0$.
In \eqref{eqn:dirac-delta:3}, $f^{h,\Lipschitza}_{\kappa,\kappa'}(\bm{x},\bm{\xi})$ is linear in $\bm{\xi}$, continuous in $\Lipschitza\in(0,1)$, satisfying 
\begin{equation}\label{eqn:dirac-delta:4}
    \lim_{\Lipschitza\rightarrow 0} \bigl(Q^h(\kappa,\kappa',\bm{x},\bm{\xi}) - Q^h(\kappa,\kappa',\lceil \bm{x} \rceil, \bm{\xi})\bigr) = 0,
\end{equation}
for all $\kappa,\kappa'\in[\mathscr{K}]$, $h\in\mathbb{N}_+$, $\bm{x}\in\mathbb{R}_0^{\mathcal{K}}$, and $\bm{\xi}\in\mathbb{R}_0^{\mathcal{K}\sum_{\kappa''\in[\mathcal{K}]}N_{i_{\kappa''}}}$.

We then consider the derivative $d f^{h,\Lipschitza}_{\kappa,\kappa'}/d\Lipschitza$ for $\Lipschitza\in(0,1)$ and $h=1$.
\begin{multline}\label{eqn:dirac-delta:5}
\frac{d f^{h,\Lipschitza}_{\kappa,\kappa'}(\bm{x},\bm{\xi})}{d\Lipschitza} 
=    -\xi_{\kappa,\kappa'}\bigl(\lceil x^-_{\kappa-1}\rceil+1\bigr)\frac{d y_{\Lipschitza}}{d\Lipschitza}\bigl(x^-_{\kappa-1} - \lfloor x^-_{\kappa-1}\rfloor\bigr) + \xi_{\kappa,\kappa'}\bigl(\lceil x^-_{\kappa-1} \rceil\bigr)\frac{d y_{\Lipschitza}}{d\Lipschitza}\bigl(\lceil x^-_{\kappa-1} \rceil - x^-_{\kappa-1}\bigr)\\
-\xi_{\kappa,\kappa'}\bigl(\lceil x^-_{\kappa} \rceil\bigr)\frac{d y_{\Lipschitza}}{d\Lipschitza}\bigl(\lceil x^-_{\kappa}\rceil - x^-_{\kappa}\bigr) + \xi_{\kappa,\kappa'}\bigl(\lceil x^-_{\kappa}\rceil+1\bigr)\frac{d y_{\Lipschitza}}{d\Lipschitza}\bigl(x^-_{\kappa} - \lfloor x^-_{\kappa} \rfloor\bigr)\\
-\indicator\{0<x(t_{\kappa}-1)<1\}\sum_{n=\lceil x^-_{\kappa-1} \rceil+1}^{\lceil x^-_{\kappa} \rceil}\xi_{\kappa,\kappa'}(n)\frac{d y_{\Lipschitza}}{d\Lipschitza}(\lceil x(t_{\kappa}-1)\rceil - x(t_{\kappa}-1)),
\end{multline}
For $u\in(0,1)$,
\begin{equation}\label{eqn:dirac-delta:6}
\frac{d y_{\Lipschitza}(u)}{d\Lipschitza} = \frac{d}{d \Lipschitza} \Phi\Bigl(\frac{\rho(u) - \frac{1}{a}}{\Lipschitza/\sqrt{2}}\Bigr) = \frac{1}{\sqrt{\pi}}\Bigl(\frac{2}{\Lipschitza^3} - \frac{\rho(u)}{\Lipschitza^2}\Bigr)e^{-\bigl(\frac{\rho(u)-\frac{1}{\Lipschitza}}{\Lipschitza}\bigr)^2},
\end{equation}
where $\Phi(x)$ is the cumulative distribution function of the standard normal distribution.
Plugging \eqref{eqn:dirac-delta:6} in \eqref{eqn:dirac-delta:5}, we have $\lim_{\Lipschitza\downarrow 0} df^{h,\Lipschitza}_{\kappa,\kappa'}/d\Lipschitza = 0$.

\section{Proof of Theorem~\ref{theorem:convergence_Z_exp}}
\label{app:theorem:convergence_Z_exp}
Recall that the following proof is based on \cite[Theorem 4.1 in Chapter 7 and Theorem 3.3 in Chapter 3]{freidlin2012random}.
It follows similar lines as the proof for \cite[Theorem 3]{fu2024patrolling} but considers a more general case and
achieves Theorem~\ref{theorem:convergence_Z_exp} that is applicable to a broader range of problems than \cite[Theorem 3]{fu2024patrolling}.

For $U\in\mathbb{R}_0$, $\bm{x},\bm{\omega}\in\mathbb{R}^{\mathcal{K}}$, $\varphi\in\Phi^1$, and $\varrho\in\mathscr{R}^{T+1}$, define
\begin{equation}\label{eqn:large_deviation:define_H}
\mathcal{G}^U(\bm{x},\bm{\omega}) \coloneqq \lim_{\bar{T}\rightarrow \infty}\frac{1}{\bar{T}}\ln \mathbb{E}\exp\Bigl\{\int_0^{\bar{T}}\bigl<\bm{\omega},b^U(\bm{x},\bm{\xi}^{\varphi,1,\varrho}_{\tau}) \bigr>d\tau\Bigr\},
\end{equation}
where $<\cdot,\cdot>$ is the dot production. 
For given $U<\infty$, $\mathcal{G}^U$ is bounded and Lipschitz continuous in both arguments.
Based on \cite[Lemma 4.1 in Chapter 7]{freidlin2012random}, such $\mathcal{G}^U$ is convex in the second argument $\bm{\omega}$.

\begin{lemma}\label{lemma:large_deviation:integral_H}
For $U\in\mathbb{R}_0\cup\{\infty\}$, $\varphi\in\Phi^1$,  $\varrho\in\mathscr{R}^{T+1}$, any compact sets $\mathscr{X}^c,\mathscr{W}^c\subset\mathbb{R}^{\mathcal{K}}$, and any $\bm{x}\in\mathscr{X}^c$ and $\bm{\omega}\in\mathscr{W}^c$, $\mathcal{G}^U(\bm{x},\bm{\omega})$ is Riemann integrable.
\end{lemma}
\proof{Proof of Lemma~\ref{lemma:large_deviation:integral_H}.}
For $\varphi\in\Phi^1$,  $\varrho\in\mathscr{R}^{T+1}$, and any $\bm{\mu}\in\mathbb{R}^{\mathcal{K}\sum_{\kappa\in[\mathcal{K}]}N^0_{i_{\kappa}}}$, define
\begin{equation}
\mathcal{G}_{\xi}(\bm{\mu})\coloneqq \lim_{\bar{T}\rightarrow\infty} \frac{1}{\bar{T}}\ln \mathbb{E} \exp\Bigl\{\int_0^{\bar{T}}\bigl<\bm{\mu},\bm{\xi}^{\varphi,1,\varrho}_{\tau}\bigr>d\tau\Bigr\}.
\end{equation}
Since each element of $\Bigl\lfloor\int_0^{\bar{T}}\bm{\xi}^{\varphi,1,\varrho}_{\tau}d\tau\Bigr\rfloor$ is Poisson distributed, we obtain 
\begin{multline}\label{eqn:lemma:integral_H:1}
\mathbb{E}\exp\Bigl\{\int_0^{\bar{T}}\bigl<\bm{\mu},\bm{\xi}^{\varphi,1,\varrho}_{\tau}\bigr>d\tau\Bigr\}
\leq \mathbb{E}\exp\Bigl\{\bigl<\bm{\mu},\bigl\lfloor\int_0^{\bar{T}}\bm{\xi}^{\varphi,1,\varrho}_{\tau}d\tau\bigr\rfloor+\bm{1}(\bm{\mu})\bigr> \Bigr\}\\
=\exp\Bigl\{\sum_{\begin{subarray}~\kappa,\kappa'\in[\mathcal{K}],\\n\in[N_{i_{\kappa}}^0]\end{subarray}}p^{\varphi,\varrho}_{i_{\kappa},t_{\kappa}}(\iota_{\kappa},\iota_{\kappa'})\bar{T}\bigl(e^{\mu_{\kappa,\kappa'}(n)}-1\bigr)+|\mu_{\kappa,\kappa'}(n)|\Bigr\},
\end{multline}
where $\bm{1}(\bm{\mu})$ is a vector with the $(\kappa,\kappa',n)$th elements equal to $\indicator\{\mu_{\kappa,\kappa'}(n)\geq 0\}-\indicator\{\mu_{\kappa,\kappa'}(n)<0\}$, and recall the parameter of the Poisson distribution $p^{\varphi,\varrho}_{i_{\kappa},t_{\kappa}}(\iota_{\kappa},\iota_{\kappa'})= \frac{1}{\bar{T}}\mathbb{E}\Bigl[\bigl\lfloor\int_0^{\bar{T}}\xi^{\varphi,1,\varrho}_{\kappa,\kappa',\tau}(n)d\tau\bigr\rfloor\Bigr]$.
Hence,
\begin{equation}\label{eqn:lemma:integral_H:2}
    \mathcal{G}_{\xi}(\bm{\mu}) \leq \sum_{\begin{subarray}~\kappa,\kappa'\in[\mathcal{K}],\\n\in[N_{i_{\kappa}}^0]\end{subarray}} p^{\varphi,\varrho}_{i_{\kappa},t_{\kappa}}(\iota_{\kappa},\iota_{\kappa'})\bigl(e^{\mu_{\kappa,\kappa'}(n)}-1\bigr).
\end{equation}

Similarly, replacing $\bm{1}(\bm{\mu})$ with $-\bm{1}(\bm{\mu})$ in \eqref{eqn:lemma:integral_H:1}, we can obtain
\begin{equation}\label{eqn:lemma:integral_H:3}
    \mathcal{G}_{\xi}(\bm{\mu}) \geq \sum_{\begin{subarray}~\kappa,\kappa'\in[\mathcal{K}],\\n\in[N_{i_{\kappa}}^0]\end{subarray}} p^{\varphi,\varrho}_{i_{\kappa},t_{\kappa}}(\iota_{\kappa},\iota_{\kappa'})\bigl(e^{\mu_{\kappa,\kappa'}(n)}-1\bigr).
\end{equation}
Together with \eqref{eqn:lemma:integral_H:2},  
\begin{equation}\label{eqn:lemma:integral_H:4}
 \mathcal{G}_{\xi}(\bm{\mu}) = \sum_{\begin{subarray}~\kappa,\kappa'\in[\mathcal{K}],\\n\in[N_{i_{\kappa}}^0]\end{subarray}} p^{\varphi,\varrho}_{i_{\kappa},t_{\kappa}}(\iota_{\kappa},\iota_{\kappa'})\bigl(e^{\mu_{\kappa,\kappa'}(n)}-1\bigr).
\end{equation}

For any $\bm{v}\in\mathbb{R}^N$, let $\min\{\bm{v},U\} \coloneqq (\min\{v_n,U\}:n\in[N])$.
Recall the definition of $b^U(\bm{x},\bm{\xi})$, for which $b^U(\bm{x},\bm{\xi}) = b^{\infty}(\bm{x},\min\{\bm{\xi},U\}) =  \tilde{\mathcal{Q}}^1(\bm{x})\min\{\bm{\xi},U\}$.
That is, 
\begin{multline}\label{eqn:lemma:integral_H:5}
    \mathcal{G}^U(\bm{x},\bm{\omega}) = \lim_{\bar{T}\rightarrow\infty}\frac{1}{\bar{T}}\ln\mathbb{E}\exp \Bigl\{\bigl<\bm{\omega}^T\tilde{\mathcal{Q}}^1(\bm{x}),\int_0^{\bar{T}}\min\{\bm{\xi}^{\varphi,1,\varrho}_{\tau},U\}d\tau\bigr>\Bigr\}
    \\
    \leq \lim_{\bar{T}\rightarrow\infty}\frac{1}{\bar{T}}\ln\mathbb{E}\exp \Bigl\{\bigl<\bigl(\bm{\omega}^T\tilde{\mathcal{Q}}^1(\bm{x})\bigr)^+,\int_0^{\bar{T}}\bm{\xi}^{\varphi,1,\varrho}_{\tau}d\tau\bigr>\Bigr\}
     = \mathcal{G}_{\xi}\bigl((\bm{\omega}^T\tilde{\mathcal{Q}}^1(\bm{x}))^+\bigr),
\end{multline}
where $(\bm{v})^+\coloneqq (\max\{v_n,0\}:n\in[N])$ for any vector $\bm{v}\in \mathbb{R}^N$.
For any compact sets $\mathscr{W}^c,\mathscr{X}^c\subset \mathbb{R}^{\mathcal{K}}$, $\bm{\omega}\in\mathscr{W}^c$, and $\bm{x}\in\mathscr{X}^c$, $\mathcal{G}^U(\bm{x},\bm{\omega})$ is bounded and is jointly continuous in both arguments.
Hence, it is Riemann integrable. The lemma is proved.

\endproof

From Lemma~\ref{lemma:large_deviation:integral_H}, there exists $\mathcal{G}^U$ defined in \eqref{eqn:large_deviation:define_H} that satisfies
\begin{equation}\label{eqn:theorem:convergence_Z_exp:condition}
    \int_0^{\bar{T}}\mathcal{G}^U(\bm{x}_\tau,\bm{\omega}_\tau)d \tau = \lim_{\epsilon\downarrow 0}\ln \mathbb{E}\exp\Bigl\{\frac{1}{\epsilon}\int_0^{\bar{T}}\bigl<\bm{\omega}_\tau,b^U(\bm{x}_\tau,\bm{\xi}^{\varphi,1,\varrho}_{\tau/\epsilon})\bigr>d\tau\Bigr\}.
\end{equation}

For $\bm{x},\bm{\beta}\in\mathbb{R}^{\mathcal{K}}$, consider the Legendre transform of $\mathcal{G}^U(\bm{x},\bm{\omega})$,
\begin{equation}\label{eqn:large_deviation:legendre_transfor}
\mathcal{L}^U(\bm{x},\bm{\beta}) \coloneq \sup_{\bm{\omega}\in\mathbb{R}^{\mathcal{K}}}\Bigl[<\bm{\omega},\bm{\beta}> - \mathcal{G}^U(\bm{x},\bm{\omega})\Bigr],
\end{equation}
where recall $\mathcal{G}^U(\bm{x},\bm{\omega})$ is convex in $\bm{\omega}$. 
Since $<\bm{0},\bm{\beta}> - \mathcal{G}^U(\bm{x},\bm{0}) = 0$, $\mathcal{L}^U(\bm{x},\bm{\beta})$ is always non-negative.

\begin{lemma}\label{lemma:large_deviation:derivative_H}
Given $\bm{x},\bm{\omega}\in\mathbb{R}^{\mathcal{K}}$, 
\begin{equation}\label{eqn:lemma:derivative_H}
\lim_{U\rightarrow \infty}\frac{\partial \mathcal{G}^U}{\partial \bm{\omega}} = \frac{\partial \mathcal{G}_{\xi}(\bm{\omega}^T\tilde{\mathcal{Q}}^1(\bm{x}))}{\partial \bm{\omega}}.
\end{equation}
\end{lemma}
\proof{Proof of Lemma~\ref{lemma:large_deviation:derivative_H}.}
For $U\in\mathbb{R}_0$, $\varphi\in\Phi^1$,  $\varrho\in\mathscr{R}^{T+1}$, and any $\bm{\mu}\in\mathbb{R}^{\mathcal{K}\sum_{\kappa\in[\mathcal{K}]}N^0_{i_{\kappa}}}$, define
\begin{equation}
\mathcal{G}_{\xi}^U(\bm{\mu})\coloneqq \lim_{\bar{T}\rightarrow\infty} \frac{1}{\bar{T}}\ln \mathbb{E} \exp\Bigl\{\int_0^{\bar{T}}\bigl<\bm{\mu},\min\{\bm{\xi}^{\varphi,1,\varrho}_{\tau},U\}\bigr>d\tau\Bigr\}.
\end{equation}
Similar to the analysis in \eqref{eqn:lemma:integral_H:1}-\eqref{eqn:lemma:integral_H:3}, we have
\begin{multline}\label{eqn:large_deviation:derivative_H:1}
\mathcal{G}_{\xi}^U(\bm{\mu})=\lim_{\bar{T}\rightarrow\infty} \frac{1}{\bar{T}}\ln \mathbb{E} \exp\Bigl\{\bigl<\bm{\mu},\int_0^{\bar{T}}\min\{\bm{\xi}^{\varphi,1,\varrho}_{\tau},U\}d\tau\bigr>\Bigr\}\\ 
\leq \lim_{\bar{T}\rightarrow\infty} \frac{1}{\bar{T}}\ln \mathbb{E} \exp\Bigl\{\bigl<\bm{\mu},\bigl\lfloor\int_0^{\bar{T}}\min\{\bm{\xi}^{\varphi,1,\varrho}_{\tau},U\}d\tau\bigr\rfloor + \bm{1}(\bm{\mu})\bigr>\Bigr\} \\
= \lim_{\bar{T}\rightarrow\infty} \frac{1}{\bar{T}}\ln \mathbb{E} \exp\Bigl\{\bigl<\bm{\mu},\bigl\lfloor\int_0^{\bar{T}}\min\{\bm{\xi}^{\varphi,1,\varrho}_{\tau},U\}d\tau\bigr\rfloor\bigr>\Bigr\},
\end{multline}
and
\begin{multline}\label{eqn:large_deviation:derivative_H:2}
\mathcal{G}_{\xi}^U(\bm{\mu})\geq \lim_{\bar{T}\rightarrow\infty} \frac{1}{\bar{T}}\ln \mathbb{E} \exp\Bigl\{\bigl<\bm{\mu},\bigl\lfloor\int_0^{\bar{T}}\min\{\bm{\xi}^{\varphi,1,\varrho}_{\tau},U\}d\tau\bigr\rfloor - \bm{1}(\bm{\mu})\bigr>\Bigr\} \\
= \lim_{\bar{T}\rightarrow\infty} \frac{1}{\bar{T}}\ln \mathbb{E} \exp\Bigl\{\bigl<\bm{\mu},\bigl\lfloor\int_0^{\bar{T}}\min\{\bm{\xi}^{\varphi,1,\varrho}_{\tau},U\}d\tau\bigr\rfloor\bigr>\Bigr\}.
\end{multline}
That is,
\begin{equation}\label{eqn:large_deviation:derivative_H:4}
\mathcal{G}_{\xi}^U(\bm{\mu})
= \lim_{\bar{T}\rightarrow\infty} \frac{1}{\bar{T}}\ln \mathbb{E} \exp\Bigl\{\bigl<\bm{\mu},\bigl\lfloor\int_0^{\bar{T}}\min\{\bm{\xi}^{\varphi,1,\varrho}_{\tau},U\}d\tau\bigr\rfloor\bigr>\Bigr\}.
\end{equation}
Recall that $\bigl\lfloor\int_0^{\bar{T}}\xi^{\varphi,\varrho}_{\kappa,\kappa',\tau}(n)d\tau\bigr\rfloor $ is Poisson distributed with parameter $p^{\varphi,\varrho}_{i_{\kappa},t_{\kappa}}(\iota_{\kappa},\iota_{\kappa'})$, and given a trajectory $\bigl\{\xi^{\varphi,\varrho}_{\kappa,\kappa',\tau}(n), 0\leq \tau\leq \bar{T}\bigr\}$, there exist a sequence of points $0=\tau_0<\tau_1<\ldots<\tau_M\leq \bar{T}$ such that, for $m\in[M]$,
\[ \int_{\tau_{m-1}}^{\tau_m}\xi^{\varphi,\varrho}_{\kappa,\kappa',\tau}(n)d\tau = 1.\]
With the imposed upper bound $U$, we have
\[
\Bigl\lfloor\int_0^{\bar{T}}\min\{\xi^{\varphi,\varrho}_{\kappa,\kappa',\tau}(n),U\}d\tau\Bigr\rfloor = \sum_{m\in[M]}\theta_m\int_{\tau_{m-1}}^{\tau_m}\xi^{\varphi,\varrho}_{\kappa,\kappa',\tau}(n)d\tau = \sum_{m\in[M]}\theta_m,
\]
where $M$ is Poisson distributed, $\theta_m \coloneqq \max\Bigl\{1,\frac{U}{\xi^{\varphi,\varrho}_{\kappa,\kappa',\tau_{m-1}}(n)}\Bigr\}$, and all such $\theta_m$ ($m\in[M]$ are independently and identically distributed.
It follows that
\begin{multline}
    \mathbb{E}\exp\Bigl\{\bigl<\bm{\mu},\bigl\lfloor\int_0^{\bar{T}}\min\{\bm{\xi}^{\varphi,1,\varrho}_{\tau},U\}d\tau\bigr\rfloor\bigr>\Bigr\}
= \prod_{\kappa,\kappa'\in[\mathcal{K}],n\in[N^0_{i_{\kappa}}]} \exp\Bigl\{\bigl(\Theta_{\kappa,\kappa'}(n)-1\bigr)p^{\varphi,\varrho}_{i_{\kappa},t_{\kappa}}(\iota_{\kappa},\iota_{\kappa'})\bar{T}\Bigr\},
\end{multline}
where 
\[\Theta_{\kappa,\kappa'}(n)\coloneqq \int_0^1 \exp\bigl\{\mu_{\kappa,\kappa'}(n)\theta\bigr\}\bbP\{\theta\}d\theta \leq \exp\{\mu_{\kappa,\kappa'}(n)\}.\]
In particular, $\lim_{U\rightarrow \infty}\bbP\{\theta\} = \delta_1(\theta)$ where $\delta_1(\theta)$ is a Dirac delta function satisfying 
\[ \delta_1(\theta) = \begin{cases}
    0,& \text{if }\theta \neq 1,\\
    \infty, &\text{otherwise},
\end{cases}\]
and $\int_{-\infty}^{\infty}\delta_1(\theta)d\theta = 1$.
Let $\bm{\Theta}\coloneqq (\Theta_{\kappa,\kappa'}(n):\kappa,\kappa'\in[\mathcal{K}],n\in[N^0_{i_{\kappa}}])$.
Obviously, $\lim_{U\rightarrow\infty} \bm{\Theta} = e^{\bm{\mu}}$.

We then obtain
\begin{equation}
    \mathcal{G}^U_{\xi}(\bm{\mu}) = \sum_{\kappa,\kappa'\in[\mathcal{K}],n\in[N^0_{i_{\kappa}}]}\bigl(\Theta_{\kappa,\kappa'}(n)-1\bigr)p^{\varphi,\varrho}_{i_{\kappa},t_{\kappa}}(\iota_{\kappa},\iota_{\kappa'}),
\end{equation}
and, for $\kappa\in[\mathcal{K}]$, 
\begin{equation}
  \frac{\partial}{\omega_{\kappa}} \mathcal{G}^U(\bm{x},\bm{\omega}) = \frac{\partial}{\omega_{\kappa}} \mathcal{G}^U_{\xi}\bigl(\bm{\omega}^T\tilde{\mathcal{Q}}^1(\bm{x})\bigr) = \Bigl<\bm{\rho}_{\kappa},\bm{q}_{\kappa}\Bigr>,
\end{equation}
where $\bm{q}_{\kappa}$ is the $\kappa$th row vector of $\tilde{\mathcal{Q}}^1(\bm{x})$, and 
\[\bm{\rho}_{\kappa} \coloneqq \biggl(p^{\varphi,\varrho}_{i_{\kappa_1},t_{\kappa_1}}(\iota_{\kappa_1},\iota_{\kappa_2})\int_0^1\theta e^{\mu_{\kappa_1,\kappa_2}(n)\theta}\bbP\{\theta\}d\theta:~\kappa_1,\kappa_2\in[\mathcal{K}],n\in[N^0_{i_{\kappa_1}}]\biggr),\]
with $\mu_{\kappa_1,\kappa_2}(n) = \bigl(\bm{\omega}^T\tilde{\mathcal{Q}}^1(\bm{x})\bigr)_{\kappa_1,\kappa_2}(n)$.
Since $\mathcal{G}(\bm{x},\bm{\omega}) = \mathcal{G}_{\xi}(\bm{\omega}^T\tilde{\mathcal{Q}}^1(\bm{x}))$ with $\mathcal{G}_{\xi}$ given by \eqref{eqn:lemma:integral_H:4}, we obtain \eqref{eqn:lemma:derivative_H}.
It proves the lemma.

\endproof

\begin{lemma}\label{lemma:large_deviation:unique_zero}
For $U\in\mathbb{R}_0\cup\{\infty\}$, $\varphi\in\Phi^1$,  $\varrho\in\mathscr{R}^{T+1}$,  $\bm{x},\bm{\beta}\in\mathbb{R}^{\mathcal{K}}$, and any $\delta > 0$, there exists $U_0<\infty$ such that, for all $U > U_0$, if $\lVert \bm{\beta} - \mathbb{E}b^U(\bm{x},\bm{\xi}^{\varphi,1,\varrho}_\tau)\rVert \geq \delta$, then $\mathcal{L}^U(\bm{x},\bm{\beta}) > 0$.
\end{lemma}
\proof{Proof of Lemma~\ref{lemma:large_deviation:unique_zero}.}
Based on \cite[Chapter 7, Section 4]{freidlin2012random}, for any $U\in\mathbb{R}_0$, if $\bm{\beta}=\bar{b}^U(\bm{x})$, then $\mathcal{L}^U(\bm{x},\bm{\beta}) = 0$.

Let $\bm{\omega}^U(\bm{x},\bm{\beta})$ represent the extreme point satisfying $\frac{\partial \mathcal{G}^U}{\partial \bm{\omega}}\Bigr|_{\bm{\omega}=\bm{\omega}^U(\bm{x},\bm{\beta})}=\bm{\beta}$. Such an extreme point may not be unique.
Based on Lemma~\ref{lemma:large_deviation:derivative_H}, when $U\rightarrow \infty$, it becomes
\begin{equation}\label{eqn:lemma:unique_zero:1}
    \bm{\beta} = \lim_{U\rightarrow \infty}\frac{\partial \mathcal{G}^U(\bm{x},\bm{\omega})}{\partial \bm{\omega}}\Bigr|_{\bm{\omega}=\bm{\omega}^U(\bm{x},\bm{\beta})} = \frac{\partial \mathcal{G}_{\xi}\bigl(\bm{\omega}\tilde{\mathcal{Q}}^1(\bm{x})\bigr)}{\partial \bm{\omega}}\Bigr|_{\bm{\omega}=\bm{\omega}^{\infty}(\bm{x},\bm{\beta})},
\end{equation}
where $\bm{\omega}^{\infty}(\bm{x},\bm{\beta}) \coloneqq \lim_{U\rightarrow \infty} \bm{\omega}^U(\bm{x},\bm{\beta})$, and $\mathcal{G}_{\xi}$ is given in \eqref{eqn:lemma:integral_H:4}.
Substituting \eqref{eqn:lemma:integral_H:4} in \eqref{eqn:lemma:unique_zero:1}, it becomes
\begin{equation}\label{eqn:lemma:unique_zero:2}
\bm{\beta} = \tilde{\mathcal{Q}}^1(\bm{x})\Bigl(p^{\varphi,\varrho}_{i_{\kappa},t_{\kappa}}(\iota_{\kappa},\iota_{\kappa'})\exp\bigl\{<\bm{\omega}^{\infty}(\bm{x},\bm{\beta}),\tilde{\bm{q}}_{\kappa,\kappa'}(n,\bm{x})>\bigr\}: \kappa,\kappa'\in[\mathcal{K}],n\in[N^0_{i_{\kappa}}]\Bigr),
\end{equation}
where $\tilde{\bm{q}}_{\kappa,\kappa'}(n,\bm{x})$ are the column vectors of matrix $\tilde{\mathcal{Q}}^1(\bm{x})$.
Moreover, if 
\begin{equation}\label{eqn:lemma:unique_zero:3}
    \lim_{U\rightarrow\infty} \mathcal{L}^U(\bm{x},\bm{\beta}) = <\bm{\omega}^{\infty}(\bm{x},\bm{\beta}),\bm{\beta}>-\mathcal{G}_{\xi}\bigl((\bm{\omega}^{\infty}(\bm{x},\bm{\beta}))^T\tilde{\mathcal{Q}}^1(\bm{x})\bigr)=0,
\end{equation}
then, substituting \eqref{eqn:lemma:unique_zero:2} and \eqref{eqn:lemma:integral_H:4} in \eqref{eqn:lemma:unique_zero:3}, we obtain
\begin{multline}\label{eqn:lemma:unique_zero:4}
<\bm{\omega}^{\infty}(\bm{x},\bm{\beta}),\bm{\beta}> - \mathcal{G}_{\xi}\bigl((\bm{\omega}^{\infty}(\bm{x},\bm{\beta}))^T\tilde{\mathcal{Q}}^1(\bm{x})\bigr)=0\\
=\bigl(\bm{\omega}^{\infty}(\bm{x},\bm{\beta})\bigr)^T\tilde{\mathcal{Q}}^1(\bm{x})\Bigl(p^{\varphi,\varrho}_{i_{\kappa},t_{\kappa}}(\iota_{\kappa},\iota_{\kappa'})\exp\bigl\{\bm{\omega}^{\infty}(\bm{x},\bm{\beta}),\tilde{\bm{q}}_{\kappa,\kappa'}(n,\bm{x})\bigr\}: \kappa,\kappa'\in[\mathcal{K}],n\in[N^0_{i_{\kappa}}]\Bigr)\\- \sum_{\begin{subarray}~\kappa,\kappa'\in[\mathcal{K}],\\n\in[N^0_{i_{\kappa}}]\end{subarray}}p^{\varphi,\varrho}_{i_{\kappa},t_{\kappa}}(\iota_{\kappa},\iota_{\kappa'})\Bigl(\exp\bigl\{<\bm{\omega}^{\infty}(\bm{x},\bm{\beta}),\tilde{\bm{q}}_{\kappa,\kappa'}(n,\bm{x})>\bigr\}-1\Bigr)\\
=\sum_{\begin{subarray}~\kappa,\kappa'\in[\mathcal{K}],\\n\in[N^0_{i_{\kappa}}]\end{subarray}}p^{\varphi,\varrho}_{i_{\kappa},t_{\kappa}}(\iota_{\kappa},\iota_{\kappa'})\Bigl(\exp\bigl\{<\bm{\omega}^{\infty}(\bm{x},\bm{\beta}),\tilde{\bm{q}}_{\kappa,\kappa'}(n,\bm{x})>\bigr\}\Bigl(<\bm{\omega}^{\infty}(\bm{x},\bm{\beta}),\tilde{\bm{q}}_{\kappa,\kappa'}(n,\bm{x})>-1\Bigr)+1\Bigr).
\end{multline}
For \eqref{eqn:lemma:unique_zero:4}, since $\exp\bigl\{<\bm{\omega},\bm{q}>\bigr\}\bigl(<\bm{\omega},\bm{q}>-1\bigr)+1 \geq 0$ for any $\bm{\omega},\bm{q} \in \mathbb{R}^{\mathcal{K}\sum_{\kappa\in[\mathcal{K}]}N^0_{i_{\kappa}}}$ and the equality holds if and only if $<\bm{\omega},\bm{q}>=0$, 
\begin{equation}
\exp\bigl\{<\bm{\omega}^{\infty}(\bm{x},\bm{\beta}),\tilde{\bm{q}}_{\kappa,\kappa'}(n,\bm{x})>\bigr\}\Bigl(<\bm{\omega}^{\infty}(\bm{x},\bm{\beta}),\tilde{\bm{q}}_{\kappa,\kappa'}(n,\bm{x})>-1\Bigr)+1 = 0,    
\end{equation}
for all $\kappa,\kappa'\in[\mathcal{K}]$ and $n\in[N^0_{i_{\kappa}}]$,
or, equivalently,
\begin{equation}\label{eqn:lemma:unique_zero:5}
  \bigl(\bm{\omega}^{\infty}(\bm{x},\bm{\beta})\bigr)^T\tilde{\mathcal{Q}}^1(\bm{x}) = \bm{0}.
\end{equation}
Based on \eqref{eqn:lemma:unique_zero:5} and \eqref{eqn:lemma:unique_zero:2}, given $\bm{x}\in\mathbb{R}^{\mathcal{K}}$,
if $\lim_{U\rightarrow \infty}\mathcal{L}^U(\bm{x},\bm{\beta}) = 0$, then $\bm{\beta}= \tilde{\mathcal{Q}}^1(\bm{x})(p^{\varphi,\varrho}_{i_{\kappa},t_{\kappa}}(\iota_{\kappa},i_{\kappa'}):\kappa,\kappa'\in[\mathcal{K}],n\in[N_{i_{\kappa}}^0])$ which is the unique solution satisfying $\lim_{U\rightarrow \infty}\mathcal{L}^U(\bm{x},\bm{\beta}) = 0$.

Recall that, for all $U\in\mathbb{R}_0$, if $\bm{\beta}=\mathbb{E}b^U(\bm{x},\bm{\xi}^{\varphi,1,\varrho}_{\tau})$, then $\mathcal{L}^U(\bm{x},\bm{\beta}) = 0$.
Since $\mathcal{L}^U(\bm{x},\bm{\beta})$ is continuous in $U$, together with the unique $\bm{\beta}$ for $\lim_{U\rightarrow \infty}\mathcal{L}^U(\bm{x},\bm{\beta}) = 0$, for any given $\bm{x}$ and $\epsilon > 0$, if $\lVert \bm{\beta} - \mathbb{E}b^U(\bm{x},\bm{\xi}^{\varphi,1,\varrho}_{\tau})\rVert \geq \delta$, then there exists $U_0 <\infty$ such that, for all $U > U_0$, $\mathcal{L}^U(\bm{x},\bm{\beta}) > 0$. It proves the lemma.

\endproof

\proof{Proof of Theorem~\ref{theorem:convergence_Z_exp}.}
For $\bm{\chi}_\tau \in \mathbb{R}^{\mathcal{K}}$ and a trajectory $\varpi\coloneqq\{\chi_\tau, 0\leq\tau\leq\bar{T}\}$, define $\Lambda^U_{0,\bar{T}}(\varpi)\coloneqq \int_0^{\bar{T}}\mathcal{L}^U(\bm{\chi}_{\tau},\dot{\bm{\chi}}_{\tau})d\tau$.
Define $\Pi_{0,\bar{T}}$ as the compact set of all such trajectories with given $\bm{\chi}_0=\sum_{i\in[I]}N_i^0\bm{\scrz}_0\in\mathbb{R}^{\mathcal{K}}$, where recall that $\sum_{i\in[I]}N_i^0\bm{\scrz}_0\in\mathbb{R}^{\mathcal{K}}$ is also the given initial state of $\{\bm{X}^{\sigma,U}_\tau,0\leq \tau\leq\bar{T}\}$ and $\{\bar{\bm{x}}^U_\tau,0\leq \tau\leq\bar{T}\}$ that are defined in Appendix~\ref{app:theorem:convergence-Z} and are the solutions of $\dot{\bm{X}}^{\sigma,U}_{\tau}=b^U(\bm{X}^{\sigma}_{\tau},\bm{\xi}^{\varphi,1,\varrho}_{\tau/\sigma})$ and $\dot{\bar{\bm{x}}}^U_{\tau} = \mathbb{E}b^U(\bm{X}^{\sigma}_{\tau},\bm{\xi}^{\varphi,1,\varrho}_{\tau/\sigma})$, respectively.
For $U\in\mathbb{R}_0$ and $\delta > 0$, define a closed set $\mathscr{C}(U,\delta)\coloneqq \{\varpi\in\Pi_{0,\bar{T}}~|~\sup_{0\leq\tau\leq\bar{T}}\lVert \bm{\chi}_{\tau}-\bar{\bm{x}}^U_\tau\rVert \geq \delta\}$.

Based on \cite[Theorem 4.1 in Chapter 7 and Theorem 3.3 in Chapter 3]{freidlin2012random}, for any $U\in\mathbb{R}_0$ and $\delta>0$, as there exists $\mathcal{G}^U$ satisfying \eqref{eqn:theorem:convergence_Z_exp:condition},
\begin{equation}\label{eqn:theorem:convergence_Z_exp:5}
    \limsup_{\sigma\downarrow 0} \sigma \ln \mathbb{P}\Bigl\{\sup_{0\leq \tau\leq \bar{T}}\lVert \bm{X}^{\sigma,U}_\tau - \bar{\bm{x}}^U_{\tau}\rVert > \delta\Bigr\}\leq -\inf_{\varpi\in\mathscr{C}(U,\delta)}\Lambda^U_{0,\bar{T}}(\varpi).
\end{equation}
Based on Lemma~\ref{lemma:large_deviation:unique_zero} and the continuity of $\bar{\bm{x}}^U_\tau$ and $\Lambda^U_{o,\bar{T}}(\varpi)$ in $U\in\mathbb{R}_+\cup\{\infty\}$, for any $\delta >0$, there exist $\sigma>0$ and $U_0<\infty$ such that, for all $U>U_0$, $\inf_{\varpi\in\mathscr{C}(U,\delta)}\Lambda^U_{0,\bar{T}}(\varpi) \geq \sigma > 0$.
Together with \eqref{eqn:theorem:convergence_Z_exp:5}, for any $\delta > 0$ and given initial state $\bm{X}^{\sigma}_0=\bar{\bm{x}}_0=\sum_{i\in[I]}N^0_i \bm{\scrz}_0$, there exist $\sigma_0,C > 0$ such that, for all $\sigma < \sigma_0$,
\begin{equation}
\mathbb{P}\Bigl\{\sup_{0\leq \tau\leq \bar{T}}\lVert \bm{X}^{\sigma}_\tau - \bar{\bm{x}}_{\tau}\rVert >  \delta\Bigr\} \leq e^{-C/\sigma},
\end{equation}
where recall $\bm{X}^{\sigma}_{\tau} = \lim_{U\rightarrow\infty} \bm{X}^{\sigma,U}_{\tau}$ and $\bar{\bm{x}}_{\tau} = \lim_{U\rightarrow\infty}\bar{\bm{x}}^U_{\tau}$.

Along the same lines as the proof of Theorem~\ref{theorem:convergence-Z}, for the solution, $\bm{\mathcal{Z}}^h_\tau$, of \[\dot{\bm{\mathcal{Z}}}^h_{\tau}=\frac{1}{h\sum_{i\in[I]}N^0_i}b^h\bigl(h\sum_{i\in[I]}N^0_i\bm{\mathcal{Z}}^h_{\tau},\bm{\xi}^{\varphi,h,\varrho}_{\tau}\bigr)\] and $\bm{\mathcal{Z}}^h_0 = \bm{\scrz}_0$ and $\sigma = \frac{1}{h}$, we have \eqref{eqn:app:convergence_Z:12}.
Hence, for any $\bar{T} <\infty$ and $\delta>0$, there exist $H <\infty$ and $C>0$ such that, for all $h>H$,
\begin{equation}
  \mathbb{P}\Bigl\{\sup_{0\leq \tau\leq \bar{T}}\lVert \bm{\mathcal{Z}}^{h}_\tau - \lim_{h\rightarrow \infty} \mathbb{E}\bm{\mathcal{Z}}^h_\tau\rVert >  \delta\Bigr\} \leq e^{-Ch}.
\end{equation}
Together with \eqref{eqn:app:convergence_Z:16}, 
we prove \eqref{eqn:theorem:convergence_Z_exp} for any given $\varphi\in\Phi^1$ and $\varrho^{\varphi,\bar{\bm{\phi}},h}(T)=\varrho\in\mathscr{R}^{T+1}$.
It proves the theorem.

\endproof

\section{Proof of Corollary~\ref{coro:converge_transition_prob}}
\label{app:coro:convergence_transition_prob}
This Corollary is proved by invoking Cram\'er's large deviation theorem, as well as Theorem~\ref{theorem:convergence_Z_exp}.
\proof{Proof of Corollary~\ref{coro:converge_transition_prob}.}
For $\varphi\in\Phi^1$, $h\in\mathbb{N}_+$, $t\in[T-1]_0$, $\iota\in[\mathcal{I}]$, and $s'\in\bS_{i_{\iota}}$, given $Z^{\varphi,h}_{\iota,t} = z_{\iota}$ and $N^0\coloneqq \sum_{i\in[I]}N^0_i$, define a random variable
\begin{equation}
Y_{\iota}(s')\coloneqq \lvert \hat{\mathscr{N}}^{\varphi,h}_{i_{\iota},t}(s_{\iota},\actionphi_{\iota},s')\rvert,
\end{equation}
which follows Binomial distribution $B(hN^0 z_{\iota}, p_{i_{\iota}}(s_{\iota},\actionphi_{\iota},s'))$.
Note that, here, $hN^0 z_{\iota}$ must be an integer.
The moment generating function of $Y_{\iota}(s')$ is
\begin{equation}
    \Psi_{\iota,s'}(t) = (1-p_{i_{\iota}}(s_{\iota},\actionphi_{\iota},s')+p_{i_{\iota}}(s_{\iota},\actionphi_{\iota},s')e^t)^{hN^0 z_{\iota}},
\end{equation}
which is strictly convex in $t\in\mathbb{R}$.
We define the Legendre transform of the moment generating function
\begin{equation}
    I_{\iota,s'}(y) \coloneqq \sup_{t\in\mathbb{R}}(yt - \Psi_{\iota,s'}(t)),
\end{equation}
for $y\in\mathbb{R}$.
Such $I_{\iota,s'}(y)$ is always non-negative and also strictly convex (in $y\in\mathbb{R}$) with the only zero point $y=\mathbb{E}Y_{\iota}(s) = hN^0 z_{\iota} p_{i_{\iota}}(s_{\iota},\actionphi_{\iota},s')$.

For any $\delta>0$, define a compact set 
\begin{equation}
\mathscr{C}_{\iota,s'}(\delta)\coloneqq \Bigl\{y\in \mathbb{R}~\Bigl|~\lvert y-\mathbb{E}Y_{\iota}(s')\rvert \geq \delta\Bigr\},
\end{equation}
Based on Cram\'er's large deviation theorem, for any $\delta>0$ and $\mathscr{C}_{\iota,s'}(\delta)$,
\begin{equation}\label{eqn:coro:transition_prob:5}
    \limsup_{M\rightarrow \infty}\frac{1}{M}\ln \mathbb{P}\Bigl\{\frac{\sum_{m\in[M]}Y^m_{\iota}(s')}{M} \in \mathscr{C}_{\iota,s'}(\delta)\Bigr\}\leq -\inf_{y\in\mathscr{C}_{\iota,s'}(\delta)}I_{\iota,s'}(y),
\end{equation}
where $Y^m_{\iota}(s')\sim Y_{\iota}(s')$ are $M$ independently and identically distributed variables.
Recall that $I_{\iota,s'}(y)$ is strictly convex with the unique zero point $y=\mathbb{E}Y_{\iota}(s')$ that is excluded from $\mathscr{C}_{\iota,s'}(\delta)$. 
That is, from \eqref{eqn:coro:transition_prob:5}, for any $\delta>0$, there exists $C>0$ such that
\begin{equation}\label{eqn:coro:transition_prob:6}
    \limsup_{\bar{h}\rightarrow\infty}\frac{1}{\bar{h}}\ln\mathbb{P}\Bigl\{\bigl\lvert\frac{\bar{Y}^{\bar{h}}_{\iota}(s')}{\bar{h}/h}-\mathbb{E}Y_{\iota}(s')\bigr\rvert \geq \delta\Bigr\} \leq -C,
\end{equation}
where $\bar{h}\coloneqq Mh$, $\bar{Y}^{\bar{h}}_{\iota}(s')\coloneqq \sum_{m\in[\bar{h}/h]}Y^m_{\iota}(s')\sim B(\bar{h}N^0z_{\iota},p_{i_{\iota}}(s_{\iota},\actionphi_{\iota},s'))$ and is in distribution equivalent to $|\hat{\mathscr{N}}^{\varphi,\bar{h}}_{i_{\iota},t}(s_{\iota},\actionphi_{\iota},s')|$ with given $Z^{\varphi,h}_{\iota}(t) = z_{\iota}$.
From \eqref{eqn:coro:transition_prob:6}, given $Z^{\varphi,\bar{h}}_{\iota}(t) = z_{\iota}$, for any $\delta>0$, there exists $C>0$ such that
\begin{multline}\label{eqn:coro:transition_prob:7}
 \lim_{\bar{h}\rightarrow\infty}\frac{1}{\bar{h}}\ln\mathbb{P}\biggl\{\Bigl\lvert\frac{|\hat{\mathscr{N}}^{\varphi,\bar{h}}_{i_{\iota},t}(s_{\iota},\actionphi_{\iota},s')|}{|\hat{\mathscr{N}}^{\varphi,\bar{h}}_{i_{\iota},t}(s_{\iota},\actionphi_{\iota})|}-p_{i_{\iota}}(s_{\iota},\actionphi_{\iota},s')\Bigr\rvert \geq \delta\biggr\} \\=    \lim_{\bar{h}\rightarrow\infty}\frac{1}{\bar{h}}\ln\mathbb{P}\Bigl\{\bigl\lvert\frac{\bar{Y}^{\bar{h}}_{\iota}(s')}{\bar{h}/h}-\mathbb{E}Y_{\iota}(s')\bigr\rvert \geq hN^0z_{\iota}\delta\Bigr\} \leq -C.
\end{multline}

For any $\delta_1,\delta_2 >0$, if $\lim_{h\rightarrow\infty}z^{\varphi,h}_{\iota}(t) > 0$, then there exist $C_1,C_2 > 0$ and $H < \infty$ such that, for all $\bar{h}>H$,
\begin{multline}\label{eqn:coro:transition_prob:8}
\mathbb{P}\biggl\{Z^{\varphi,\bar{h}}_{\iota}(t)\Bigl\lvert\frac{|\hat{\mathscr{N}}^{\varphi,\bar{h}}_{i_{\iota},t}(s_{\iota},\actionphi_{\iota},s')|}{|\hat{\mathscr{N}}^{\varphi,\bar{h}}_{i_{\iota},t}(s_{\iota},\actionphi_{\iota})|}-p_{i_{\iota}}(s_{\iota},\actionphi_{\iota},s')\Bigr\rvert \geq \delta_1\biggr\}\\
=\mathbb{P}\biggl\{Z^{\varphi,\bar{h}}_{\iota}(t)\Bigl\lvert\frac{|\hat{\mathscr{N}}^{\varphi,\bar{h}}_{i_{\iota},t}(s_{\iota},\actionphi_{\iota},s')|}{|\hat{\mathscr{N}}^{\varphi,\bar{h}}_{i_{\iota},t}(s_{\iota},\actionphi_{\iota})|}-p_{i_{\iota}}(s_{\iota},\actionphi_{\iota},s')\Bigr\rvert \geq \delta_1~\biggl|~Z^{\varphi,\bar{h}}_{\iota,t} \geq \delta_2\biggr\} \mathbb{P}\Bigl\{Z^{\varphi,\bar{h}}_{\iota,t} \geq \delta_2\Bigr\}\\
+\mathbb{P}\biggl\{Z^{\varphi,\bar{h}}_{\iota}(t)\Bigl\lvert\frac{|\hat{\mathscr{N}}^{\varphi,\bar{h}}_{i_{\iota},t}(s_{\iota},\actionphi_{\iota},s')|}{|\hat{\mathscr{N}}^{\varphi,\bar{h}}_{i_{\iota},t}(s_{\iota},\actionphi_{\iota})|}-p_{i_{\iota}}(s_{\iota},\actionphi_{\iota},s')\Bigr\rvert \geq \delta_1~\biggl|~Z^{\varphi,\bar{h}}_{\iota,t} < \delta_2\biggr\}\mathbb{P}\Bigl\{Z^{\varphi,\bar{h}}_{\iota,t} < \delta_2\Bigr\}\\
\stackrel{(a)}{\leq}\mathbb{P}\biggl\{Z^{\varphi,\bar{h}}_{\iota}(t)\Bigl\lvert\frac{|\hat{\mathscr{N}}^{\varphi,\bar{h}}_{i_{\iota},t}(s_{\iota},\actionphi_{\iota},s')|}{|\hat{\mathscr{N}}^{\varphi,\bar{h}}_{i_{\iota},t}(s_{\iota},\actionphi_{\iota})|}-p_{i_{\iota}}(s_{\iota},\actionphi_{\iota},s')\Bigr\rvert \geq \delta_1~\biggl|~Z^{\varphi,\bar{h}}_{\iota,t} \geq \delta_2\biggr\} 
+e^{-C_2 \bar{h}}\\
\stackrel{(b)}{\leq}e^{-C_1 \bar{h}}+e^{-C_2 \bar{h}},
\end{multline}
where Inequality (a) comes from Theorem~\ref{theorem:convergence_Z_exp}, and Inequality (b) is based on \eqref{eqn:coro:transition_prob:7}.

From Theorem~\ref{theorem:convergence_Z_exp}, for $\delta>0$, if $\lim_{h\rightarrow\infty}z^{\varphi,h}_{\iota}(t) = 0$, then there exist $C>0$ and $H<\infty$ such that, for all $h>H$,
\begin{equation}\label{eqn:coro:transition_prob:9}
\mathbb{P}\biggl\{Z^{\varphi,h}_{\iota}(t)\Bigl\lvert\frac{|\hat{\mathscr{N}}^{\varphi,h}_{i_{\iota},t}(s_{\iota},\actionphi_{\iota},s')|}{|\hat{\mathscr{N}}^{\varphi,h}_{i_{\iota},t}(s_{\iota},\actionphi_{\iota})|}-p_{i_{\iota}}(s_{\iota},\actionphi_{\iota},s')\Bigr\rvert \geq \delta\biggr\}\\
\leq \mathbb{P}\bigl\{Z^{\varphi,h}_{\iota}(t)\geq \delta\bigr\} \leq e^{-Ch}.
\end{equation}
We have proven \eqref{eqn:convergence_transition_prob:1}.


Equality~\eqref{eqn:convergence_transition_prob:2} is led by \eqref{eqn:convergence_transition_prob:1} and the law of large numbers.
More precisely, for $\varphi\in\Phi^1$, $k\in[K]$, $\bar{\phi}_k\in\oPhiLocal$, $t\in[T]_0$, and $\iota\in[\mathcal{I}]$,
\begin{multline}\label{eqn:coro:transition_prob:10}
    Z^{\varphi,h}_{\iota}(t)\lvert w^{k,\bar{\phi}_k}_{i_{\iota},t}(s_{\iota},\actionphi_{\iota},\bm{Q})\rvert
    \leq Z^{\varphi,h}_{\iota}(t)\Bigl\lvert \frac{\sum_{n\in\hat{\mathscr{N}}^{\varphi,h}_{i_{\iota},t}(s_{\iota},\actionphi_{\iota})}\bar{R}^{k}_{i_{\iota},n}(s_{\iota},\actionphi_{\iota})}{|\hat{\mathscr{N}}^{\varphi,h}_{i_{\iota},t}(s_{\iota},\actionphi_{\iota})|}-\bar{r}^k_{i_{\iota}}(s_{\iota},\actionphi_{\iota})\Bigr\rvert \\+ Z^{\varphi,h}_{\iota}(t)\sum_{s'\in\bS_{i_{\iota}}\backslash\{s_0\}}\biggl\lvert\frac{|\hat{\mathscr{N}}^{\varphi,h}_{i_{\iota},t}(s_{\iota},\actionphi_{\iota},s')|}{|\hat{\mathscr{N}}^{\varphi,h}_{i_{\iota},t}(s_{\iota},\actionphi_{\iota})|}-p_{i_{\iota}}(s_{\iota},\actionphi_{\iota},s') \biggr\rvert \biggl\lvert \functionQ_{i_{\iota}}(s',\actionbarphi_{k,i_{\iota}}(s'))\biggr\rvert,
\end{multline}
where, as $h\rightarrow\infty$, because of the law of large numbers and \eqref{eqn:convergence_transition_prob:1}, the first and the second terms at the right hand side diminish to zero. It proves \eqref{eqn:convergence_transition_prob:2}.

If \eqref{eqn:assumption:R_exp} holds, then for any $\epsilon>0$, there exist $C_1,C_2>0$ and $H<\infty$ such that, for all $h>H$,
\begin{multline}
\mathbb{P}\Bigl\{Z^{\varphi,h}_{\iota}(t)\lvert w^{k,\bar{\phi}_k}_{i_{\iota},t}(s_{\iota},\actionphi_{\iota},\bm{Q})\rvert > \epsilon\Bigr\}
\leq \mathbb{P}\biggl\{Z^{\varphi,h}_{\iota}(t)\Bigl\lvert \frac{\sum_{n\in\hat{\mathscr{N}}^{\varphi,h}_{i_{\iota},t}(s_{\iota},\actionphi_{\iota})}\bar{R}^{k}_{i_{\iota},n}(s_{\iota},\actionphi_{\iota})}{|\hat{\mathscr{N}}^{\varphi,h}_{i_{\iota},t}(s_{\iota},\actionphi_{\iota})|}-\bar{r}^k_{i_{\iota}}(s_{\iota},\actionphi_{\iota})\biggr\rvert > \epsilon/2\Bigr\}\\
+ \mathbb{P}\biggl\{Z^{\varphi,h}_{\iota}(t)\sum_{s'\in\bS_{i_{\iota}}\backslash\{s_0\}}\biggl\lvert\frac{|\hat{\mathscr{N}}^{\varphi,h}_{i_{\iota},t}(s_{\iota},\actionphi_{\iota},s')|}{|\hat{\mathscr{N}}^{\varphi,h}_{i_{\iota},t}(s_{\iota},\actionphi_{\iota})|}-p_{i_{\iota}}(s_{\iota},\actionphi_{\iota},s') \biggr\rvert \biggl\lvert \functionQ_{i_{\iota}}(s',\actionbarphi_{k,i_{\iota}}(s'))\biggr\rvert > \epsilon/2\biggr\}
\leq e^{-C_1 h} + e^{-C_2 h},
\end{multline}
where the last inequality comes from \eqref{eqn:assumption:R_exp} and \eqref{eqn:convergence_transition_prob:1}.
It proves \eqref{eqn:convergence_transition_prob:3}.

\endproof

\section{Proof of Proposition~\ref{prop:stimulating_algo:convergence_span}}
\label{app:stimulating_algo:convergence_span}

We explain here that the ergodic state $s_0$ can be interpreted as a \emph{termination state}.

For each $i\in[I]$, we define a virtual, dummy state $\bar{s}_0$ and its associated transition kernels $\tilde{\mathcal{P}}_i(a) = \bigl[\tilde{p}_i(s,a,s')\bigr]_{(|\bS_i|+1)\times(|\bS_i|+1)}$ ($a\in\bA_i$), for which the transition probabilities  $\tilde{p}_i(s,a,\bar{s}_0) = \hat{p}_i(s,a,s_0)$, $\tilde{p}_i(s,a,s')=\hat{p}_i(s,a,s')$ and $\tilde{p}_i(s,a,s_0) = 0$ for all $a\in\bA_i$, $s\in\bS_i$ and $s'\in\bS_i\backslash\{s_0\}$, $\tilde{p}_i(\bar{s}_0,a,s') = 0$, and $\tilde{p}_i(\bar{s}_0,a,\bar{s}_0)=1$ for all $(a,s')\in\bA_i\times\bS_i$.
The state $\bar{s}_0$ is a duplicated version of $s_0$ but used to distinguish the first and the second entrance to $s_0$ of the underlying MDP.
With respect to transition kernels $\tilde{\mathcal{P}}_i(a)$, we are only interested in the process that starts from any state $s\in\bS_i\backslash\{s_0\}$ until it enters state $s_0$ (or, equivalently, $\bar{s}_0$); and, if the process starts in state $s_0$, then we are also interested in the process until it re-enters $s_0$ (that is, $\bar{s}_0$) again.
These periods of the underlying MDP process are suffix for obtaining the Q factors.
The operation $\hat{\mathcal{T}}^{k,\phi}_i$ defined in \eqref{eqn:define_hat_T} can be rewritten as 
\begin{equation}\label{eqn:re_define_hat_T}
(\hat{\mathcal{T}}^{k,\phi}_i\functionQ)(s,a) = \hat{r}^k_i(s,a) + \sum_{s'\in\{\bar{s}_0\}\cup\bS_i}\tilde{p}_i(s,a,s')Q(s',\actionphi_i(s')),
\end{equation}
with $Q(\bar{s}_0,a)\equiv 0$ for all $a\in\bA_i$.
In this context, the virtual state $\bar{s}_0$ is a termination state (the process runs into $\bar{s}_0$ and never gets out).
For any $i\in[I], m\in[M]$, replacing \eqref{eqn:define_hat_T} with \eqref{eqn:re_define_hat_T} will not change the operation $\hat{\mathcal{T}}^{k,\phi}_i$ on all the states $s\in\bS_i$, nor the values of $\tilde{\functionQ}^k_{i,m}$ in Step~\ref{step:11}.
Similarly, we can attach the termination state $\bar{s}_0$ for the operation $\mathcal{T}^{k,\phi}_i$ described in \eqref{eqn:define_H_phi}.

\begin{condition}{Good Case}\label{condition:good_case}
Given $h\in\mathbb{N}_+$, initial states $\bm{z}_0\in\Delta_{[\mathcal{I}]}$, $\bm{Q}_0\in\prod_{i\in[I]}\bF(\bS_i\times\bA_i)$, secondary policies $\bar{\bm{\phi}}\in(\oPhiLocal)^K$, primary policy $\varphi\in \Phi^2(\bm{z}_0,\bm{Q}_0,\bar{\bm{\phi}})$, and any SA pairs $\iota,\iota'\in[\mathcal{I}]$, if $p_{i_{\iota}}(s_{\iota},\actionphi_{\iota},s_{\iota'}) >0$, then, by time $T^*$ and after Step~\ref{step:9}, the estimate probability $\hat{p}_{i_{\iota}}(s_{\iota},\actionphi_{\iota},s_{\iota'}) > 0$.
\end{condition}

\begin{lemma}\label{lemma:good_case}
Given $h\in\mathbb{N}_+$, initial states $\bm{z}_0\in\Delta_{]\mathcal{I}]}$, $\bm{Q}_0\in\prod_{i\in[I]}\bF(\bS_i\times\bA_i)$, secondary policies $\bar{\bm{\phi}}\in(\oPhiLocal)^K$, and primary policy $\varphi\in \Phi^2(\bm{z}_0,\bm{Q}_0,\bar{\bm{\phi}})$,
there exists a constant $C>0$ such that the WCG-learning system falls into the \partialref{condition:good_case}{good case} is with probability at least $1-e^{-Ch}$.
\end{lemma}
\proof{Proof of Lemma~\ref{lemma:good_case}.}
Based on Corollary~\ref{coro:converge_transition_prob}, since the employed primary policy $\varphi$ is in $\Phi^2(\bm{z}_0,\bm{Q}_0,\bar{\bm{\phi}})$, for $k\in[K]$, $i\in[I]$, and any $\epsilon>0$, there exist $H<\infty$ and $C>0$ such that, for all $h>H$, $\iota\in[\mathcal{I}]$ and $s'\in\bS_i$, if $p_{i_{\iota}}(s_{\iota},\actionphi_{\iota},s')>0$, then
\begin{equation}\label{eqn:stimulating_algo:convergence_span:1}
    \mathbb{P}\Bigl\{\hat{p}_{i_{\iota}}(s_{\iota},\actionphi_{\iota},s')<p_{i_{\iota}}(s_{\iota},\actionphi_{\iota},s')-\epsilon \Bigr\} < e^{-Ch}.
\end{equation}
In other words, when $p_i(s_{\iota},\actionphi_{\iota},s') >0$, we can always take sufficiently small $\epsilon$, such that $\hat{p}_{i_{\iota}}(s_{\iota},\actionphi_{\iota},s')$ in \eqref{eqn:stimulating_algo:convergence_span:1} is ensured to be positive with probability $1-e^{-Ch}$.
The probability of having the \partialref{condition:good_case}{good case} is at least $1-e^{-Ch}$ for some constant $C>0$.
It proves the lemma.

\endproof

\proof{Proof of Proposition~\ref{prop:stimulating_algo:convergence_span}.}
For $k\in[K]$ and $i\in[I]$, in the \partialref{condition:good_case}{good case}, since $\bar{\phi}_k\in\oPhiLocal$ (satisfying \partialref{condition:ergodic}{ergodic condition}), there exists $M<\infty$ such that, for any initial probability distribution $\bm{\pi}\in\Delta_{[\mathcal{I}]}$ of all the SA pairs, 
\begin{equation}\label{eqn:stimulating_algo:convergence_span:2}
    \Bigl(\bigl(\hat{\mathcal{T}}^{k,\bar{\phi}_k}_i\bigr)^M\bm{\pi}\Bigr)\bigl(\bar{s}_0,\actionbarphi_{k,i}(s_0)\bigr) > 0.
\end{equation}
For such $M<\infty$, consider
\begin{multline}\label{eqn:stimulating_algo:convergence_span:3}
    \Delta^{k,M}_i(\bar{\phi}_k)\coloneqq \min_{(\iota_1,\iota_2)\in[\mathcal{I}]^2}\sum_{\iota_3\in[\mathcal{I}]}\min\biggl\{\Bigl(\bigl(\hat{\mathcal{T}}^{k,\bar{\phi}_k}_i\bigr)^M\bm{\pi}(\iota_1)\Bigr)(s_{\iota_3},\actionphi_{\iota_3}),\bigl(\hat{\mathcal{T}}^{k,\bar{\phi}_k}_i\bigr)^M\bm{\pi}(\iota_2)\Bigr)(s_{\iota_3},\actionphi_{\iota_3})\biggr\}\\
    \geq \min_{(\iota_1,\iota_2)\in[\mathcal{I}]^2}\min\biggl\{\Bigl(\bigl(\hat{\mathcal{T}}^{k,\bar{\phi}_k}_i\bigr)^M\bm{\pi}(\iota_1)\Bigr)(\bar{s}_0,\actionbarphi_{k,i}(s_0)),\Bigl(\bigl(\hat{\mathcal{T}}^{k,\bar{\phi}_k}_i\bigr)^M\bm{\pi}(\iota_2)\Bigr)(\bar{s}_0,\actionbarphi_{k,i}(s_0))\biggr\} > 0,
\end{multline}
where $\bm{\pi}(\iota)$, for any $\iota\in[\mathcal{I}]$, is a vector with all element zero but $\pi_{\iota}=1$, and the last inequality is based on \eqref{eqn:stimulating_algo:convergence_span:2}. 
From \cite[Proposition 6.6.1]{puterman2005markov}, we obtain that, for any $k\in[K]$ and $i\in[I]$, if it is in the \partialref{condition:good_case}{good case}, then, for any $\functionQ\in\bF\bigl((\bS_i\cup\{\bar{s}_0)\times \bA_i\bigr)$,
\begin{equation}\label{eqn:stimulating_algo:convergence_span:4}
    sp\Bigl(\bigl(\hat{\mathcal{T}}^{k,\bar{\phi}_k}_i\bigr)^M\functionQ\Bigr)\leq (1-\Delta^{k,M}_i(\bar{\phi}_k))sp (\functionQ),
\end{equation}
where $\functionQ\in\bF\bigl((\bS_i\cup\{\bar{s}_0)\times \bA_i\bigr)$ includes the SA pairs $(\bar{s}_0,a)$ for $a\in\bA_i$.
Let $\tilde{\functionQ}^{k,+}_{i,m}\in\bF\bigl((\bS_i\cup\{\bar{s}_0)\times \bA_i\bigr)$ represent the Q factors for all $(s,a)\in \bigl(\bS_i\cup\{\bar{s}_0\}\bigr)\times \bA_i$, for which $\tilde{\functionQ}^k_{i,m} = \tilde{\functionQ}^{k,+}_{i,m}(s,a)$ for all $(s,a)\in\bS_i\times \bA_i$.
It follows that
\begin{equation}\label{eqn:stimulating_algo:convergence_span:5}
    sp(\tilde{\functionQ}^{k,+}_{i,M+1} - \tilde{\functionQ}^{k,+}_{i,M}) = sp\Bigl(\bigl(\hat{\mathcal{T}}^{k,\bar{\phi}_k}_i\bigr)^M(\tilde{\functionQ}^{k,+}_{i,1}-\tilde{\functionQ}^{k,+}_{i,0})\Bigr)\leq (1-\Delta^{k,M}_i(\bar{\phi}_k))(\tilde{\functionQ}^{k,+}_{i,1}-\tilde{\functionQ}^{k,+}_{i,0}),
\end{equation}
where $1-\Delta^{k,M}_i(\bar{\phi}_k)<1$.
Together with Lemma~\ref{lemma:good_case}, we proves the proposition.

\endproof

\section{Proof of Proposition~\ref{prop:stimulating_algo:convergence_norm}}
\label{app:stimulating_algo:convergence_norm}

\begin{lemma}\label{lemma:value_interation:convergence}
For $\bm{z}_0\in[0,1]^{\mathcal{I}}$,  $\bm{Q}_0\in\prod_{i\in[I]}\bF(\bS_i\times\bA_i)$, secondary polices $\bar{\bm{\phi}}\in\bigl(\oPhiLocal\bigr)^K$,  and primary policy $\varphi\in\Phi^2(\bm{z}_0,\bm{Q}_0,\bar{\bm{\phi}})$,
if the WCG-learning system is in the \partialref{condition:good_case}{good case}, then
the limit $\tilde{\functionQ}^{k,+}_{i,\infty}\coloneqq \lim_{m\rightarrow \infty} \tilde{\functionQ}^{k,+}_{i,m}$ exists and satisfies $\tilde{\functionQ}^{k,+}_{i,\infty} = \hat{\mathcal{T}}^{k,\bar{\phi}_k}_i\tilde{\functionQ}^{k,+}_{i,\infty}$.
\end{lemma}
\proof{Proof of Lemma~\ref{lemma:value_interation:convergence}.}
For each $(k,i)\in[K]\times[I]$, we consider the value iteration operation $\hat{\mathcal{T}}^{k,\bar{\phi}_k}_i$ in the way that described in \eqref{eqn:re_define_hat_T} with the additional virtual state $\bar{s}_0$.
For $k\in[K]$ and $i\in[I]$, in the \partialref{condition:good_case}{good case}, based on \eqref{eqn:stimulating_algo:convergence_span:5}, 
we have $\lim_{m\rightarrow \infty} sp(\tilde{\functionQ}^{k,+}_{i,m+1}-\tilde{\functionQ}^{k,+}_{i,m}) = 0$; that is,
\begin{equation}\label{eqn:convergence:value_iteration:1}
\lim_{m\rightarrow \infty} \bigl(\tilde{\functionQ}^{k,+}_{i,m+1}-\tilde{\functionQ}^{k,+}_{i,m}\bigr) = x \bm{1},
\end{equation}
where $x\in\mathbb{R}$ and $\bm{1}$ is a vector with all elements equal to 1. 
Since $\tilde{Q}^{k,+}_{i,m}(\bar{s}_0,a)\equiv 0$ for all $m\in\mathbb{N}_+$, we have $x=0$ for \eqref{eqn:convergence:value_iteration:1}.
That is, in the \partialref{condition:good_case}{good case}, the limit $\tilde{\functionQ}^{k,+}_{i,\infty}\coloneqq \lim_{m\rightarrow \infty} \tilde{\functionQ}^{k,+}_{i,m}$ exists and hence, $\tilde{\functionQ}^{k,+}_{i,\infty} = \hat{\mathcal{T}}^{k,\bar{\phi}_k}_i\tilde{\functionQ}^{k,+}_{i,\infty}$.
It proves the lemma.
\endproof

\begin{lemma}\label{lemma:convergence:value_iteration}
For $\bm{z}_0\in\Delta_{[\mathcal{I}]}$,  $\bm{Q}_0\in\prod_{i\in[I]}\bF(\bS_i\times\bA_i)$, secondary polices $\bar{\bm{\phi}}\in\bigl(\oPhiLocal\bigr)^K$,  and primary policy $\varphi\in\Phi^2(\bm{z}_0,\bm{Q}_0,\bar{\bm{\phi}})$,
if the WCG-learning system is in the \partialref{condition:good_case}{good case}, then
there exist $C,H<\infty$ such that, for all $h>H$, $k\in[K]$, $i\in[I]$, and any $\epsilon > 0$ used for the stopping condition in Step~\ref{step:11},
\begin{equation}\label{eqn:convergence:value_iteration}
\bigl\lVert \tilde{\functionQ}^{k,+}_i - \tilde{\functionQ}^{k,+}_{i,\infty}\bigr\rVert\leq C\epsilon,
\end{equation}
where recall $\tilde{\functionQ}^k_i$ is the output of the value iteration in \rm{Step~\ref{step:11}}, and $\tilde{\functionQ}^{k,+}_{i,\infty}\coloneqq \lim_{m\rightarrow \infty}\tilde{\functionQ}^{k,+}_{i,m}$.
\end{lemma}
\proof{Proof of Lemma~\ref{lemma:convergence:value_iteration}.}
In the \partialref{condition:good_case}{good case}, for $m>M$ with $M<\infty$  the smallest integer satisfying \eqref{eqn:stimulating_algo:convergence_span:2},
\begin{multline}\label{eqn:convergence:value_iteration:2}
sp(\tilde{\functionQ}^{k,+}_{i,m+1} - \tilde{\functionQ}^{k,+}_{i,\infty})\leq \sum_{u=0}^{\infty}sp\bigl(\tilde{\functionQ}^{k,+}_{i,m+u+1} - \tilde{\functionQ}^{k,+}_{i,m+u+2}\bigr)\\
\leq \sum_{u=0}^{\infty}(1-\Delta^{k,M}_i(\bar{\phi}_k))^u M\max_{m'=m+1}^{m+M} sp\bigl(\tilde{\functionQ}^{k,+}_{i,m'} - \tilde{\functionQ}^{k,+}_{i,m'+1}\bigr)\\
=\frac{M}{\Delta^{k,M}_i(\bar{\phi}_k)}\max_{m'=m+1}^{m+M} sp\bigl(\tilde{\functionQ}^{k,+}_{i,m'} - \tilde{\functionQ}^{k,+}_{i,m'+1}\bigr),
\end{multline}
where the second inequality is based on \eqref{eqn:stimulating_algo:convergence_span:4}.
If the stopping condition in Step~\ref{step:11} is satisfied with $sp(\tilde{\functionQ}^k_{i,m+1} - \tilde{\functionQ}_{i,m})< \epsilon$, then, together with \eqref{eqn:stimulating_algo:convergence_span:4} and the boundness of $\tilde{\functionQ}^{k,+}_{i,0}$, there exists $B < \infty$ (independent from $\epsilon$) such that $\max_{m' = m+1}^{m+M}sp(\tilde{\functionQ}^{k,+}_{i,m'} - \tilde{\functionQ}^{k,+}_{i,m'+1}) \leq B\epsilon$.
For such $m,M,B$ and $\epsilon$, from \eqref{eqn:convergence:value_iteration:2}, we obtain that 
\begin{equation}\label{convergence:value_iteration:3}
sp(\tilde{\functionQ}^{k,+}_{i,m+1} - \tilde{\functionQ}^{k,+}_{i,\infty}) \leq \frac{MB}{\Delta^{k,M}_i(\bar{\phi}_k)}\epsilon.
\end{equation}
Since $\tilde{Q}^{k,+}_{i,m+1} (\bar{s}_0,a)= \tilde{Q}^{k,+}_{i,\infty}(\bar{s}_0,a) = 0 $, we have $\lVert\tilde{\functionQ}^{k,+}_{i,m+1} - \tilde{\functionQ}^{k,+}_{i,\infty}\rVert\leq \frac{MB\sqrt{|\bS_i||\bA_i|}}{\Delta^{k,M}_i(\bar{\phi}_k)}\epsilon $.
It proves the lemma.

\endproof

\begin{lemma}\label{lemma:stimulating_algo:convergence_norm}
For $\bm{z}_0\in\Delta_{[\mathcal{I}]}$,  $\bm{Q}_0\in\prod_{i\in[I]}\bF(\bS_i\times\bA_i)$, secondary polices $\bar{\bm{\phi}}\in\bigl(\oPhiLocal\bigr)^K$,  and primary policy $\varphi\in\Phi^2(\bm{z}_0,\bm{Q}_0,\bar{\bm{\phi}})$,
if the WCG-learning system is in the \partialref{condition:good_case}{good case}, then 
there exist $C_1,C_2 < \infty$ such that, for $k\in[K]$, $i\in[I]$, and any $\epsilon$ used for the stopping condition of the value iteration in Step~\ref{step:11},
\begin{equation}\label{eqn:stimulating_algo:convergence_norm}
\lVert \tilde{\functionQ}^k_i - \functionQ^{k,\bar{\phi}_k}_i\rVert \leq C_1\epsilon + C_2\lVert \bm{w}^{k,\bar{\phi}_k}_t(\tilde{\functionQ}^k_i)\rVert,
\end{equation}
where recall $\tilde{\functionQ}^k_i$ is the output of the value iteration process in Step~\ref{step:11}.
\end{lemma}
\proof{Proof of Lemma~\ref{lemma:stimulating_algo:convergence_norm}.}
Let $\bar{\bm{r}}^k_i\coloneqq (\bar{r}^k_i(s,a):s\in\bS_i,a\in\bA_i)$, $\mathcal{P}_i^{\bar{\phi}_k}\coloneqq \bigl[p_i^{\bar{\phi}_k}(s,a,s',a')\bigr]_{|\bS_i||\bA_i|\times|\bS_i||\bA_i|}$ and $\hat{\mathcal{P}}^{\bar{\phi}_k}_i=\bigl[\hat{p}^{\bar{\phi}_k}_i(s,a,s',a')\bigr]_{|\bS_i||\bA_i|\times |\bS_i||\bA_i|}$ with
\begin{equation}
 p^{\bar{\phi}_k}_i(s,a,s',a')=\begin{cases}
  p_i(s,a,s') , &\text{if } a' = \actionbarphi_{k,i}(s'), s'\neq s_0,\\
  0,&\text{otherwise},
 \end{cases}   
\end{equation}
and
\begin{equation}
  \hat{p}^{\bar{\phi}_k}_i(s,a,s',a')=\begin{cases}
  \hat{p}_i(s,a,s') , &\text{if } a' = \actionbarphi_{k,i}(s'), s'\neq s_0,\\
  0,&\text{otherwise},
 \end{cases}    
\end{equation}
for all $(s,a,s',a')\in(\bS_i\times\bA_i)^2$,
and $\bm{\Delta}^k_i\coloneqq \functionQ^{k,\bar{\phi}_k}_i - \tilde{\functionQ}^k_{i,\infty}$.

For $k\in[K]$ and $i\in[I]$, from Proposition~\ref{prop:stimulating_algo:convergence_span} and Lemma~\ref{lemma:value_interation:convergence}, in the \partialref{condition:good_case}{good case}, $\tilde{\functionQ}^k_{i,\infty} \coloneqq \lim_{m\rightarrow \infty} \tilde{\functionQ}^k_{i,m}$ exists and satisfies $\tilde{\functionQ}^k_{i,\infty}= \hat{\mathcal{T}}^{k,\bar{\phi}_k}_i\tilde{\functionQ}^k_{i,\infty} = \bar{\bm{r}}^k_i + \mathcal{P}^{\bar{\phi}_k}_i\tilde{\functionQ}^k_{i,\infty} + \bm{w}^{k,\bar{\phi}_k}_t(\tilde{\functionQ}^k_{i,\infty})$, where $\bm{w}^{k,\bar{\phi}_k}_t(\cdot)$ is defined in \eqref{eqn:define_general_w}.
We further obtain
\begin{equation}
(\mathbb{I}-\mathcal{P}_i^{\bar{\phi}_k})(\functionQ^{k,\bar{\phi}_k}_i-\bm{\Delta}^k_i) = \bar{\bm{r}}^k_i+\bm{w}^{k,\bar{\phi}_k}_t(\tilde{\functionQ}^k_{i,\infty}), 
\end{equation}
where $\mathbb{I}$ is the identity matrix.
Since $(\mathbb{I} - \mathcal{P}_i^{\bar{\phi}_k})\functionQ^{k,\bar{\phi}_k}_i = \bar{\bm{r}}^k_i$, it becomes 
\begin{equation}\label{eqn:stimulating_algo:convergence_norm:5}
(\mathcal{P}_i^{\bar{\phi}_k}-\mathbb{I})\bm{\Delta}^k_i = \bm{w}^{k,\bar{\phi}_k}_t(\tilde{\functionQ}^k_{i,\infty}).
\end{equation}
Since, for any $\bar{\phi}_k\in\oPhiLocal$, $s_0$ is reachable, with positive probability, from any other state in $\bS_i$ within a finite time period,
and the state space $\bS_i$ is finite, $\functionQ^{k,\bar{\phi}_k}_i$ is the unique solution to $(\mathbb{I} - \mathcal{P}_i^{\bar{\phi}_k})\functionQ^{k,\bar{\phi}_k}_i = \bar{\bm{r}}^k_i$ ($|\bS_i||\bA_i|$ Bellman equations) for any bounded $\bar{\bm{r}}^k_i\in\bF(\bS_i\times\bA_i)$; that is, $(\mathbb{I} - \mathcal{P}_i^{\bar{\phi}_k})$ is invertible.
Equation \eqref{eqn:stimulating_algo:convergence_norm:5} leads to 
\begin{equation}\label{eqn:stimulating_algo:convergence_norm:6}
    \bm{\Delta}^k_i = (\mathcal{P}^{\bar{\phi}_k}_i-\mathbb{I})^{-1}\bm{w}^{k,\bar{\phi}_k}_t(\tilde{\functionQ}^k_{i,\infty}).
\end{equation}

Together with Lemma~\ref{lemma:convergence:value_iteration}, in the \partialref{condition:good_case}{good case}, there exists $C_1,C_2<\infty$ such that, for any $\epsilon$ used for the stopping condition of the value iteration in Step~\ref{step:11},
\begin{multline}
\lVert \tilde{\functionQ}^k_i - \functionQ^{k,\bar{\phi}_k}_i\rVert \leq \lVert \tilde{\functionQ}^k_i - \tilde{\functionQ}^k_{i,\infty}\rVert +\lVert \tilde{\functionQ}^k_{i,\infty} - \functionQ^{k,\bar{\phi}_k}_i\rVert
= \lVert \tilde{\functionQ}^k_i - \tilde{\functionQ}^k_{i,\infty}\rVert + \lVert\bm{\Delta}^k_i\rVert\\
\leq \lVert \tilde{\functionQ}^k_i - \tilde{\functionQ}^k_{i,\infty}\rVert + C_2 \lVert\bm{w}^{k,\bar{\phi}_k}_t(\tilde{\functionQ}^k_i)\rVert + C_2\bigl\lVert\hat{\mathcal{P}}^{\bar{\phi}_k}_i - \mathcal{P}^{\bar{\phi}_k}_i\bigr\rVert\bigl\lVert\tilde{\functionQ}^k_{i,\infty} - \tilde{\functionQ}^k_i\bigr\rVert
\leq C_1\epsilon + C_2 \lVert\bm{w}^{k,\bar{\phi}_k}_t(\tilde{\functionQ}^k_i)\rVert. 
\end{multline}
It proves the lemma.
\endproof

\proof{Proof of Proposition~\ref{prop:stimulating_algo:convergence_norm}.}
It is a straightforward result of Lemma~\ref{lemma:good_case} and Lemma~\ref{lemma:stimulating_algo:convergence_norm}.
\endproof

\section{Proof of Proposition~\ref{prop:asym_opt:LP}}\label{app:asym_opt:LP}
\proof{Proof of Proposition~\ref{prop:asym_opt:LP}.}
We start with showing that \eqref{eqn:assumption:action} is sufficient for \eqref{eqn:prop:asym_opt:LP:1}.
For any policy $\dbphi\in\dblPhi$ and $t\in[T]_0$, define an action matrix $\mathcal{A}^{\dbphi,h}(t) \coloneqq [\alpha^{\dbphi,h}_{i,s,\iota}(t)]_{(\sum_{i\in[I]}|\bS_i|)\times\mathcal{I}}$,  where
\begin{equation}
    \alpha^{\dbphi,h}_{i,s,\iota}(t) \coloneqq \begin{cases}
        \mathbb{P}\bigl\{\actiondbphi^{h}_{i,n}(t)=\actionphi_{\iota}~\bigl|~s^{\dbphi}_{i,n}(t) = s,\bh(t)\bigr\},&\text{if }i_{\iota}=i,s_{\iota}=s,\\
        0,&\text{otherwise},
    \end{cases}
\end{equation}
where $\mathbb{P}\bigl\{\actiondbphi^{h}_{i,n}(t)=\actionphi_{\iota}~\bigl|~s^{\dbphi}_{i,n}(t) = s,\bh(t)\bigr\}$ is equal for any $n\in[N_i]$.
Note that, in general, $\mathcal{A}^{\dbphi,h}(t)$ is dependent on $\bh(t)$ and is a random variable; and, in the special case with $\dbphi\in\dblPhiz$, $\mathcal{A}^{\dbphi,h}(t)$ becomes deterministic and identical for all $h\in\mathbb{N}_+$.
For $\dbpsi\in\dblPhi$ and $\dbphi\in\dblPhiz$, define $\mathcal{P}^{\dbpsi,h}(t)\coloneqq [p^{\dbpsi,h}_{\iota,\iota'}(t)]_{\mathcal{I}\times\mathcal{I}}\coloneqq \mathcal{P}\mathcal{A}^{\dbpsi,h}(t+1)$ and $\mathcal{P}^{\dbphi}(t)\coloneqq [p^{\dbphi}_{\iota,\iota'}(t)]_{\mathcal{I}\times\mathcal{I}} \coloneqq \mathcal{P}\mathcal{A}^{\dbphi,h}(t+1)$, where the former is random and the latter is deterministic and identical for all $h\in\mathbb{N}_+$.




For $\dbphi\in\dblPhi$ and $\dbpsi\in\dblPhiz$, if we assume that, for any $\epsilon>0$, there exist $C_1<\infty$, $C_2 > 0$ and  $H<\infty$ such that for all $h > H$, 
\begin{equation}\label{eqn:prop:asym_opt:LP:3}
    \lVert\bm{z}^{\dbpsi,h}(t) - \bm{z}^{\dbphi,h}(t)\rVert\leq C_1 e^{-C_2 h} + \epsilon,
\end{equation}
then 
\begin{equation}\label{eqn:prop:asym_opt:LP:4}
    \lVert\bm{\upsilon}^{\dbpsi,h}(t+1) - \bm{\upsilon}^{\dbphi,h}(t+1)\rVert= \lVert(\bm{z}^{\dbpsi,h}(t) - \bm{z}^{\dbphi,h}(t))^T\mathcal{P}\rVert\leq C_1 e^{-C_2 h} + \epsilon.
\end{equation}
For any $\epsilon''>0$, there exist $0<\epsilon''' < \epsilon''$, $0<\epsilon'<\epsilon'''/K$, $0<\epsilon<\epsilon'''/K-\epsilon'$, $C>0$ and $H<\infty$ such that, for all $h>H$,
\begin{multline}\label{eqn:prop:asym_opt:LP:5}
    \mathbb{P}\Bigl\{\bigl\lVert \mathcal{A}^{\dbpsi,h}(t+1) - \mathcal{A}^{\dbphi,h}(t+1)\bigr\rVert > \epsilon''\Bigr\}\\
    \leq  \mathbb{P}\biggl\{\Bigl\lVert\mathcal{A}^{\dbpsi,h}(t+1) - \mathcal{A}^{\dbphi,h}(t+1)\Bigr\rVert>\epsilon''~\biggl|~\lVert \bm{\Upsilon}^{\dbpsi,h}(t+1) - \bm{\upsilon}^{\dbpsi,h}(t+1)\rVert \leq \epsilon'\biggr\} + e^{-C h}\\
    \leq \mathbb{P}\biggl\{K\Bigl\lVert\bm{\Upsilon}^{\dbpsi,h}(t+1) - \bm{\upsilon}^{\dbphi,h}(t+1)\Bigr\rVert>\epsilon'''~\biggl|~\lVert \bm{\Upsilon}^{\dbpsi,h}(t+1) - \bm{\upsilon}^{\dbpsi,h}(t+1)\rVert \leq \epsilon'\biggr\}  + e^{-Ch}\\
    \leq \mathbb{P}\Bigl\{\bigl\lVert \bm{\upsilon}^{\dbpsi,h}(t+1)-\bm{\upsilon}^{\dbphi,h}(t+1)\bigr\rVert > \epsilon'''/K - \epsilon'\Bigr\} + e^{-Ch}
    \leq e^{-Ch},
\end{multline}
where the first inequality is from Theorem~\ref{theorem:convergence_Z_exp}, and the second inequality comes from \eqref{eqn:assumption:action} with $K$ the constant in \eqref{eqn:assumption:action}.
We further obtain that, for any $\epsilon>0$, there exist $H<\infty$ and $C >0$ such that, for all $h>H$,
\begin{multline}\label{eqn:prop:asym_opt:LP:6}
\bigl\lVert \bm{z}^{\dbpsi,h}(t+1) - \bm{z}^{\dbphi,h}(t+1)\bigr\rVert = \Bigl\lVert \mathbb{E}\bigl[\bm{Z}^{\dbpsi,h}(t) \mathcal{P}^{\dbpsi,h}(t)\bigr] - \bm{z}^{\dbphi,h}(t)\mathcal{P}^{\dbphi}(t)\Bigr\rVert\\
\stackrel{(a)}{\leq} \mathbb{E}\Bigl[ \bigl\lVert \bm{Z}^{\dbpsi,h}(t) \mathcal{P}^{\dbpsi,h}(t) - \bm{z}^{\dbphi,h}(t)\mathcal{P}^{\dbphi}(t)\bigr\rVert~\Bigl|~\lVert\mathcal{P}^{\dbpsi,h}(t)-\mathcal{P}^{\dbphi}(t) \rVert\leq \epsilon\Bigr]\mathbb{P}\bigl\{\lVert\mathcal{P}^{\dbpsi,h}(t)-\mathcal{P}^{\dbphi}(t) \rVert\leq \epsilon\bigr\} + e^{-Ch}
\\ 
\leq \mathbb{E}\Bigl[\bigl\lVert \bm{Z}^{\dbpsi,h}(t)\bigl(\mathcal{P}^{\dbpsi,h}(t)-\mathcal{P}^{\dbphi}(t)\bigr)~\bigl|~\bigl\lVert \mathcal{P}^{\dbpsi,h}(t)-\mathcal{P}^{\dbphi}(t)\bigr\rVert \leq \epsilon\bigr\rVert\Bigr] + \Bigl\lVert(\bm{z}^{\dbpsi,h}(t) - \bm{z}^{\dbphi,h}(t))\mathcal{P}^{\dbphi}(t)\Bigr\rVert \\
~~~~~~~~~~~~~~~~~~~~~~~~~~~~~~~~~~~~~~~~~~~~~~~~~~~~~~~~~~~~~~~~~~~~~~~~~~+\mathbb{E}\Bigl[\bigl\lVert (\bm{Z}^{\dbpsi,h}(t) - \bm{z}^{\dbpsi,h}(t))\mathcal{P}^{\dbpsi,h}(t)\bigr\rVert\Bigr]
+ e^{-Ch}
\\\stackrel{(b)}{\leq}
\epsilon\mathbb{E}\Bigl[\bigl\lVert \bm{Z}^{\dbpsi,h}(t)\bigr\lVert\Bigr] + \bigl\lVert \bm{z}^{\dbpsi,h}(t)-\bm{z}^{\dbphi,h}(t)\bigr\rVert +\epsilon + 2e^{-Ch}
\\
\leq 2e^{-Ch} + 3\epsilon + C_1e^{-C_2h},
\end{multline}
where Inequality (a) is based on \eqref{eqn:prop:asym_opt:LP:5}, and Inequality (b) relies on Theorem~\ref{theorem:convergence_Z_exp} and \eqref{eqn:prop:asym_opt:LP:3}.
Together with the initial condition $\bm{Z}^{\dbpsi,h}(0) = \bm{Z}^{\dbphi,h}(0)$, which satisfies \eqref{eqn:prop:asym_opt:LP:3}, we prove the proposition.

\endproof

\section{Proof of Proposition~\ref{prop:converge_LP_opt}}
\label{app:prop:convergence_LP_opt}
\proof{Proof of Proposition~\ref{prop:converge_LP_opt}.}
By time $T^*$, based on the stopping condition in Step~\ref{step:10} of the stimulate process and Corollary~\ref{coro:converge_transition_prob}, for any $\epsilon > 0$, there exist $C>0$ and $H<\infty$ such that, for all $h>H$,
\begin{equation}\label{eqn:prop:convergence_LP_opt:2}
    \mathbb{P}\Bigl\{\bigl\lVert \mathcal{P} -\hat{\mathcal{P}}\bigr\rVert >\epsilon \Bigr\} \leq e^{-Ch}.
\end{equation}

For $\bm{x}_0 = \mathcal{X}\bm{z}^{\dbvtheta,h}(T^*)$, let $\Delta\mathcal{P}(t) \coloneqq \bar{\mathcal{P}}^*\mathcal{A}^{\dbphi^*(\varepsilon,\bm{x}_0),h}(t)-\mathcal{P}\mathcal{A}^{\dbpsi(\varepsilon,\bm{x}_0),h}(t)$ and $\Delta \bm{z}(t) = \bm{x}^*_{t-T^*}(\varepsilon,\bm{x}_0)-\bm{z}^{\dbpsi(\varepsilon,\bm{x}_0),h}(t) $,
where recall $\dbphi^*(\varepsilon,\bm{x}_0)$ is led by $\bm{x}^*(\varepsilon,\bm{x}_0)$, which, together with $\bar{\mathcal{P}}^*$ and $\bar{\bm{r}}^*$, is an optimal solution to \eqref{eqn:obj:linear-programming:optimistic}-\eqref{eqn:constraint:linear-programming:optimistic:5}, and $\dbpsi(\varepsilon,\bm{x}_0)\in\dblPhi$ is applicable to the original problem and satisfies \eqref{eqn:assumption:action} for $\dbpsi=\dbpsi(\varepsilon,\bm{x}_0)$ and $\phi = \dbphi^*(\varepsilon,\bm{x}_0)$.
Note that, here, $\bar{\mathcal{P}}^*$ and $\bar{\bm{r}}^*$ are random variables, for which the randomness comes from $\hat{\mathcal{P}}$ and $\hat{\bm{r}}$ in the $(\varepsilon,\bm{x}_0)$-LP problem in \eqref{eqn:obj:linear-programming:optimistic}-\eqref{eqn:constraint:linear-programming:optimistic:5}.

For $t \geq T^*$, we assume that, for any $\epsilon>0$, 
\begin{equation}\label{eqn:assumption:proof_LP}
    \lVert \bm{z}^{\dbvtheta,h}(t)-\bm{x}^*_{t-T^*}(\varepsilon,\bm{x}_0)\rVert = C_1 e^{-C_2h}+\epsilon.
\end{equation}
Given an instance of $\hat{\mathcal{P}}$ and $\hat{\bm{r}}$,
\begin{multline}\label{eqn:prop:convergence_LP_opt:3}
\bigl\lVert \bm{z}^{\dbvtheta,h}(t+1) - \bm{x}^*_{t+1-T^*}(\varepsilon,\bm{x}_0)\bigr\rVert  \stackrel{(a)}{=} \bigl\lVert \bm{z}^{\dbpsi(\varepsilon,\bm{x}_0),h}(t+1) - \bm{x}^*_{t+1-T^*}(\varepsilon,\bm{x}_0)\bigr\rVert\\
\leq \Bigl\lVert \mathbb{E}\bigl[\bm{Z}^{\dbpsi(\varepsilon,\bm{x}_0),h}(t)\mathcal{P}\mathcal{A}^{\dbpsi(\varepsilon,\bm{x}_0),h}(t)\bigr] - \bm{x}^*_{t-T^*}(\varepsilon,\bm{x}_0)\bar{\mathcal{P}}^*\mathcal{A}^{\dbphi^*(\varepsilon,\bm{x}_0),h}(t)\Bigr\rVert~~~~~~~~~~~~~~~~~~~~~~~~\\
\leq \Bigl\lVert \mathbb{E}\bigl[\bm{Z}^{\dbpsi(\varepsilon,\bm{x}_0),h}(t)\Delta\mathcal{P}(t)\bigr]\Bigr\rVert +\Bigl\lVert \mathbb{E}\bigl[\bm{Z}^{\dbpsi(\varepsilon,\bm{x}_0),h}(t)-\bm{x}^*_{t-T^*}(\varepsilon,\bm{x}_0)\bigr]\bar{\mathcal{P}}^*\mathcal{A}^{\dbphi^*(\varepsilon,\bm{x}_0),h}(t)\Bigr\rVert\\
\leq  \mathbb{E}\Bigl[\bigl\lVert \Delta\mathcal{P}(t)\bigr\rVert\Bigr] + \bigl\lVert \Delta \bm{z}(t)\bigr\rVert \Bigl\lVert \bar{\mathcal{P}}^*\mathcal{A}^{\dbphi^*(\varepsilon,\bm{x}_0),h}(t)\Bigr\rVert,
\end{multline}
where Equality~(a) holds under the condition $\bm{z}^{\dbpsi(\varepsilon,\bm{x}_0),h}(T^*) = \bm{z}^{\dbvtheta,h}(T^*)$.
Similar to \eqref{eqn:prop:asym_opt:LP:5}, for such $t\geq T^*$ and any $\epsilon > 0$, there exist $0<\varepsilon<\epsilon'<\epsilon$,  $H < \infty$ and $C,C'>0$ such that, for all $h>H$,
\begin{multline}\label{eqn:prop:convergence_LP_opt:4}
    \mathbb{P}\Bigl\{\bigl\lVert \Delta \mathcal{P}(t)\bigr\rVert>\epsilon\Bigr\} 
    \leq \mathbb{P}\Bigl\{\bigl\lVert (\mathcal{P}-\bar{\mathcal{P}}^*) \mathcal{A}^{\dbpsi(\varepsilon,\bm{x}_0),h}(t)\bigr\rVert + \bigl\lVert \bar{\mathcal{P}}^*(\mathcal{A}^{\dbpsi(\varepsilon,\bm{x}_0),h}(t) - \mathcal{A}^{\dbphi^*(\varepsilon,\bm{x}_0),h}(t))\bigr\rVert>\epsilon\Bigr\} 
    \\
    \stackrel{(a)}{\leq} \mathbb{P}\Bigl\{\bigl\lVert \mathcal{P}-\bar{\mathcal{P}}^*\bigr\rVert\bigl\lVert \mathcal{A}^{\dbpsi(\varepsilon,\bm{x}_0),h}(t)\bigr\rVert >\epsilon/2\Bigr\} + \mathbb{P}\Bigl\{\bigl \lVert\bm{\upsilon}^{\dbpsi(\varepsilon,\bm{x}_0),h}(t) - \mathcal{U}\bm{x}^*_{t-T^*}(\varepsilon,\bm{x}_0) \bigr\rVert > \frac{\epsilon'}{K\lVert\bar{\mathcal{P}}^*\rVert} \Bigr\} + e^{-C'h}\\
    \stackrel{(b)}{\leq} e^{-Ch} + \mathbb{P}\Bigl\{\bigl \lVert\Delta\bm{z}(t) \bigr\rVert > \frac{\epsilon'}{K\lVert\mathcal{P}\rVert\lVert\mathcal{U}\rVert} \Bigr\}+e^{-C'h}\\
    \stackrel{(c)}{\leq} e^{-Ch} + e^{-C' h}\leq e^{-\min\{C,C'\}h},
\end{multline}
where Inequality (a) is along the same lines as those for \eqref{eqn:prop:asym_opt:LP:5} (which is based on Theorem~\ref{theorem:convergence_Z_exp} and \eqref{eqn:assumption:action}), Inequality (b) holds because \eqref{eqn:prop:convergence_LP_opt:2} and $\lVert \bar{\mathcal{P}}^*-\hat{\mathcal{P}}\rVert \leq \varepsilon$, and Inequality (c) is based on 
\eqref{eqn:assumption:proof_LP}.
Plugging \eqref{eqn:assumption:proof_LP} and \eqref{eqn:prop:convergence_LP_opt:4} in \eqref{eqn:prop:convergence_LP_opt:3}, for any $\epsilon>0$, there exist $C_3,H <\infty$ and $C_4 >0$ such that, for all $h>H$,
\begin{equation}
    \bigl\lVert \bm{z}^{\dbvtheta,h}(t+1)-\bm{x}^*_{t+1-T^*}(\varepsilon,\bm{x}_0)\bigr\rVert \leq C_3 e^{-C_4h} + \epsilon.
\end{equation}
We now go back to the initial condition $\bm{z}^{\dbvtheta,h}(T^*) = \bm{z}^{\dbphi^*(\varepsilon,\bm{x}_0),h}(T^*)$, which satisfies \eqref{eqn:assumption:proof_LP}.
It proves the proposition.

\endproof

\section{Proof of Corollary~\ref{coro:asym_opt:OALP}}
\label{app:coro:asym_opt:OALP}

\proof{Proof of Corollary~\ref{coro:asym_opt:OALP}.}
Let $\mathfrak{d}(h) \coloneqq \lVert \mathcal{P}-\hat{\mathcal{P}}\rVert + \lVert \bm{r} - \hat{\bm{r}}\rVert$.
From Corollary~\ref{coro:converge_transition_prob}, for any $\epsilon > 0$, we obtain
\begin{equation}\label{eqn:app:coro:asym_opt:OALP:0}
\lim_{h\rightarrow \infty} \mathbb{P}\bigl\{\scrd(h)> \epsilon\bigr\} = 0.
\end{equation}
That is, for any $\varepsilon > 0$ and initial condition $\bm{x}_0$, there exists $H<\infty$ such that, for all $h>H$, an optimal solution $\bm{x}^*$ to the linear programming problem \eqref{eqn:obj:linear-programming}-\eqref{eqn:constraint:linear-programming:3}, is feasible to the $(\varepsilon,\bm{x}_0)$-LP problem \eqref{eqn:obj:linear-programming:optimistic}-\eqref{eqn:constraint:linear-programming:optimistic:5}, and, in this case, 
we have
\begin{equation}\label{eqn:app:coro:asym_opt:OALP:1}
    \frac{1}{h}\Gamma^{\dbpsi(\varepsilon,\bm{x}_0)}(T,\bs_0)\leq\frac{1}{h}\Gamma^{\phi^*}(T,\bs_0)\leq \mathbb{E}\biggl[\sum_{t\in[T]_0}\sum_{\iota\in[\mathcal{I}]}N^0_{i_{\iota}}\bar{r}^*_{i_{\iota}}(s_{\iota},\actionphi_{\iota})x^*_{\iota,t}(\varepsilon,\bm{x}_0)\biggr],
\end{equation}
where recall $\bm{x}^*(\varepsilon,\bm{x}_0)$ is an optimal solution to the $(\varepsilon,\bm{x}_0)$-LP problem, 
$\dbphi^*,\dbphi^*(\varepsilon,\bm{x}_0)\in\dblPhiz$ are the policies led by $\bm{x}^*,\bm{x}^*(\varepsilon,\bm{x}_0)$, respectively, and $\bs_0$ is the initial state for the WCG process that can be interpreted to $\bm{z}^{\dbphi^*,h}(0)=\bm{z}^{\dbpsi(\varepsilon,\bm{x}_0),h}(0)=\bm{z}_0$ through the relationship defined in \eqref{eqn:define_Z}, for which $\bm{x}_0=\mathcal{X}\bm{z}_0$.
The first inequality in \eqref{eqn:app:coro:asym_opt:OALP:1} holds because $\dbphi^*$ is optimal to the relaxed problem \eqref{eqn:objective:finite-time:h} and \eqref{eqn:constraint:relax:h}, for which $\Gamma^{\dbphi^*}(T,\bs_0)$ is an upper bound for any policy $\dbphi\in\dblPhi$ applicable to the original problem (described in \eqref{eqn:objective:finite-time:h} and \eqref{eqn:constraint:linear:h}), including $\dbpsi(\varepsilon,\bm{x}_0)$.
From \eqref{eqn:app:coro:asym_opt:OALP:1}, 
\begin{multline}\label{eqn:app:coro:asym_opt:OALP:2}
\Bigl\lvert \Gamma^{\dbphi^*}(T,\bs_0) - \Gamma^{\dbpsi(\varepsilon,\bm{x}_0)}(T,\bs_0)\Bigr\rvert
\leq \mathbb{E}\biggl[\Bigl\lvert \sum_{t\in[T]_0}\sum_{\iota\in[\mathcal{I}]}N^0_{i_{\iota}}\bar{r}^*_{i_{\iota}}(s_{\iota},\actionphi_{\iota})x^*_{\iota,t}(\varepsilon,\bm{x}_0) - \Gamma^{\dbpsi(\varepsilon,\bm{x}_0)}(T,\bs_0)\Bigr\rvert \biggr]\\
\leq \sum_{t\in[T]_0}\sum_{i\in[I]}N^0_i\mathbb{E}\Bigl[\bigl\lvert\bar{\bm{r}}^*\cdot \bm{x}^*_{t-T^*}(\varepsilon,\bm{x}_0) - \bm{r}\cdot \bm{z}^{\dbpsi(\varepsilon,\bm{x}_0)}(t)\bigr\rvert\Bigr]\\
\leq \sum_{t\in[T]_0}\sum_{i\in[I]}N^0_i\mathbb{E}\Bigl[\bigl\lVert \bar{\bm{r}}^* - \bm{r} \bigr\rVert\bigl\lVert \bm{z}^{\dbpsi(\varepsilon,\bm{x}_0)}(t)\bigr\rVert + \bigl\lVert \bm{z}^{\dbpsi(\varepsilon,\bm{x}_0)}(t)-\bm{x}^*_{t-T^*}(\varepsilon,\bm{x}_0)\bigr\rVert \bigl\lVert \bar{\bm{r}}^*\bigr\rVert \Bigr].
\end{multline}
Together with Proposition~\ref{prop:converge_LP_opt},  for any $\epsilon > 0$, there exist $0<\varepsilon<\epsilon$, $C_1,C_3,H<\infty$ and $C_2,C_4>0$ such that
\begin{multline}\label{eqn:app:coro:asym_opt:OALP:3}
  \Bigl\lvert \Gamma^{\dbphi^*}(T,\bs_0) - \Gamma^{\dbpsi(\varepsilon,\bm{x}_0)}(T,\bs_0)\Bigr\rvert\\
\leq   \sum_{t\in[T]_0}\sum_{i\in[I]}N^0_i\mathbb{E}\Bigl[\bigl\lVert \bar{\bm{r}}^* - \bm{r} \bigr\rVert + \bigl\lVert \bm{z}^{\dbpsi(\varepsilon,\bm{x}_0)}(t)-\bm{x}^*_{t-T^*}(\varepsilon,\bm{x}_0)\bigr\rVert \max_{\iota\in[\mathcal{I}]}(|\hat{r}_{i_{\iota}}(s_{\iota},\actionphi_{\iota})|+\varepsilon) \Bigr]\\
\stackrel{(a)}{\leq}\sum_{t\in[T]_0}\sum_{i\in[I]}N^0_i\mathbb{E}\Bigl[\scrd(h) + \varepsilon +(C_1 e^{-C_2h}+\epsilon)\bigl(\max_{\iota\in[\mathcal{I}]}|\hat{r}_{i_{\iota}}(s_{\iota},\actionphi_{\iota})|+\epsilon \bigr)\Bigr]\\
\leq C_3 e^{-C_4 h} + O(\epsilon),
\end{multline}
where Inequality (a) is based on $\scrd(h) + \varepsilon \geq \lVert \bar{\bm{r}}^*-\bm{r}\rVert $ and $\varepsilon < \epsilon$.
It proves \eqref{eqn:coro:asym_opt:OALP}.
 

\endproof

\bibliographystyle{alpha}
\bibliography{references/IEEEabrv,references/bib1}
\end{document}

%% file: proof_propn1.tex
\section{Proof of Proposition~\ref{prop:decomposition}}\label{app:prop:decomposition}
First we recall the relaxed  subproblem optimizer:
\begin{equation}\label{eqn:sub-problem1}
    D_i(\bm{\gamma})\coloneqq \max_{\dbphi\in\dbsPhiz}\limsup_{T\rightarrow \infty}\frac{1}{T} \sum_{t\in[T]}\mathbb{E}\biggl[ R^{\dbphi}_{i,n}(t)- \sum_{\ell\in[L]}\gamma_{\ell} f_{i,\ell}\bigl((s^{\dbphi}_{i,n}(t),\dba^{\dbphi}_{i,n}(t)\bigr)\biggr].
\end{equation}
By \eqref{eq:4ddd}, this is equal to
\begin{equation}
  \label{eq:1app}
    D_i(\bm{\gamma})= \max_{\dbphi\in\dbsPhiz}\sum_{\begin{subarray}~s\in \bS_i\\a\in \bA_{i}\end{subarray}} \bigl(R^{\dbphi}_{i,n}(t)- \sum_{\ell\in[L]}\gamma_{\ell} f_{i,\ell}\bigl((s^{\dbphi}_{i,n}(t),\dba^{\dbphi}_{i,n}(t)\bigr)\pi^{{\dbphi}}_{i}(s,a).
  \end{equation}
  Here $\pi^{{\dbphi}}_{i}(s,a)$ is the steady state distribution for the randomized policy $\dbphi$ and so satisfies
  \begin{equation}
        \label{eq:3app}
    \begin{aligned}
      \sum_{\substack{s\in \bS_{i}\\a\in \bA_{i}}}p(s,a,s')\pi_{i}^{\dbphi}(s,a)&= \sum_{\substack{s\in \bS_{i}\\a\in \bA_{i}}}\pi_{i}^{\dbphi}(s,a)\delta_{s'}(s)\\
      \sum_{\substack{s\in \bS_{i}\\a\in \bA_{i}}}\pi^{\dbphi}(s,a)&=1 \qquad  \\
    \end{aligned}
      \end{equation}
where $\delta_{s'}(s)= 1$ if $s=s'$ and zero otherwise.
We write $\Xi_{i}$ for the set of all $\pi$ satisfying \eqref{eq:3app}.

For any $i\in[I]$, given $\pi\in \Xi_{i}$, we can obtain a randomized policy with that $\pi$ as its steady state distribution, by just normalizing:
\begin{equation}\label{eqn:prop:decomposition:4}
\alpha^{\dbphi}_{i}(a,s) = \frac{\pi_i(s,a)}{\sum_{a'\in\bA_i}\pi_i(s,a')}.
\end{equation}
It follows that
\begin{equation}
  \label{eq:1appaa}
    D_i(\bm{\gamma})= \max_{\pi\in \Xi_{i}}\sum_{\begin{subarray}~s\in \bS_i\\a\in \bA_{i}\end{subarray}} \bigl(r_{i}(s,a)- \sum_{\ell\in[L]}\gamma_{\ell} f_{i,\ell}(s,a)\bigr)\pi(s,a).
  \end{equation}
  We note here, by the same ergodic argument, that for $D(\gamma)$ in \eqref{eqn:long-run:dual-func},
\begin{equation}
  \label{eq:1appaa}
    D(\bm{\gamma})= \max_{\pi\in \Xi}\sum_{i\in [I]}\sum_{n\in [N_{i}]}\sum_{\begin{subarray}~s\in \bS_i\\a\in \bA_{i}\end{subarray}} \bigl(r_{i}(s,a)- \sum_{\ell\in[L]}\gamma_{\ell} f_{i,\ell}(s,a)\bigr)\pi_{i}(s,a),
  \end{equation}
  where $\Xi=\prod_{i\in [I]}\Xi_{i}$, and $\pi=(\pi_{i})_{i}$. Since there is no constraint that couples the indices $i$, and the processes in gang $i$ are identical,  it is clear that
  \begin{equation}
    \label{eq:4}
    D(\gamma)=\sum_{i\in[I]}N_{i}D_{i}(\gamma).
  \end{equation}
This proves part 2. of Proposition $1$.

It remains to show that $D_{i}(\gamma)$ is achieved by a deterministic policy; that is, that we can always find a $\pi_{i}^{\phi}$ arising from $\phi\in \oPhiLocal$.  To do this, we recall the maximizing  $Q$-factor in \eqref{eqn:bellman} that follows from Bellman's equation; in fact, we use the randomized form of it:
\begin{multline}\label{eqn:decomposition:1}
    D_i(\bm{\gamma}) + Q^{\dbphi,\bm{\gamma}}_i(s,a) \\
    =
    r_i(s,a)-\sum_{\ell\in[L]}\gamma_{\ell}f_{i,\ell}(s,a)+ \sum_{s'\in\bS_i\backslash\{s_0\}}p_i(s,a,s') \max_{\bm{\alpha}:\bS_{i}\to \Delta_{\bA_{i}}}\sum_{a'\in\bA_i}\alpha(s',a')Q^{\dbphi,\bm{\gamma}}_i(s',a'),
  \end{multline}
  where we regard $\alpha:\bS_{i}\to \Delta_{\bA_{i}}$ as a function on $\bS_{i}\times \bA_{i}$.
  Now, by convexity of $\bA_{i}$,  it is clear that
  \[\max_{\bm{\alpha}:\bS_{i}\to \Delta_{\bA_{i}}}\sum_{a'\in\bA_i}\alpha(s',a')
      Q^{\dbphi,\bm{\gamma}}_i(s',a')
    \]
    is maximized by the extreme functions
     $\alpha:\bS_{i}\to \ext(\Delta_{\bA_{i}})$,
    which correspond to deterministic actions. The result follows. \endproof

%% file: 20241201.bbl
\newcommand{\etalchar}[1]{$^{#1}$}
\begin{thebibliography}{WRMS20}

\bibitem[AB16]{avrachenkov2016whittle}
Konstantin~E Avrachenkov and Vivek~S Borkar.
\newblock {Whittle} index policy for crawling ephemeral content.
\newblock {\em IEEE Transactions on Control of Network Systems}, 5(1):446--455,
  2016.

\bibitem[AB22]{avrachenkov2022whittle}
Konstantin~E Avrachenkov and Vivek~S Borkar.
\newblock Whittle index based q-learning for restless bandits with average
  reward.
\newblock {\em Automatica}, 139:110186, 2022.

\bibitem[AM08]{adelman2008relaxations}
Daniel Adelman and Adam~J Mersereau.
\newblock Relaxations of weakly coupled stochastic dynamic programs.
\newblock {\em Oper. Res.}, 56(3):712--727, 2008.

\bibitem[AM23]{akbarzadeh2023learning}
Nima Akbarzadeh and Aditya Mahajan.
\newblock On learning whittle index policy for restless bandits with scalable
  regret.
\newblock {\em IEEE Transactions on Control of Network Systems}, 2023.

\bibitem[BS20]{brown2020index}
David~B Brown and James~E Smith.
\newblock Index policies and performance bounds for dynamic selection problems.
\newblock {\em Management Science}, 66(7):3029--3050, 2020.

\bibitem[BT96]{bertsekas1996neuro}
D.~P. Bertsekas and John~N. Tsitsiklis.
\newblock {\em Neuro-dynamic programming}.
\newblock Athena Scientific, Belmont, MA, 1996.

\bibitem[BT15]{bertsekas2015parallel}
Dimitri Bertsekas and John Tsitsiklis.
\newblock {\em Parallel and distributed computation: numerical methods}.
\newblock Athena Scientific, 2015.

\bibitem[BZ23]{brown2023fluid}
David~B Brown and Jingwei Zhang.
\newblock Fluid policies, reoptimization, and performance guarantees in dynamic
  resource allocation.
\newblock {\em Operations Research}, 2023.

\bibitem[CL55]{coddington1955theory}
Earl~A Coddington and Norman Levinson.
\newblock {\em Theory of ordinary differential equations}.
\newblock Tata McGraw-Hill Education, 1955.

\bibitem[FM20]{fu2020energy}
Jing Fu and Bill Moran.
\newblock Energy-efficient job-assignment policy with asymptotically guaranteed
  performance deviation.
\newblock {\em IEEE/ACM Transactions on Networking}, 28(3):1325--1338, 2020.

\bibitem[FMG{\etalchar{+}}16]{fu2016asymptotic}
Jing Fu, Bill Moran, Jun Guo, Eric W~M Wong, and Moshe Zukerman.
\newblock Asymptotically optimal job assignment for energy-efficient
  processor-sharing server farms.
\newblock {\em IEEE Journal on Selected Areas in Communications},
  34(12):4008--4023, 2016.

\bibitem[FMT22]{fu2018restless}
Jing Fu, Bill Moran, and Peter~G. Taylor.
\newblock A restless bandit model for resource allocation, competition and
  reservation.
\newblock {\em Operations Research}, 70(1):416--431, 2022.

\bibitem[FMTX]{fu2020resource}
Jing Fu, Bill Moran, Peter~G. Taylor, and Chenchen Xing.
\newblock Resource competition in virtual network embedding with released
  physical resources.
\newblock {\em Stochastic Models}, pages 231 -- 263.
\newblock Dec. 2020.

\bibitem[FNMT19]{fu2019towards}
Jing Fu, Yoni Nazarathy, Sarat Moka, and Peter~G Taylor.
\newblock Towards {Q-learning} the whittle index for restless bandits.
\newblock In {\em 2019 Australian \& New Zealand Control Conference (ANZCC)},
  pages 249--254, Auckland, 2019. IEEE.

\bibitem[FW12]{freidlin2012random}
Mark~I. Freidlin and Alexander~D. Wentzell.
\newblock {\em Random perturbations of dynamical systems}.
\newblock Springer Science \& Business Media, 2012.
\newblock translated by J. Sz{\"u}cs.

\bibitem[FWC24]{fu2024patrolling}
Jing Fu, Zengfu Wang, and Jie Chen.
\newblock Coordinated multi-agent patrolling with state-dependent cost rates --
  asymptotically optimal policies for large-scale systems.
\newblock {\em arXiv preprint arXiv:2309.13388}, 2024.

\bibitem[FZL25]{fu2025restless}
Jing Fu, Lele Zhang, and Zhiyuan Liu.
\newblock A restless bandit model for dynamic ride matching with reneging
  travelers.
\newblock {\em European Journal of Operational Research}, 320(3):581--592,
  2025.

\bibitem[GGW11]{gittins2011multiarmed}
J~C. Gittins, K~Glazebrook, and R~R. Weber.
\newblock {\em Multi-armed bandit allocation indices: 2nd edition}.
\newblock Wiley, Mar. 2011.

\bibitem[GGY23]{gast2023linear}
Nicolas Gast, Bruno Gaujal, and Chen Yan.
\newblock Linear program-based policies for restless bandits: Necessary and
  sufficient conditions for (exponentially fast) asymptotic optimality.
\newblock {\em Mathematics of Operations Research}, 2023.

\bibitem[GGY24]{gast2024reoptimization}
Nicolas Gast, Bruno Gaujal, and Chen Yan.
\newblock Reoptimization nearly solves weakly coupled markov decision
  processes.
\newblock {\em HAL: hal-04570177 [PREPRINT]}, 2024.

\bibitem[JJL{\etalchar{+}}23]{jiang2023online}
Bowen Jiang, Bo~Jiang, Jian Li, Tao Lin, Xinbing Wang, and Chenghu Zhou.
\newblock Online restless bandits with unobserved states.
\newblock In {\em International Conference on Machine Learning}, pages
  15041--15066. PMLR, 2023.

\bibitem[KD07]{krishnamurthy2007structured}
Vikram Krishnamurthy and Dejan~V Djonin.
\newblock Structured threshold policies for dynamic sensor scheduling—a
  partially observed {Markov}decision process approach.
\newblock {\em {IEEE} Trans. Signal Process.}, 55(10):4938--4957, 2007.

\bibitem[NM01a]{nino2001restless}
Jos{\'e} Ni{\~n}o-Mora.
\newblock Restless bandits, partial conservation laws and indexability.
\newblock {\em Advances in Applied Probability}, 33(1):76--98, 2001.

\bibitem[NM01b]{ninomora2001restless}
Jos{\'e} Ni{\~n}o-Mora.
\newblock Restless bandits, partial conservation laws and indexability.
\newblock {\em Advances in Applied Probability}, pages 76--98, 2001.

\bibitem[NM02]{nino2002dynamic}
Jos{\'e} Ni{\~n}o-Mora.
\newblock Dynamic allocation indices for restless projects and queueing
  admission control: a polyhedral approach.
\newblock {\em Mathematical programming}, 93(3):361--413, 2002.

\bibitem[NM06]{nino2006restless}
Jos{\'e} Ni{\~n}o-Mora.
\newblock Restless bandit marginal productivity indices, diminishing returns,
  and optimal control of make-to-order/make-to-stock m/g/1 queues.
\newblock {\em Mathematics of Operations Research}, 31(1):50--84, 2006.

\bibitem[NM07]{nino2007dynamic}
Jos{\'e} Ni{\~n}o-Mora.
\newblock Dynamic priority allocation via restless bandit marginal productivity
  indices.
\newblock {\em TOP}, 15(2):161--198, Sep. 2007.

\bibitem[NM20]{ninomora2020verification}
Jos{\'e} Ni{\~n}o-Mora.
\newblock A verification theorem for threshold-indexability of real-state
  discounted restless bandits.
\newblock {\em Mathematics of Operations Research}, 45(2):465--496, 2020.

\bibitem[NM22]{nino2022multigear}
Jos{\'e} Ni{\~n}o-Mora.
\newblock Multi-gear bandits, partial conservation laws, and indexability.
\newblock {\em Mathematics}, 10(14):2497, 2022.

\bibitem[NM23]{ninomora2023survey}
Jos{\'e} Ni{\~n}o-Mora.
\newblock Markovian restless bandits and index policies: {A} review.
\newblock {\em Mathematics}, 11(7):1639, 2023.

\bibitem[Nn08]{nino2008index}
J.~Ni\~{n}o{-}Mora.
\newblock An index policy for multiarmed multimode restless bandits.
\newblock In J.~Baras and C.~Courcoubetis, editors, {\em 3rd International
  Conference on Performance Evaluation Methodologies and Tools (ValueTools),
  Athens, Greece}, ACM International Conference Proceedings Series, Brussels,
  Belgium, 2008. ICST.

\bibitem[Nor98]{norris1998markov}
James~R Norris.
\newblock {\em Markov chains}.
\newblock Number~2. Cambridge university press, 1998.

\bibitem[PT99]{papadimitriou1999complexity}
Christos~H. Papadimitriou and John~N. Tsitsiklis.
\newblock The complexity of optimal queuing network control.
\newblock {\em Mathematics of Operations Research}, 24(2):293--305, May 1999.

\bibitem[Put05]{puterman2005markov}
Martin~L Puterman.
\newblock {\em Markov decision processes: discrete stochastic dynamic
  programming}.
\newblock John Wiley \& Sons, 2005.

\bibitem[Ros92]{ross1992applied}
Sheldon~M. Ross.
\newblock {\em Applied probability models with optimization applications}.
\newblock Dover Publications (New York), 1992.

\bibitem[Ver16]{verloop2016asymptotically}
Ina~Maria Verloop.
\newblock Asymptotically optimal priority policies for indexable and
  non-indexable restless bandits.
\newblock {\em Ann. Appl. Probab.}, 26(4):1947--1995, Aug. 2016.

\bibitem[Vil16]{villar2016indexability}
Sof{\'\i}a~S Villar.
\newblock Indexability and optimal index policies for a class of reinitialising
  restless bandits.
\newblock {\em Probability in the engineering and informational sciences},
  30(1):1--23, 2016.

\bibitem[Web07]{weber2007comments}
Richard Weber.
\newblock Comments on: Dynamic priority allocation via restless bandit marginal
  productivity indices.
\newblock {\em Transactions on Operations Research}, 15(2):211--216, 2007.

\bibitem[Whi88]{whittle1988restless}
Peter Whittle.
\newblock Restless bandits: Activity allocation in a changing world.
\newblock {\em Journal of Applied Probability}, 25:287--298, 1988.

\bibitem[WRMS20]{wang2020whittle}
Jiazheng Wang, Xiaoqiang Ren, Yilin Mo, and Ling Shi.
\newblock {Whittle} index policy for dynamic multi-channel allocation in remote
  state estimation.
\newblock {\em {IEEE} Trans. Autom. Control}, 65(2), Feb. 2020.

\bibitem[WW90]{weber1990index}
Richard~R. Weber and Gideon Weiss.
\newblock On an index policy for restless bandits.
\newblock {\em Journal of Applied Probability}, (3):637--648, Sep. 1990.

\bibitem[WXTT23]{wang2023optimistic}
Kai Wang, Lily Xu, Aparna Taneja, and Milind Tambe.
\newblock Optimistic whittle index policy: Online learning for restless
  bandits.
\newblock In {\em Proceedings of the AAAI Conference on Artificial
  Intelligence}, volume~37, pages 10131--10139, 2023.

\end{thebibliography}
